\documentclass[10pt,a4paper]{article}

\usepackage{a4wide}
\usepackage[english]{babel} 
\usepackage{amsmath,amssymb}
\usepackage[amsmath,thmmarks]{ntheorem}
\usepackage{booktabs, longtable, pdflscape, multirow}
\usepackage{afterpage}
\usepackage[font=small]{caption}
\usepackage{pstricks, pst-plot, picins, wrapfig}
\usepackage{calc}
\usepackage{subcaption}
\usepackage[T1]{fontenc}
\usepackage{microtype}

\theoremstyle{plain}
\theoremheaderfont{\normalfont\bfseries}
\theorembodyfont{\itshape}
\theoremseparator{.}
\theoremnumbering{arabic}
\theoremsymbol{}
\newtheorem{thm}{Theorem}[section]
\newtheorem{cor}[thm]{Corollary}
\newtheorem{lem}[thm]{Lemma}
\newtheorem{prop}[thm]{Proposition}

\theoremstyle{plain}
\theoremheaderfont{\normalfont\bfseries}
\theorembodyfont{\normalfont}
\theoremseparator{.}
\theoremnumbering{arabic}
\theoremsymbol{}
\newtheorem{defn}[thm]{Definition}
\newtheorem{rmk}[thm]{Remark}

\theoremstyle{plain}
\theoremheaderfont{\normalfont\bfseries}
\theorembodyfont{\normalfont}
\theoremseparator{.}
\theoremnumbering{arabic}
\theoremsymbol{\ensuremath{\lozenge}}

\theoremstyle{nonumberplain}
\theoremheaderfont{\itshape}
\theorembodyfont{\normalfont}
\theoremseparator{.}
\theoremnumbering{arabic}
\theoremsymbol{\ensuremath{\square}}
\newtheorem{proof}{Proof}

\theoremstyle{empty}
\theoremheaderfont{\itshape}
\theorembodyfont{\normalfont}
\theoremseparator{.}
\theoremnumbering{arabic}
\theoremsymbol{\ensuremath{\square}}
\newtheorem{proofoptname}{}

\newcommand{\N}{\mathbb{N}}
\newcommand{\Z}{\mathbb{Z}}
\newcommand{\Q}{\mathbb{Q}}
\newcommand{\R}{\mathbb{R}}
\newcommand{\C}{\mathbb{C}}

\newcommand{\gauss}{\mbox{$_2 F_1$}}

\newcommand{\nfn}{\mbox{$_n F_{n-1}$}}
\newcommand{\A}{\mathcal{A}}
\newcommand{\HA}{H_\A(\bbeta)}
\newcommand{\HAi}[1]{H_{\A_{#1}}(\bbeta)}

\newcommand{\sign}{\sigma_\A(\bbeta)}
\newcommand{\signk}[1]{\sigma_\A(#1 \bbeta)}
\newcommand{\trian}{\mathcal{T}}
\newcommand{\lat}{\mathbb{L}}

\newcommand{\fp}[1]{\{ #1 \}}
\newcommand{\entier}[1]{\lfloor #1 \rfloor}
\newcommand{\inv}[1]{\frac{1}{#1}}
\newcommand{\half}{\inv{2}}

\newcommand{\la}{\lambda}
\newcommand{\subs}{\subseteq}

\newcommand{\poch}[2]{(#1)_{#2}}
\newcommand{\vol}{\tr{Vol}}

\newcommand{\ld}{\ldots}
\newcommand{\tr}{\textrm}

\newcommand{\vect}[1]{\boldsymbol{#1}}
\newcommand{\va}{\vect{a}}

\newcommand{\ve}{\vect{e}}
\newcommand{\vl}{\vect{l}}

\newcommand{\vp}{\vect{p}}
\newcommand{\vq}{\vect{q}}

\newcommand{\vx}{\vect{x}}

\newcommand{\bbeta}{\vect{\beta}}

\newcommand{\blam}{\vect{\la}}
\newcommand{\bmu}{\vect{\mu}}

\newcommand{\ahypergeometric}{$\A$-hyper\-geo\-me\-tric }

\title{A classification of the irreducible algebraic \ahypergeometric functions associated to planar point configurations}
\author{Esther Bod\thanks{esther.bod@gmail.com. Department of Mathematics, Utrecht University, P.O.\ Box 80010, 3508 TA Utrecht, The Netherlands.
This work was supported by the Netherlands Organisation for Scientific Research (NWO) by grant number OND1331860.}}
\date{\today}

\begin{document}

\maketitle

\begin{abstract}
We consider \ahypergeometric functions associated to normal sets in the plane.
We give a classification of all point configurations for which there exists a parameter vector such that the associated hypergeometric function is algebraic.
In particular, we show that there are no irreducible algebraic functions if the number of boundary points is sufficiently large and $\A$ is not a pyramid.
\end{abstract}


\section{Introduction}\label{sec:polygons-introduction}

\ahypergeometric functions were introduced in the 1980's by Gelfand, Graev, Kapranov and Zelevinsky in a series of papers~\cite{ggz_holonomic_systems_equations_series_hypergeometric_type,gzk_equations_hypergeometric_newton_polyhedra,gzk_hypergeometric_functions_and_toral_manifolds,gkz_hypergeometric_functions_and_toral_manifolds_correction}.
They are generalizations of, among others, the hypergeometric functions called after Gauss, Appell and Horn.
It is a classical question for which values of the parameters these functions are algebraic.
In 1873, Schwarz gave a list of all $(a,b,c)$ such that the Gauss function $\gauss(a,b,c|z)$ is irreducible and algebraic.
This list is reproduced in Table~\ref{tab:gauss_abc_orbits}, which contains all such tuples up to permutations of $a$ and $b$, translations over $\Z^3$ and multiplications by integers coprime with the smallest common denominator of $a$, $b$ and $c$.
This list has been extended to the general hypergeometric function $\nfn$~\cite{beukers_heckman_monodromy_nfn-1} and all Appell-Lauricella and Horn functions in~\cite{sasaki_lauricella_fd,beazley_cohen_wolfart_algebraic_appell_lauricella,kato_appell_f4,kato_appell_f2,bod_algebraic_appell_horn}.

\afterpage{
\begin{table} 
\centering
\caption{The tuples $(a,b,c)$ such that $\gauss(a,b,c|z)$ is irreducible and algebraic} 
\label{tab:gauss_abc_orbits}
\renewcommand{\arraystretch}{1.5}
\begin{tabular}{l@{\hspace{0.55cm}}l@{\hspace{0.55cm}}l@{\hspace{0.55cm}}l@{\hspace{0.55cm}}l@{\hspace{0.55cm}}l@{\hspace{0.55cm}}l} 
\toprule

\multicolumn{7}{l}{$(r, -r, \half)$, 
$(r, r+\half, \half)$ \tr{and} 
$(r, r+\half, 2r)$
with $2r \not\in \Z$} \\

$(\half, 	\inv{6}, 	\inv{3})$ &
$(\inv{6}, 	\frac{5}{6}, 	\inv{5})$ &
$(\inv{10}, 	\frac{9}{10}, 	\inv{3})$ &
$(\inv{15}, 	\frac{7}{15}, 	\inv{3})$ &
$(\inv{20}, 	\frac{13}{20}, 	\half)$ &
$(\inv{24}, 	\frac{19}{24}, 	\half)$ &
$(\inv{60}, 	\frac{41}{60}, 	\half)$ \\

$(\inv{4}, 	\frac{3}{4}, 	\inv{3})$ & 
$(\inv{6}, 	\frac{5}{12}, 	\inv{3})$ &
$(\inv{10}, 	\frac{9}{10}, 	\inv{5})$ &
$(\inv{15}, 	\frac{7}{15}, 	\inv{5})$ &
$(\inv{20}, 	\frac{13}{20}, 	\inv{5})$ &
$(\inv{24}, 	\frac{19}{24}, 	\inv{3})$ &
$(\inv{60}, 	\frac{41}{60}, 	\inv{5})$ \\

$(\inv{4}, 	\frac{7}{12}, 	\half)$ & 
$(\inv{6}, 	\frac{5}{12}, 	\inv{4})$ &
$(\inv{10}, 	\frac{13}{30}, 	\inv{3})$ &
$(\inv{15}, 	\frac{11}{15}, 	\inv{5})$ &
$(\inv{24}, 	\frac{13}{24}, 	\inv{3})$ &
$(\inv{30}, 	\frac{11}{30}, 	\inv{5})$ &
$(\inv{60}, 	\frac{49}{60}, 	\half)$ \\

$(\inv{4}, 	\frac{7}{12}, 	\inv{3})$ &
$(\inv{6}, 	\frac{11}{30}, 	\inv{3})$ &
$(\inv{10}, 	\frac{13}{30}, 	\inv{5})$ &
$(\inv{15}, 	\frac{11}{15}, 	\frac{3}{5})$ &
$(\inv{24}, 	\frac{13}{24}, 	\inv{4})$ &
$(\inv{30}, 	\frac{19}{30}, 	\inv{3})$ &
$(\inv{60}, 	\frac{49}{60}, 	\inv{3})$ \\

$(\inv{6}, 	\frac{5}{6}, 	\inv{3})$ &
$(\inv{6}, 	\frac{11}{30}, 	\inv{5})$ &
$(\inv{12}, 	\frac{5}{12}, 	\inv{4})$ &
$(\inv{20}, 	\frac{11}{20}, 	\inv{5})$ &
$(\inv{24}, 	\frac{17}{24}, 	\half)$ &
$(\inv{60}, 	\frac{31}{60}, 	\inv{3})$ \\

$(\inv{6}, 	\frac{5}{6}, 	\inv{4})$ &
$(\inv{10}, 	\frac{3}{10}, 	\inv{5})$ &
$(\inv{12}, 	\frac{7}{12}, 	\inv{3})$ &
$(\inv{20}, 	\frac{11}{20}, 	\frac{2}{5})$ &
$(\inv{24}, 	\frac{17}{24}, 	\inv{4})$ &
$(\inv{60}, 	\frac{31}{60}, 	\inv{5})$ \\
\bottomrule
\end{tabular}
\renewcommand{\arraystretch}{1}
\end{table}}

In this paper, we extend this list to all \ahypergeometric functions where $\A$ is a subset of $\Z^2$ or $\Z^3$.
We start by recalling the definition and some basic facts about hypergeometric functions. 

\begin{defn}\label{defn:Ahypergeometric_function}
Let $\A=\{\va_1, \ld, \va_N\}$ be a finite subset of $\Z^r$ such that $\Z \A = \Z^r$ and there exists a linear form $h$ on $\R^r$ such that $h(\va_i)=1$ for all $i$.
The \emph{lattice of relations} of $\A$ is $\lat = \{(l_1, \ld, l_N) \in \Z^N \ | \ l_1 \va_1 + \ld + l_N \va_N = 0 \}$.
Let $\bbeta \in \C^r$ and denote by $\partial_i$ the differential operator $\frac{\partial}{\partial z_i}$.
The \emph{\ahypergeometric system associated to $\A$ and $\bbeta$}, denoted $\HA$, consists of two sets of differential equations:
\begin{itemize}
\item
the \emph{structure equations}: for all $\vl = (l_1, \ld, l_N) \in \lat$ 
\begin{equation*}
\square_{\vl} \Phi = 
\left( \prod_{l_i>0} \partial_i^{l_i} \right) \Phi - \left( \prod_{l_i<0} \partial_i^{-l_i} \right) \Phi = 0. 
\end{equation*}
\item 
the \emph{homogeneity} or \emph{Euler equations}: for $1 \leq i \leq r$ 
\begin{equation*}
a_{1i} z_1 \partial_1 \Phi + \ld + a_{Ni} z_N \partial_N \Phi = \beta_i \Phi.
\end{equation*}
\end{itemize}
The solutions $\Phi(z_1, \ld, z_N)$ of this system are called \emph{\ahypergeometric functions}.
\end{defn}

We call $\A$ a \emph{pyramid} if $\A = \A' \cup \{\va_i\}$ where all elements of $\A'$ lie in an $r-2$ dimensional hyperplane.
If this is the case, one can easily show that an \ahypergeometric function is the product of an $\A'$-hypergeometric function and a monomial factor $z_i^{\gamma_i}$.
This monomial factor does not influence algebraicity, so we will restrict ourselves to non-pyramidal sets $\A$.

The convex hull and real non-negative cone spanned by the elements of $\A$ will be denoted $Q(\A)$ and $C(\A)$, respectively. 
Throughout this paper, we make the assumption that $\A$ is \emph{normal}, i.e., $C(\A) \cap \Z^r = \Z_{\geq 0} \A$.
In this paper we only consider planar point configurations, i.e., subsets of $\Z^2$ and $\Z^3$.
For these $\A$, normality is equivalent to the condition $Q(\A) \cap \Z^r = \A$.

We will also assume that $\A$ lies in the hyperplane $x_r = 1$.
This is possible by the following proposition:  

\begin{prop}\label{prop:Aheight1}
Let $\A$ be as in Definition~\ref{defn:Ahypergeometric_function}.
Then there is an isomorphism $f$ of $\Z^r$ such that $f(\A) \subs \Z^{r-1} \times \{1\}$. 
\end{prop}

\begin{proof}
Write $\A = \{\va_1, \ld, \va_N\}$ and $h(\vx) = h_1 x_1 + \ld + h_r x_r$.
For each $i$ there exist $\la_1, \ld, \la_N \in \Z$ such that $\ve_i = \la_1 \va_1 + \ld + \la_N \va_N$.
It follows from $h(\va_j)=1$ that $h_i = h(\ve_i) = \la_1 + \ld + \la_N \in \Z$.
Furthermore, $\gcd(h_1, \ld, h_r) = 1$ because $h(\va_1)=1$.
There exists a basis of $\Z^r$ containing the vector $(h_1, \ld, h_r)$.
The matrix whose rows are the elements of this basis, with $r^{\tr{th}}$ row $(h_1, \ld, h_r)$, gives the desired isomorphism.
\end{proof}

We will only consider irreducible systems $\HA$, i.e.\ systems for which the monodromy group acts irreducibly on the solution space.
In order to check irreducibility without computing the monodromy group, we use the almost equivalent concept of non-resonance.

\begin{defn}\label{defn:resonance}
$\HA$ is called \emph{resonant} if $\bbeta+\Z^r$ contains a point in a face of $C(\A)$.
Otherwise $\HA$ is \emph{non-resonant}.
$\HA$ is called \emph{totally non-resonant} if $\bbeta+\Z^r$ contains no point in any hyperplane spanned by $r-1$ independent elements of $\A$.
\end{defn}

Note that total non-resonance implies non-resonance, and hence irreducibility. \\

\begin{thm}[{{\cite[Theorem~2.11]{gkz_euler_integrals_and_hypergeometric_functions}},{\cite[Theorem~4.1]{schulze_walther_resonance_reducibility}}}]\label{thm:irrreducibility}
Suppose that $\A$ is not a pyramid.
Then $\HA$ is non-resonant if and only if it is irreducible.
\end{thm}

To sets $\A \subs \Z^3$ we associate the polygon $P(\A) = \{ \vx \in \R^2 \ | \ (\vx,1) \in Q(\A) \}$.
This is a convex lattice polygon in the sense of the following definition:

\begin{defn}\label{defn:polygons}
A \emph{simple polygon} $P$ is a subset of $\R^2$ bounded by a closed chain of line segments, that does not intersect itself.
If the points in which the line segments meet are integral, then $P$ is called a \emph{lattice polygon}.
Two lattice polygons $P$ and $P'$ are \emph{isomorphic} if $P'$ is the image of $P$ under a linear isomorphism of $\Z^2$.
\end{defn}

Conversely, to each convex lattice polygon $P$ there is an associated set $\A$: by defining $\A = (P \cap \Z^2) \times \{1\}$, we have $P=P(\A)$.
Then $\A$ is clearly normal and it spans $\Z^3$ over $\Z$ if $P$ contains 3 integral points forming a triangle with area 1. 
It will follow from Corollary~\ref{cor:existence_triangulation_polygon} that such triangles exist, so $\A$ satisfies the conditions of Definition~\ref{defn:Ahypergeometric_function}.
Note that isomorphisms of polygons correspond to isomorphisms of the associated point configurations, and hence preserve algebraicity of the \ahypergeometric functions. \\

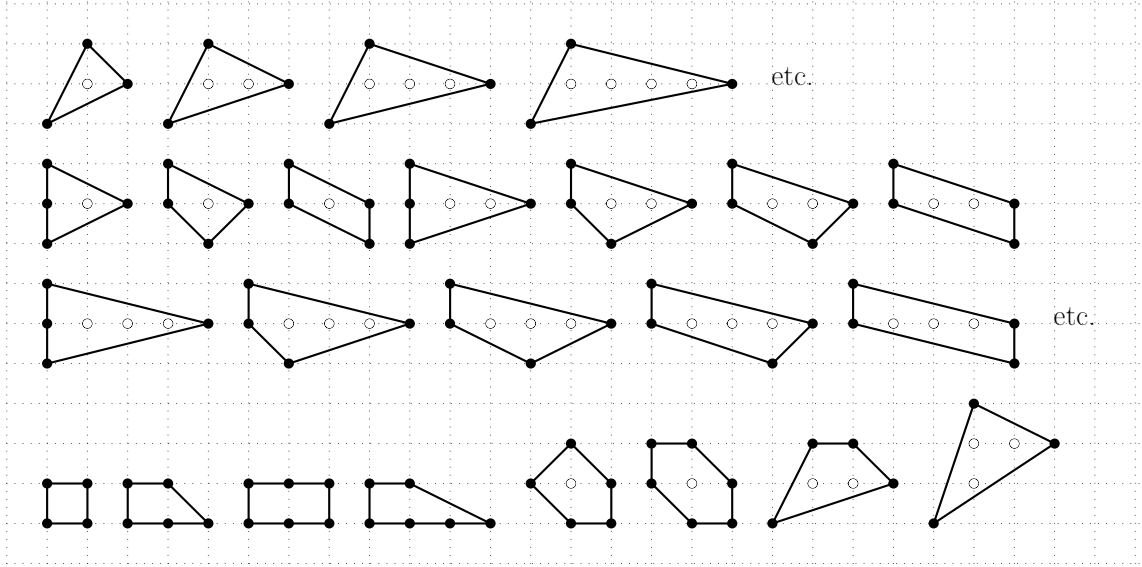
\begin{figure}
\centering
\psscalebox{0.53}{
\begin{pspicture}(28,14)
\psgrid[subgriddiv=1,gridlabels=0,gridwidth=0.2pt,griddots=5](28,14)

\psdots[dotscale=2](1,11)(3,12)(2,13)
\psdots[dotscale=2,dotstyle=o](2,12)
\psline[linewidth=1.5pt](1,11)(3,12)(2,13)(1,11)

\psdots[dotscale=2](4,11)(7,12)(5,13)
\psdots[dotscale=2,dotstyle=o](5,12)(6,12)
\psline[linewidth=1.5pt](4,11)(7,12)(5,13)(4,11)

\psdots[dotscale=2](8,11)(12,12)(9,13)
\psdots[dotscale=2,dotstyle=o](9,12)(10,12)(11,12)
\psline[linewidth=1.5pt](8,11)(12,12)(9,13)(8,11)

\psdots[dotscale=2](13,11)(18,12)(14,13)
\psdots[dotscale=2,dotstyle=o](14,12)(15,12)(16,12)(17,12)
\psline[linewidth=1.5pt](13,11)(18,12)(14,13)(13,11)

\rput(19.5,12.2){\huge{etc.}}

\psdots[dotscale=2](1,8)(1,9)(3,9)(1,10)
\psdots[dotscale=2,dotstyle=o](2,9)
\psline[linewidth=1.5pt](1,8)(3,9)(1,10)(1,8)

\psdots[dotscale=2](5,8)(4,9)(6,9)(4,10)
\psdots[dotscale=2,dotstyle=o](5,9)
\psline[linewidth=1.5pt](5,8)(6,9)(4,10)(4,9)(5,8)

\psdots[dotscale=2](9,8)(7,9)(9,9)(7,10)
\psdots[dotscale=2,dotstyle=o](8,9)
\psline[linewidth=1.5pt](9,8)(9,9)(7,10)(7,9)(9,8)

\psdots[dotscale=2](10,8)(10,9)(13,9)(10,10)
\psdots[dotscale=2,dotstyle=o](11,9)(12,9)
\psline[linewidth=1.5pt](10,8)(13,9)(10,10)(10,8)

\psdots[dotscale=2](15,8)(14,9)(17,9)(14,10)
\psdots[dotscale=2,dotstyle=o](15,9)(16,9)
\psline[linewidth=1.5pt](15,8)(17,9)(14,10)(14,9)(15,8)

\psdots[dotscale=2](20,8)(18,9)(21,9)(18,10)
\psdots[dotscale=2,dotstyle=o](19,9)(20,9)
\psline[linewidth=1.5pt](20,8)(21,9)(18,10)(18,9)(20,8)

\psdots[dotscale=2](25,8)(22,9)(25,9)(22,10)
\psdots[dotscale=2,dotstyle=o](23,9)(24,9)
\psline[linewidth=1.5pt](25,8)(25,9)(22,10)(22,9)(25,8)

\psdots[dotscale=2](1,5)(1,6)(5,6)(1,7)
\psdots[dotscale=2,dotstyle=o](2,6)(3,6)(4,6)
\psline[linewidth=1.5pt](1,5)(5,6)(1,7)(1,5)

\psdots[dotscale=2](7,5)(6,6)(10,6)(6,7)
\psdots[dotscale=2,dotstyle=o](7,6)(8,6)(9,6)
\psline[linewidth=1.5pt](7,5)(10,6)(6,7)(6,6)(7,5)

\psdots[dotscale=2](13,5)(11,6)(15,6)(11,7)
\psdots[dotscale=2,dotstyle=o](12,6)(13,6)(14,6)
\psline[linewidth=1.5pt](13,5)(15,6)(11,7)(11,6)(13,5)

\psdots[dotscale=2](19,5)(16,6)(20,6)(16,7)
\psdots[dotscale=2,dotstyle=o](17,6)(18,6)(19,6)
\psline[linewidth=1.5pt](19,5)(20,6)(16,7)(16,6)(19,5)

\psdots[dotscale=2](25,5)(21,6)(25,6)(21,7)
\psdots[dotscale=2,dotstyle=o](22,6)(23,6)(24,6)
\psline[linewidth=1.5pt](25,5)(25,6)(21,7)(21,6)(25,5)

\rput(26.5,6.2){\huge{etc.}}

\psdots[dotscale=2](1,1)(2,1)(1,2)(2,2)
\psline[linewidth=1.5pt](1,1)(1,2)(2,2)(2,1)(1,1)

\psdots[dotscale=2](3,1)(4,1)(5,1)(3,2)(4,2)
\psline[linewidth=1.5pt](3,1)(5,1)(4,2)(3,2)(3,1)

\psdots[dotscale=2](6,1)(7,1)(8,1)(6,2)(7,2)(8,2)
\psline[linewidth=1.5pt](6,1)(8,1)(8,2)(6,2)(6,1)

\psdots[dotscale=2](9,1)(10,1)(11,1)(12,1)(9,2)(10,2)
\psline[linewidth=1.5pt](9,1)(12,1)(10,2)(9,2)(9,1)

\psdots[dotscale=2](14,1)(15,1)(13,2)(15,2)(14,3)
\psdots[dotscale=2,dotstyle=o](14,2)
\psline[linewidth=1.5pt](14,1)(15,1)(15,2)(14,3)(13,2)(14,1)

\psdots[dotscale=2](17,1)(18,1)(16,2)(18,2)(16,3)(17,3)
\psdots[dotscale=2,dotstyle=o](17,2)
\psline[linewidth=1.5pt](17,1)(18,1)(18,2)(17,3)(16,3)(16,2)(17,1)

\psdots[dotscale=2](19,1)(22,2)(20,3)(21,3)
\psdots[dotscale=2,dotstyle=o](20,2)(21,2)
\psline[linewidth=1.5pt](19,1)(22,2)(21,3)(20,3)(19,1)

\psdots[dotscale=2](23,1)(26,3)(24,4)
\psdots[dotscale=2,dotstyle=o](24,2)(24,3)(25,3)
\psline[linewidth=1.5pt](23,1)(26,3)(24,4)(23,1)
\end{pspicture}}
\caption{The polygons $P(\A)$ such that there exists $\bbeta$ such that $\HA$ has irreducible algebraic solutions \label{fig:classification_planar_algebraic_functions}}
\end{figure}

The main result of this paper is the following:

\begin{thm}\label{thm:classification_planar_algebraic_functions}
If $\A$ is a finite normal subset of $\Z^2$, lying on an affine line, 
then there exists $\bbeta \in \Q^2$ such that $H_\A(\bbeta)$ is non-resonant and has algebraic solutions.

If $\A$ is a finite normal non-pyramidal subset of $\Z^3$, lying in an affine plane, 
then there exists $\bbeta \in \Q^3$ such that $H_\A(\bbeta)$ is non-resonant and has algebraic solutions if and only if $P(\A)$ is isomorphic to one of the polygons in Figure~\ref{fig:classification_planar_algebraic_functions}.
\end{thm} 

In the next section we describe how to determine the algebraic functions using a combinatorial criterion and triangulations of polygons.
Then we analyze one-dimensional sets in Section~\ref{sec:polygons-collinear_configurations}.
For the two-dimensional sets, we will mostly use the polygons $P(\A)$ instead of the sets $\A$.
In Section~\ref{sec:polygons-012_interior_points} we will compute all polygons with at most 2 interior lattice points and determine the irreducible algebraic functions.
Then we show in Section~\ref{sec:polygons-many_boundary_points} that all polygons with sufficiently many interior and boundary lattice points contain smaller subpolygons, to which we can reduce the problem of finding irreducible algebraic functions.
Finally, the polygons with many interior points and few boundary points will be treated in Section~\ref{sec:polygons-34_boundary_points}. \\


\section{Algebraic \ahypergeometric functions}\label{sec:polygons-algebraic_hypergeometric_functions}

To determine all algebraic \ahypergeometric functions, we use a combinatorial proven by Beukers.
To state this, we need the concept of apex points.

\begin{defn}\label{defn:apexpoint}
An \emph{apex point} is a point $\vect{p} \in (\bbeta + \Z^r) \cap C(\A)$ such that $\vect{p}-\va \not\in C(\A)$ for all $\va \in \A$.
The number of apex points is called the \emph{signature} of $\A$ and $\bbeta$ and is denoted $\sign$.
\end{defn}

We will use the following combinatorial criterion to determine all algebraic functions:

\begin{thm}[{{\cite[Theorem~1.10]{beukers_algebraic_ahypergeometric_functions}}}]\label{thm:algebraic_solutions}
Suppose that $\HA$ is non-resonant. 
If $\HA$ has algebraic solutions, then $\bbeta \in \Q^r$.
Let $D$ be the smallest common denominator of the coordinates of $\bbeta \in \Q^r$. 
Then the solutions of $\HA$ are algebraic over $\C(\vect{z})$ if and only if $\signk{k} = \vol(Q(\A))$ for all $k \in \Z$ with $1 \leq k < D$ and $\gcd(k,D)=1$.
\end{thm}

Here and everywhere else in this paper, by `volume' we mean the \emph{simplex volume}.
It is a normalized version of the Euclidean volume, such that the simplex spanned by the standard basis has volume 1.
For $\A \subs \Z^2$, the convex hull of $\A$ is a segment of the line $x_2=1$ and the volume is the length of this line segment.
For $\A \subs \Z^3$, the convex hull is a polygon in the plane $x_3=1$ and the volume is twice the area of this polygon. 

Note that $\sign$, and hence algebraicity, only depends on the fractional part $\fp{\bbeta}$ of $\bbeta$ (where $\fp{\bbeta}_i = \fp{\beta_i} = \beta_i - \entier{\beta_i}$).
Throughout this paper, we will assume that $\bbeta \in \Q^r$ and $0 \leq \beta_i < 1$ for all $i$.
A consequence of Theorem~\ref{thm:algebraic_solutions} is the fact that either the solution set of $H_\A(k \bbeta)$ consists of algebraic functions for all $k$ coprime to the smallest common denominator of the coordinates of $\bbeta$, or the solutions are transcendental for all $k$. 
We will call the parameter vectors $\bbeta$ and $k \bbeta$ \emph{conjugated}. \\

It is well-known that the Gauss function $\gauss(a,b,c|z)$ is irreducible if and only if $a$, $b$, $c-a$ and $c-b$ are non-integral.
If it is irreducible, then it is algebraic if and only if for every $k$ coprime with the denominators of $a, b$ and $c$, either $\fp{ka}  \leq \fp{kc} < \fp{kb}$ or $\fp{kb} \leq \fp{kc} < \fp{ka}$ for all $k$ coprime with the smallest common denominator of $a$, $b$ and $c$.
This condition can be found in~\cite{katz_alg_sol_diff_eq} and is called the \emph{interlacing condition}.
We can view $\gauss(a,b,c|z)$ as an \ahypergeometric function with $\A = \{\ve_1, \ve_2, \ve_3, \ve_1+\ve_2-\ve_3\}$ and $\bbeta=(-a,-b,c-1)$ (see~\cite[Section~2.6]{stienstra_gkz_hypergeometric_structures}).
It can easily be computed that the above condition for irreducibility is equivalent to the condition that $\HA$ is non-resonant.
We can use Theorem~\ref{thm:algebraic_solutions} to find similar interlacing conditions for other hypergeometric functions. 
We will use the fact that the number of apex points is constant on certain parts of $[0,1)^r$, when we let $\bbeta$ vary:

\begin{lem}\label{lem:signature_entier_linear_forms}
Let $\HA$ be non-resonant and 
$C(\A) = \{ \vx \in \R^r \ | \ m_1(\vx) \geq 0, \ld, m_n(\vx) \geq 0 \}$
where $m_1, \ld, m_n$ are linear forms with integral coefficients with greatest common divisor 1.
Then $\sign$ depends on $(\entier{m_1(\bbeta)}, \ld, \entier{m_n(\bbeta)})$, but not on $\bbeta$ itself.
\end{lem}

\begin{proof}
Let $\vx \in \Z^r$.
Then $\vx + \bbeta$ is an apex point if and only if $\vx+\bbeta \in C(\A)$ and for all $i$, $\vx-\va_i+\bbeta \not\in C(\A)$.
Equivalently, $m_j(\vx) \geq -m_j(\bbeta$) for all $j$, and for all $i$ there exists $j$ such that $m_j(\vx) < m_j(\va_i)-m_j(\bbeta)$.
Since $m_j(\vx)$ and $m_j(\va_i)$ are integral, whereas $m_j(\bbeta)$ is non-integral (because $\HA$ is non-resonant), the apex points are those $\vx+\bbeta$ such that $m_j(\vx) \geq -\entier{m_j(\bbeta)}$ for all $j$, and for all $i$ there exists $j$ such that $m_j(\vx) \leq m_j(\va_i)-\entier{m_j(\bbeta)}-1$.
Hence the conditions on $\vx+\bbeta$ to be an apex point only depend on $\entier{m_j(\bbeta)}$.
\end{proof}

We now describe how to find a condition on $\bbeta$ to have maximal signature, similar to the interlacing condition for $\gauss$.
As input we need the linear forms $m_i$ that determine the faces of the cone $C(\A)$.
Write $m_i(\vx) = \sum_j m_{ij} x_j$.
Since we only consider $\bbeta$ such that $\beta_i \in [0,1)$ for all $i$, $\entier{m_i(\bbeta)}$ can only take integral values between $\sum_j \min(m_{ij},0)$ and $\sum_j \max(m_{ij},0)$ (both boundaries are excluded, unless they are zero).
Hence $(\entier{m_1(\bbeta)}, \ld, \entier{m_n(\bbeta)})$ takes only finitely many values.
For each of those, it suffices to find one corresponding $\bbeta$ and compute the number of apex points.
Finding such $\bbeta$ boils down to solving a linear system of inequalities.
This can easily be done by hand or using a computer algebra system, which will also detect the values of $(\entier{m_1(\bbeta)}, \ld, \entier{m_n(\bbeta)})$ for which no $\bbeta$ exists.
Having found $\bbeta$, finding apex points can again be done by solving a system of linear inequalities, in this case over the integers.

Using this algorithm, finding the interlacing condition can entirely be done by computer.
However, this algorithm is very slow, although it can be used for small sets $\A$. \\

We will use Theorem~\ref{thm:algebraic_solutions} to find the algebraic functions for small sets $\A$.
Then we will reduce large sets $\A$ to smaller subsets.
For this reduction, $\A$ and its subset need to have compatible triangulations.

\begin{defn}\label{defn:triangulation}
Let $\A$ be a finite subset of $\Z^r$ in an affine hyperplane.
A \emph{triangulation} of $Q(\A)$ is a finite set $\trian = \{Q(V_1), \ld, Q(V_l)\}$ such that each $V_i$ is a subset of $\A$ consisting of $r$ linearly independent elements, $Q(V_i) \cap Q(V_j) = Q(V_i \cap V_j)$ for all $i$ and $j$ and $Q(\A) = \cup_{i=1}^l Q(V_i)$.
If all $Q(V_i)$ have simplex volume 1, then the triangulation is called \emph{unimodular}.
We will call $\{V_1, \ld, V_l\}$ a (unimodular) triangulation of $\A$.
\end{defn}

If $\A \subs \Z^2$ and $Q(\A)$ is a segment of the line $x_2=1$, a unimodular triangulation is a division of this segment into intervals of length 1.
For $\A \subs \Z^3$, triangulations of $\A$ correspond to triangulations of the polygon $P(\A)$.
Such triangulations are dissections of $P(\A)$ into triangles, of which the vertices are lattice points.
A triangulation of $\A$ is unimodular if all triangles in $P(\A)$ have simplex area 1. 
To check that a triangulation is unimodular, one can use Pick's formula:

\begin{lem}[Pick's formula]\label{lem:pick}
Let $P$ be a convex lattice polygon with $i$ interior lattice points and $b$ lattice points on the boundary.
Then its simplex area is equal to $2i+b-2$.
\end{lem}

It follows that a triangle with integral vertices has area 1 if and only if the only integral points are the vertices.

\begin{prop}\label{prop:subsets_algebraic}
Let $\A' \subs \A \subs \Z^r$ be as in Definition~\ref{defn:Ahypergeometric_function}.
Suppose that $\A'$ has a unimodular triangulation that can be extended to a unimodular triangulation of $\A$.
If $\HA$ is non-resonant and has algebraic solutions, then $H_{\A'}(\bbeta)$ is also non-resonant with algebraic solutions.
In particular, if there are no $\bbeta$ such that $H_{\A'}(\bbeta)$ is non-resonant and has algebraic solutions, then are are also no irreducible algebraic functions associated to $\A$. 
\end{prop}

To prove this proposition, we need the following two lemmas:

\begin{lem}\label{lem:unimodular_trian_1apexpoint_simplex}
Let $\{V_1, \ld, V_l\}$ be a unimodular triangulation of $\A$.
Then each cone $C(V_i)$ contains at most one apex point.
\end{lem}

\begin{proof}
Let $\vp$ and $\vq$ be apex points in $C(V_i)$.
By reordering the vectors in $\A$ if necessary, we can assume that $V_i = \{\va_1, \ld, \va_r\}$.
There exist $\la_j, \mu_j \geq 0$ such that $\vp = \la_1 \va_1 + \ld + \la_r \va_r$ and $\vq = \mu_1 \va_1 + \ld + \mu_r \va_r$.
Since $\vp$ is an apex point, we have $\vp-\va \not\in C(\A)$ for all $\va \in \A$.
In particular, $\vp-\va_j \not\in C(V_i)$, so $\la_j < 1$ for $j=1, \ld, r$.
Similarly, $\mu_j<1$.

Note that $\vp-\vq = (\va_1, \ld, \va_r) \cdot (\blam-\bmu)$, so $\blam-\bmu = (\va_1, \ld, \va_r)^{-1} (\vp-\vq)$, where $(\va_1, \ld, \va_r)$ is viewed as an invertible $r \times r$-matrix.
As $\vp$ and $\vq$ are both apex points, we have $\fp{\vp}=\fp{\vq}=\fp{\bbeta}$, so $\vp-\vq \in \Z^r$.
The matrix $(\va_1, \ld, \va_r)$ has determinant $\pm 1$ because the triangulation is unimodular.
It follows that $\blam-\bmu \in \Z^r$.
This implies that $\blam=\bmu$, and hence $\vp=\vq$.
\end{proof}

\begin{lem}\label{lem:unimodular_trian_resonant_transcendental}
Suppose that $\HA$ is non-resonant but not totally non-resonant.
Then there exist $r-1$ independent elements of $\A$ such that $\bbeta + \Z^r$ contains a point of the hyperplane $\mathcal{F}$ through these elements of $\A$.
Suppose that $\A$ has a unimodular triangulation $\{V_1, \ld, V_l\}$ such that $\mathcal{F}$ is a face of one of the cones $C(V_i)$.
Then $\HA$ has transcendental solutions.
\end{lem}

\begin{proof}
As the algebraicity of solutions only depends on $\fp{\bbeta}$, we can assume that $\bbeta$ itself lies on $\mathcal{F}$.
Since $\HA$ is non-resonant, $\mathcal{F}$ is not a face of $C(\A)$ and hence there is a $V_j \neq V_i$ such that $\mathcal{F}$ is also a face of $C(V_j)$.
The sets $V_i$ and $V_j$ have $r-1$ points in common, say $V_i = \{\va_1, \ld, \va_r\}$ and $V_j = \{\va_1, \ld, \va_{r-1}, \va_{r+1}\}$. 
There exists a vector $\blam = (\la_1, \ld, \la_{r-1}, 0)$ such that $\bbeta = V_i \blam = V_j \blam$, where we identify $V_i$ and $V_j$ with the matrices whose columns are the vectors $\va_i$.
By again translating $\bbeta$ if necessary, we can assume that $0 \leq \la_i < 1$ for all $i$.
Suppose that $\vx+\bbeta \in C(V_i)$ is an apex point, with $\vx \in \Z$.
Then we can write $\vx+\bbeta = V_i \bmu$ with $0 \leq \mu_i < 1$.
But then $\bmu - \blam = V_i^{-1} \vx \in \Z^r$ with $-1 < \mu_i - \la_i < 1$, so $\bmu = \blam$ and $\vx = \vect{0}$.
Hence if there is an apex point in $C(V_i)$, then it must be $\bbeta$.
Similarly, the only possible apex point in $C(V_j)$ is $\bbeta$.
By Lemma~\ref{lem:unimodular_trian_1apexpoint_simplex}, the other $C(V_k)$ contain at most $l-2$ apexpoints, so $\sign \leq l-1$.
Theorem~\ref{thm:algebraic_solutions} now implies that the solutions of $\HA$ are transcendental.
\end{proof}

\begin{rmk}\label{rmk:resonant_transcendental}
One can show that systems that are not totally non-resonant have solutions involving logarithms.
Hence the above lemma also holds if $\A$ has no unimodular triangulation, but the proof is more involved.
\end{rmk}

\begin{proofoptname}[Proof of Proposition~\ref{prop:subsets_algebraic}.]
Suppose that $H_{\A'}(\bbeta)$ is resonant.
Then $\HA$ is not totally non-resonant.
Extending the triangulation of $\A'$ to $\A$ gives a triangulation satisfying the hypothesis of Lemma~\ref{lem:unimodular_trian_resonant_transcendental}.
Hence $\HA$ has transcendental solutions.
It follows that $H_{\A'}(\bbeta)$ is non-resonant.

Let $k \in \Z$ be coprime with the smallest common denominator of the coordinates of $\bbeta$.
By Theorem~\ref{thm:algebraic_solutions}, $\signk{k} = \vol(Q(\A))$.
Lemma~\ref{lem:unimodular_trian_1apexpoint_simplex} shows that each cone in the triangulation of $\A$ contains exactly one apex point.
In particular, this holds for the cones in the triangulation of $\A'$.
Hence $\sigma_{\A'}(k \bbeta) = \vol(Q(\A'))$ and $H_{\A'}(\bbeta)$ has algebraic solutions.
\end{proofoptname}

To apply Proposition~\ref{prop:subsets_algebraic}, we now show that each normal subset $\A'$ of $\A$ has a unimodular triangulation that can be extended to a unimodular triangulation of $\A$.
If $\A \subs \Z^2$, this is trivial: 
$Q(\A)$ is a line segment and $Q(\A')$ is a subsegment.
Dividing both segments into parts of length 1 gives the desired triangulations.
For $\A \subs \Z^3$ the proof is a bit more complicated and requires two lemmas.

\begin{lem}\label{lem:existence_triangulation_triangle}
Every lattice triangle has a unimodular triangulation.
\end{lem}

\begin{proof}
Let $P$ be a lattice triangle.
We show that $P$ can be divided into smaller triangles.
Since the area of such triangles is positive and integral, after finitely many steps we will find a unimodular triangulation of $P$.
If $P$ doesn't have interior lattice points or lattice points on the boundary except for the vertices, then the area is 1 and $P$ is triangulated already.
If $P$ has an interior lattice point, then connecting this interior point with the three vertices divides the triangle into three smaller triangles.
If there is a lattice point on an edge of $P$, which is not a vertex, then connecting this with the opposite vertex will divide $P$ into two smaller triangles.
\end{proof}

\begin{lem}\label{lem:extension_triangulation_polygon}
Let $P, P'$ be convex lattice polygons with $P' \subs P$.
Then every unimodular triangulation of $P'$ can be extended to a unimodular triangulation of $P$.
\end{lem}

\begin{figure}
\centering
\centering
\psscalebox{0.8}{
\begin{pspicture}(-1.5,-0.5)(5,4.5)
\psdots(0,0)(2,0)(3,1)(3,2)(1,4)(-1,1.5)(4,1.5)
\psline(0,0)(2,0)(3,1)(3,2)(1,4)(-1,1.5)(0,0)
\psline[linestyle=dashed](1,4)(4,1.5)(3,2)
\psline[linestyle=dashed](3,1)(4,1.5)(2,0)
\rput(4.4,1.5){\large{$v_0$}}
\rput(2.5,0){\large{$v_{-1}$}}
\rput(1.4,4){\large{$v_1$}}
\rput(2.5,1){\large{$w_1$}}
\rput(2.5,2){\large{$w_2$}}
\end{pspicture}}
\caption{Adding a vertex to a triangulated polygon
\label{fig:extended_triangulation}}
\end{figure}
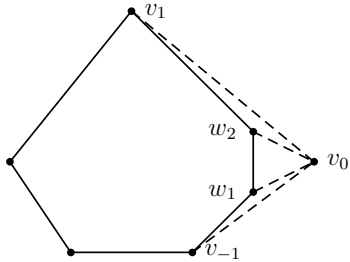

\begin{proof}
Let $V$ be the set of lattice points in $P$ that are not in $P'$.
It is clear that it suffices to show that the lemma holds for $|V|=1$, since we can then add the points in $V$ one by one to $P'$ while preserving the triangulation.
So suppose that $V=\{v_0\}$.
Then $v_0$ is a vertex of $P$, since otherwise it would be contained in $P'$, which is the convex hull of all other lattice points in $P$.
Let $v_{-1}$ and $v_1$ be the previous and next lattice point on the boundary of $P$ (in counterclockwise order; see Figure~\ref{fig:extended_triangulation}).
There can be lattice points on the boundary of $P'$ in between $v_{-1}$ and $v_1$; call them $w_1, \ld, w_k$.
By connecting $v_0$ to $v_{-1}, v_1, w_1, \ld, w_k$, we clearly get a triangulation of $P$ that extends the triangulation of $P'$.
Since the only lattice points in each of the triangles are the vertices, all triangles have area 1 and the triangulation is unimodular.
\end{proof}

\begin{cor}\label{cor:existence_triangulation_polygon}
Let $P$ be a convex lattice polygon and $P'$ a convex subpolygon.
Then $P'$ has a unimodular triangulation that can be extended to a unimodular triangulation of $P$.
\end{cor}

\begin{proof}
By Lemma~\ref{lem:extension_triangulation_polygon}, it suffices to show that $P'$ has a unimodular triangulation.
We prove this by induction on the number of vertices of $P'$.
If $P'$ is a triangle, then the statement follows from Lemma~\ref{lem:existence_triangulation_triangle}.
Suppose that $P'$ has at least 4 vertices.
Then there is a diagonal dividing $P'$ in two smaller polygons.
By induction, we can find a unimodular triangulation of one of these.
Lemma~\ref{lem:extension_triangulation_polygon} shows that we can extend this to a unimodular triangulation of $P'$.
\end{proof}

\begin{cor}\label{cor:extension_less_algebraic_functions_resonant}
Let $P' \subs P$ be convex lattice polygons, with associated sets $\A' \subs \A \subs \Z^3$.
If $\HA$ has irreducible algebraic solutions, then $H_{\A'}(\bbeta)$ is irreducible and its solutions are algebraic.
If there are no $\bbeta$ such that $H_{\A'}(\bbeta)$ has irreducible algebraic solutions, then such $\bbeta$ also don't exist for $\HA$. 
\end{cor}

\begin{proof}
This follows from Corollary~\ref{cor:existence_triangulation_polygon} and Proposition~\ref{prop:subsets_algebraic}.
\end{proof}


\section{Collinear point configurations}\label{sec:polygons-collinear_configurations}

\begin{figure}
\centering
\begin{subfigure}[t]{.6\textwidth}
\centering
\begin{pspicture}(-2,-0.7)(5,2)
\pspolygon[linestyle=none,fillstyle=solid,fillcolor=lightgray](-2,2)(0,0)(5,1.25)(5,2)
\pspolygon[linestyle=none,fillstyle=solid,fillcolor=darkgray](-1,1)(-0.2,1.2)(0,1)(0.8,1.2)(1,1)(1.8,1.2)(2,1)(2.8,1.2)(3,1)(3.8,1.2)(4,1)(0,0)
\psgrid[subgriddiv=1,griddots=1,gridwidth=2pt,gridlabels=0](0,0)(-2,-0.5)(5,2)
\psaxes[ticks=none,linewidth=0.5pt,labels=none](0,0)(-2,-0.5)(5,2)
\psdots[dotsize=5pt](-1,1)(0,1)(1,1)(2,1)(3,1)(4,1)
\psline[linewidth=2pt](-1,1)(4,1)
\psline(-2,2)(0,0)(5,1.25)
\psline(0,0)(0,2)
\psline(0,0)(2,2)
\psline(0,0)(4,2)
\psline(0,0)(5,1.6666667)
\psline[linestyle=dotted](-1,1)(-0.2,1.2)(0,1)(0.8,1.2)(1,1)(1.8,1.2)(2,1)(2.8,1.2)(3,1)(3.8,1.2)(4,1)
\end{pspicture}
\caption{The convex hull (thick line), positive cone (light gray) and apex points (dark gray)
\label{fig:collinear_point_configurations_apexpoints}}
\end{subfigure}
\hfill
\begin{subfigure}[t]{.3\textwidth}
\centering
\psscalebox{0.9}{
\begin{pspicture}(-1,-0.5)(3.2,4)
\pspolygon[fillstyle=solid,fillcolor=lightgray,linestyle=none](0,0)(2.4,0.6)(3,0)
\pspolygon[fillstyle=solid,fillcolor=lightgray,linestyle=none](0,3)(0.6,2.4)(3,3)
\psaxes[Dx=1,Dy=1,dx=3,dy=3](3,3)
\psline[linestyle=dotted,linewidth=1.5pt](3,0)(3,3)
\psline[linestyle=dotted,linewidth=1.5pt](0,3)(3,3)
\psline(0,3)(0.6,2.4)(3,3)
\psline[linestyle=dashed](0,0)(2.4,0.6)(3,0)
\psline[linestyle=dotted,linewidth=1pt](0,2.25)(0.6,2.4)(2.4,0.6)(3,0.75)
\psline(-0.15,0.75)(0.15,0.75)
\psline(-0.15,2.25)(0.15,2.25)
\rput(-0.7,0.75){$\inv{N-2}$}
\rput(-0.7,2.25){$\frac{N-3}{N-2}$}
\end{pspicture}}
\caption{The interlacing condition
\label{fig:interlacing_condition_collinear_points}}
\end{subfigure}
\caption{A collinear point configuration ($N=6$)
\label{fig:collinear_point_configuration}}
\end{figure}
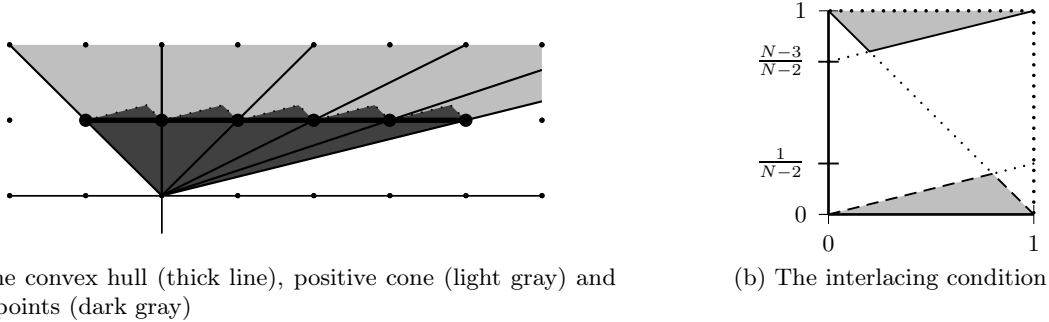

In this section we assume that $\A$ is a finite normal subset of the line $x_2=1$ in $\Z^2$.
We can shift $\A$ so that $\A=\{(k,1) | -1 \leq k \leq N-2\}$.
Since $\A$ spans $\Z^2$, we have $N \geq 2$.
For $N=4$, $\HA$ is the \ahypergeometric system correspronding to the Horn $G_3$ function.
We compute the irreducible algebraic functions for each $N$.
Note that $\HA$ is irreducible if and only if $\beta_1+\beta_2, -\beta_1+(N-2)\beta_2 \not\in \Z$, as can be seen in Figure~\ref{fig:collinear_point_configurations_apexpoints}.
As $Q(\A)$ is an interval of length $N-1$, we determine under which conditions there are $N-1$ apex points.

\begin{lem}\label{lem:collinear_interlacing}
For $0 \leq i \leq N-2$, let $V_i = \{(i-1,1), (i,1)\}$ and let $(x,y) \in C(V_i)$.
Then $(x,y)$ is an apex point if and only if $x+y<i+1$ and $x>(N-2)y-N+i+1$.
There are $N-1$ apex points if and only if $(\entier{\beta_1+\beta_2}, \entier{-\beta_1+(N-2)\beta_2}) \in \{(-1,0), (N-3,1)\}$.
\end{lem}

\begin{proof}
Suppose that $(x,y) \in C(V_i)$ is an apex point.
Since $(x,y) \in C(V_i)$, we have $(i-1)y \leq x \leq iy$.
Then $(x,y)-\va_j \not\in C(\A)$ for all $\va_j$, so in particular $(x-i,y-1), (x-i+1,y-1) \not\in C(\A)$, i.e., $x+y<i+1$ or $x-i>(N-2)(y-1)$, and $x+y<i$ or $x>(N-2)y-N+i+1$.
We consider two cases: $y<1$ and $y \geq 1$.
If $y<1$, then $x+y \leq (i+1)y < i+1$.
Furthermore, if $x+y < i$, then $(N-2)y-N+i+1 = (N-2)(y-1)+i-1 < (i-1)(y-1)+i-1 = (i-1)y \leq x$.
If $y \geq 1$, then $x-i \leq i(y-1) \leq (N-2)(y-1)$, so we again have $x+y<i+1$.
It also holds that $x+y \geq iy \geq i$, so $x>(N-2)y-N+i+1$.

On the other hand, for points $(x,y) \in C(V_i)$ satisfying $x+y<i+1$ and $x>(N-2)y-N+i+1$, it is easily checked that $(x,y)-(k,1) \not\in C(\A)$ for all $-1 \leq k \leq N-2$.

It follows that the apex points are the points with fractional part $\bbeta$ lying in the dark gray area in Figure~\ref{fig:collinear_point_configurations_apexpoints}.
Hence there are $N-1$ apex points if and only if $-\beta_1+(N-2)\beta_2<0$ and $\beta_1+\beta_2<1$, or $-\beta_1+(N-2)\beta_2 \geq N-3$ and $\beta_1+\beta_2 \geq 1$.
\end{proof}

A graphical interpretation of the interlacing condition is given in Figure~\ref{fig:interlacing_condition_collinear_points}.
Now we have found the interlacing condition, we can determine the irreducible algebraic functions:

\begin{thm}\label{thm:collinear_point_configurations}
Let $N \geq 2$ and $\A=\{(k,1) | -1 \leq k \leq N-2\}$.
Then $\HA$ is irreducible and has algebraic solutions if and only if $\bbeta \pmod \Z$ is a conjugate of one of the tuples in Table~\ref{tab:collinear_point_configurations}.
\end{thm}

\begin{table}
\centering
\caption{The parameters $\bbeta$ such that $\HA$ is irreducible and has algebraic solutions if $\A=\{(k,1) | -1 \leq k \leq N-2\}$} 
\label{tab:collinear_point_configurations} 
\renewcommand{\arraystretch}{1.5}
\begin{tabular}{l@{\hspace{0.8cm}}l@{\hspace{0.6cm}}l@{\hspace{0.6cm}}l@{\hspace{0.6cm}}l@{\hspace{0.6cm}}l@{\hspace{0.6cm}}l@{\hspace{0.6cm}}l@{\hspace{0.6cm}}l} 
\toprule
$N$ & \multicolumn{8}{@{}l}{$\bbeta$} \\
\midrule
All $N \geq 2$ &
\multicolumn{8}{@{}l}{$(r, 0)$ with $r \not\in \Z$} \\

$2$ & \multicolumn{8}{@{}l}{All irreducible functions are algebraic} \\
$3$ & 
$(\inv{3}, \inv{4})$ & 
$(\inv{3}, \frac{5}{6})$ & 
$(\inv{4}, \inv{6})$ & 
$(\inv{5}, \inv{10})$ &

$(\inv{3}, \inv{6})$ & 
$(\inv{3}, \inv{10})$ &
$(\inv{5}, \inv{6})$ &
$(\inv{5}, \frac{9}{10})$ \\

$4$ &
$(\half, \inv{6})$ & 
$(\inv{3}, \frac{5}{6})$ \\ 

$5$ &
$(\inv{3}, \frac{5}{6})$ \\
\bottomrule
\end{tabular}
\renewcommand{\arraystretch}{1}
\end{table}

Except for $N=3$, these results can also be found in~\cite{schipper_thesis}.

\begin{proof}
If $N=2$, then the lattice $\lat$ is trivial, so the $\Gamma$-series solutions of $H_\A(\bbeta)$ are monomials $z_1^{-\beta_1} z_2^{\beta_1+\beta_2}$, and hence algebraic.

If $N=3$, the Gauss function $\gauss(\frac{-\beta_1-\beta_2}{2}, \frac{-\beta_1-\beta_2-1}{2}, -\beta_1+1|4z)$ is a solution.
It is irreducible if and only if $H_\A(\bbeta)$ is irreducible.
Checking the interlacing condition for all parameters such that the Gauss function is algebraic (Table~\ref{tab:gauss_abc_orbits}) shows that $H_\A(\bbeta)$ has irreducible algebraic solutions if and only if $\bbeta$ is of the form $(\half,r)$ or $(r,0)$ with $r \not\in \Z$ or is one of the 36 conjugates of $(\inv{3}, \inv{4})$, $(\inv{3}, \inv{6})$, $(\inv{3}, \frac{5}{6})$, $(\inv{3}, \inv{10})$, $(\inv{4}, \inv{6})$, $(\inv{5}, \inv{6})$, $(\inv{5}, \inv{10})$ and $(\inv{5}, \frac{9}{10})$.

For $N=4$, we have the Horn $G_3$ function, which is irreducible and algebraic if and only if $\bbeta \in \{(r,0), \pm(\half, \inv{6}), \pm(\inv{3}, \frac{5}{6})\}$ with $r \not\in \Z$ (see~\cite{schipper_thesis} or~\cite{bod_algebraic_appell_horn}).

Now let $N=5$.
By Proposition~\ref{prop:subsets_algebraic}, if $H_\A(\bbeta)$ has irreducible algebraic solutions, then it is also irreducible for $N=4$ and has algebraic solutions.
Hence we only have to check $(r,0), \pm(\half, \inv{6})$ and $\pm(\inv{3}, \frac{5}{6})$.
It turns out that $(r,0)$ is a solution for all $r \not\in \Z$, as well as $\pm (\inv{3},\frac{5}{6})$, but $\pm(\half, \inv{6})$ is resonant. 

For $N=6$, we only have to check the solutions $(r,0)$ and $\pm (\inv{3}, \frac{5}{6})$ of $N=5$.
Now $\pm (\inv{3}, \frac{5}{6})$ is resonant, but $(r,0)$ is a solution.

Finally, for $N > 6$, the tuple $(r,0)$ gives irreducible functions if $r \not\in \Z$ and it satisfies the interlacing condition.
\end{proof}


\section{Polygons with at most 2 interior points}\label{sec:polygons-012_interior_points}

In the remainder of this paper, we will assume that $\A$ is a normal subset of $\Z^3$ and identify sets $\A$ with the corresponding polygons $P(\A)$.
In this section, we will determine the algebraic functions corresponding to polygons with at most 2 lattice interior points.
In the remainder of this paper, we will abbreviate `interior lattice point' and `lattice point on the boundary' to `interior point' and `boundary point', respectively.

\subsection*{Polygons without interior points}\label{subsec:polygons-012_interior_points-0_interior_points}

Rabinowitz has given a classification of the polygons without interior points~\cite{rabinowitz_census_lattice_polygons_one_point}.
Up to isomorphism, there are three types:
triangles with vertices $(0,0), (p,0)$ and $(0,1)$, the triangle with vertices $(0,0), (2,0)$ and $(0,2)$, and trapezoids with vertices $(0,0),(p,0),(q,1)$ and $(0,1)$.
The first type of triangles gives pyramidal sets $\A$, which we excluded from our considerations.
For the other polygons, we denote the corresponding sets $\A$ by $\A_1$ and $\A_{p,q}^{(1)}$, respectively.

\noindent
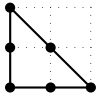
\begin{minipage}[c]{0.1\linewidth}
\psscalebox{0.53}{
\begin{pspicture}(2,2)
\psgrid[subgriddiv=1,gridlabels=0,gridwidth=0.2pt,griddots=5](0,0)(2,2)
\psdots[dotscale=2](0,0)(1,0)(2,0)(0,1)(1,1)(0,2)
\psline[linewidth=1.5pt](0,0)(2,0)(0,2)(0,0)
\end{pspicture}}
\end{minipage}
\hfill
\begin{minipage}[c]{0.85\linewidth}
\begin{lem}\label{lem:exceptional_triangle}
Let $\A_1=\begin{pmatrix} 0 & 1 & 2 & 0 & 1 & 2 \\ 0 & 0 & 0 & 1 & 1 & 1 \\ 1 & 1 & 1 & 1 & 1 & 1 \end{pmatrix}$.
Then there are no $\bbeta$ such that $\HAi{1}$ has irreducible algebraic solutions.
\end{lem}
\end{minipage}

\begin{proof}
Using the algorithm described in Section~\ref{sec:polygons-algebraic_hypergeometric_functions} one can easily show that there are never four apexpoints.
\end{proof}

\noindent
\begin{minipage}[c]{0.25\linewidth}
\psscalebox{0.53}{
\begin{pspicture}(6,1)
\psgrid[subgriddiv=1,gridlabels=0,gridwidth=0.2pt,griddots=5](0,0)(6,1)
\psdots[dotscale=2](0,0)(1,0)(3,0)(4,0)(6,0)(0,1)(1,1)(3,1)
\psline[linewidth=1.5pt](1.5,1)(0,1)(0,0)(1.5,0)
\psline[linewidth=1.5pt](2.5,0)(4.5,0)
\psline[linewidth=1.5pt](5.5,0)(6,0)(3,1)(2.5,1)
\psdots[dotsize=1pt](1.7,0)(1.9,0)(2.1,0)(2.3,0)(4.7,0)(4.9,0)(5.1,0)(5.3,0)(1.7,1)(1.9,1)(2.1,1)(2.3,1)
\end{pspicture}}
\end{minipage}
\hfill
\begin{minipage}[c]{0.7\linewidth}
\begin{lem}\label{lem:trapezoid}
Let 
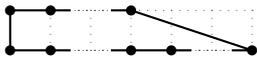
\begin{equation*}
\A_{p,q}^{(1)}=\begin{pmatrix} 0 & 1 & 2 & \ld & p & 0 & 1 & 2 & \ld & q \\ 0 & 0 & 0 & \ld & 0 & 1 & 1 & 1 & \ld & 1 \\ 1 & 1 & 1 & \ld & 1 & 1 & 1 & 1 & \ld & 1 \end{pmatrix}
\end{equation*}
with $p \geq q$.
Then $H_{\A_{p,q}^{(1)}}(\bbeta)$ has irreducible algebraic solutions if and only if, up to conjugation and equivalence modulo $\Z$, $\bbeta$ is one of the tuples in Table~\ref{tab:solutions_trapezoid}.
\end{lem}
\end{minipage}

\begin{table} 
\centering
\caption{The parameters $\bbeta$ such that $H_{\A_{p,q}^{(1)}}(\bbeta)$ has irreducible algebraic solutions} 
\label{tab:solutions_trapezoid}
\renewcommand{\arraystretch}{1.5}
\begin{tabular}{l@{\hspace{0.7cm}}l@{\hspace{0.6cm}}l@{\hspace{0.6cm}}l@{\hspace{0.6cm}}l@{\hspace{0.6cm}}l@{\hspace{0.6cm}}l} 
\toprule
$(p,q)$ & \multicolumn{6}{@{}l}{$\bbeta$} \\
\toprule
$(1,1)$ & 
\multicolumn{6}{@{}l}{Up to permutations of $\{\beta_1, \beta_2\}$:} \\

 & \multicolumn{6}{@{}l}{$(r, -r, \frac{1}{2})$,
$(r, r+\frac{1}{2}, \frac{1}{2})$ \tr{and}
$(r, r+\frac{1}{2}, -2r)$ 
with $2r \not\in \Z$} \\

 & $(\frac{1}{2}, \frac{1}{6}, \frac{2}{3})$ &
$(\frac{1}{6}, \frac{5}{12}, \frac{2}{3})$ & 
$(\frac{1}{10}, \frac{13}{30}, \frac{2}{3})$ &
$(\frac{1}{15}, \frac{11}{15}, \frac{4}{5})$ &
$(\frac{1}{24}, \frac{17}{24}, \frac{1}{2})$ &
$(\frac{1}{60}, \frac{31}{60}, \frac{4}{5})$ \\

 & $(\frac{1}{4}, \frac{3}{4}, \frac{1}{3})$ &
$(\frac{1}{6}, \frac{5}{12}, \frac{3}{4})$ &
$(\frac{1}{10}, \frac{13}{30}, \frac{4}{5})$ &
$(\frac{1}{20}, \frac{11}{20}, \frac{3}{5})$ &
$(\frac{1}{24}, \frac{17}{24}, \frac{3}{4})$ &
$(\frac{1}{60}, \frac{41}{60}, \frac{1}{2})$ \\

 & $(\frac{1}{4}, \frac{7}{12}, \frac{1}{2})$ &
$(\frac{1}{6}, \frac{11}{30}, \frac{2}{3})$ &
$(\frac{1}{12}, \frac{5}{12}, \frac{3}{4})$ &
$(\frac{1}{20}, \frac{11}{20}, \frac{4}{5})$ &
$(\frac{1}{24}, \frac{19}{24}, \frac{1}{2})$ &
$(\frac{1}{60}, \frac{41}{60}, \frac{4}{5})$ \\

 & $(\frac{1}{4}, \frac{7}{12}, \frac{2}{3})$ &
$(\frac{1}{6}, \frac{11}{30}, \frac{4}{5})$ &
$(\frac{1}{12}, \frac{7}{12}, \frac{2}{3})$ &
$(\frac{1}{20}, \frac{13}{20}, \frac{1}{2})$ &
$(\frac{1}{24}, \frac{19}{24}, \frac{2}{3})$ &
$(\frac{1}{60}, \frac{49}{60}, \frac{1}{2})$ \\

 & $(\frac{1}{6}, \frac{5}{6}, \frac{1}{3})$ &
$(\frac{1}{10}, \frac{3}{10}, \frac{4}{5})$ &
$(\frac{1}{15}, \frac{7}{15}, \frac{2}{3})$ &
$(\frac{1}{20}, \frac{13}{20}, \frac{4}{5})$ &
$(\frac{1}{30}, \frac{11}{30}, \frac{4}{5})$ &
$(\frac{1}{60}, \frac{49}{60}, \frac{2}{3})$ \\

 & $(\frac{1}{6}, \frac{5}{6}, \frac{1}{4})$ &
$(\frac{1}{10}, \frac{9}{10}, \frac{1}{3})$ &
$(\frac{1}{15}, \frac{7}{15}, \frac{4}{5})$ &
$(\frac{1}{24}, \frac{13}{24}, \frac{2}{3})$ &
$(\frac{1}{30}, \frac{19}{30}, \frac{2}{3})$ \\ 

 & $(\frac{1}{6}, \frac{5}{6}, \frac{1}{5})$ &
$(\frac{1}{10}, \frac{9}{10}, \frac{1}{5})$ &
$(\frac{1}{15}, \frac{11}{15}, \frac{2}{5})$ &
$(\frac{1}{24}, \frac{13}{24}, \frac{3}{4})$ &
$(\frac{1}{60}, \frac{31}{60}, \frac{2}{3})$ \\

$(2,1)$ &  
$(\frac{1}{3}, \frac{5}{6}, \frac{2}{3})$ & 
$(\frac{1}{6}, \frac{2}{3}, \half)$ & 
$(\frac{1}{6}, \frac{5}{6}, \frac{2}{3})$ \\

$(2,2)$ &  
$(\frac{1}{6}, \frac{5}{6}, \frac{2}{3})$ \\

$(3,1)$ &  
$(\frac{1}{6}, \frac{5}{6}, \frac{2}{3})$ \\
\bottomrule
\end{tabular}
\renewcommand{\arraystretch}{1}
\end{table}

\begin{proof}
Note that $\A_{1,1}^{(1)}$ is isomorphic to the set $\A$ for the Gauss function $\gauss$.
The parameters such that $\gauss$ is irreducible and algebraic are well-known and can for example be found in Table~\ref{tab:gauss_abc_orbits}.
In~\cite{bod_algebraic_appell_horn}, the irreducible algebraic Horn $G_1$ functions are determined.
For this function, the set $\A$ is isomorphic to $\A_{2,1}^{(1)}$.
This gives the first two cases in Table~\ref{tab:solutions_trapezoid}.

Let $p=q=2$.
The interlacing condition is $(\entier{-\beta_1+2\beta_3}, \entier{-\beta_2+\beta_3}) \in \{(-1,0), (1,-1)\}$.
It is clear that $\A_{2,1}^{(1)}$ is included in $\A_{2,2}^{(1)}$.
Hence we only have to check the interlacing condition for $\pm(\frac{1}{3}, \frac{5}{6}, \frac{2}{3})$, $\pm(\frac{1}{6}, \frac{2}{3}, \half)$ and $\pm(\frac{1}{6}, \frac{5}{6}, \frac{2}{3})$.
It turns out that this condition is satisfied only for $\pm(\frac{1}{6}, \frac{5}{6}, \frac{2}{3})$.

For $p=3$ and $q=1$, we again have the inclusion $\A_{2,1}^{(1)} \subseteq \A_{3,1}^{(1)}$.
The interlacing condition is $(\entier{-\beta_1-2\beta_2+3\beta_3}, \entier{-\beta_2+\beta_3}) \in \{(-1,0), (0,-1)\}$.
Checking $\pm(\frac{1}{3}, \frac{5}{6}, \frac{2}{3})$, $\pm(\frac{1}{6}, \frac{2}{3}, \half)$ and $\pm(\frac{1}{6}, \frac{5}{6}, \frac{2}{3})$ shows that only $\pm(\frac{1}{6}, \frac{5}{6}, \frac{2}{3})$ gives an irreducible algebraic function. 

For $(p,q)=(3,2)$, the interlacing condition is $(\entier{-\beta_1-\beta_2+3\beta_3}, \entier{-\beta_2+\beta_3}) \in \{(-1,0), (1,-1)\}$.
Due to the inclusion $\A_{3,1}^{(1)} \subseteq \A_{3,2}^{(1)}$, we only have to check this condition for $\pm(\frac{1}{6}, \frac{5}{6}, \frac{2}{3})$.
This tuple doesn't satisfy the condition, so there are no irreducible algebraic functions.
It follows immediately from this that there are also no irreducible algebraic functions for any $\A_{p,q}^{(1)}$ with $p \geq 3$ and $q \geq 2$.

This leaves us with the case $p \geq 4$ and $q=1$.
Let $p=4$.
The interlacing condition is $(\entier{-\beta_1-3\beta_2+4\beta_3}, \entier{-\beta_2+\beta_3}) \in \{(-1,0), (0,-1)\}$ and we only have to check $\pm(\frac{1}{6}, \frac{5}{6}, \frac{2}{3})$.
Again there are no irreducible algebraic functions.
Hence there are also no irreducible algebraic functions with $p>4$ and $q=1$.
\end{proof}

\subsection*{Polygons with exactly one interior point}\label{subsec:polygons-012_interior_points-1_interior_point}

\begin{figure}
\centering
\psscalebox{0.53}{
\begin{pspicture}(28,8)
\psgrid[subgriddiv=1,gridlabels=0,gridwidth=0.2pt,griddots=5](28,8)

\pspolygon[fillstyle=solid,linewidth=0pt,fillcolor=lightgray](1,5)(2,7)(3,6)
\psdots[dotscale=2](1,5)(3,6)(2,7)
\psdots[dotscale=2,dotstyle=o](2,6)
\psline[linewidth=1.5pt](1,5)(2,7)(3,6)(1,5)

\pspolygon[fillstyle=solid,linewidth=0pt,fillcolor=lightgray](4,5)(6,6)(4,7)
\psdots[dotscale=2](4,5)(4,6)(6,6)(4,7)
\psdots[dotscale=2,dotstyle=o](5,6)
\psline[linewidth=1.5pt](4,5)(6,6)(4,7)(4,5)

\pspolygon[fillstyle=solid,linewidth=0pt,fillcolor=lightgray](8,5)(9,6)(7,7)(7,6)
\psdots[dotscale=2](8,5)(7,6)(9,6)(7,7)
\psdots[dotscale=2,dotstyle=o](8,6)
\psline[linewidth=1.5pt](8,5)(9,6)(7,7)(7,6)(8,5)

\pspolygon[fillstyle=solid,linewidth=0pt,fillcolor=lightgray](12,5)(12,6)(10,7)(10,6)
\psdots[dotscale=2](12,5)(10,6)(12,6)(10,7)
\psdots[dotscale=2,dotstyle=o](11,6)
\psline[linewidth=1.5pt](12,5)(12,6)(10,7)(10,6)(12,5)

\pspolygon[fillstyle=solid,linewidth=0pt,fillcolor=lightgray](13,5)(15,5)(15,6)(14,7)
\psdots[dotscale=2](13,5)(14,5)(15,5)(15,6)(14,7)
\psdots[dotscale=2,dotstyle=o](14,6)
\psline[linewidth=1.5pt](13,5)(15,5)(15,6)(14,7)(13,5)
\psline[linewidth=1.5pt,linestyle=dashed](13,5)(15,6)

\pspolygon[fillstyle=solid,linewidth=0pt,fillcolor=lightgray](17,5)(18,5)(18,6)(17,7)(16,6)
\psdots[dotscale=2](17,5)(18,5)(16,6)(18,6)(17,7)
\psdots[dotscale=2,dotstyle=o](17,6)
\psline[linewidth=1.5pt](17,5)(18,5)(18,6)(17,7)(16,6)(17,5)
\psline[linewidth=1.5pt,linestyle=dashed](18,5)(17,7)

\pspolygon[fillstyle=solid,linewidth=0pt,fillcolor=lightgray](19,5)(21,5)(21,6)(20,7)(19,6)
\psdots[dotscale=2](19,5)(20,5)(21,5)(19,6)(21,6)(20,7)
\psdots[dotscale=2,dotstyle=o](20,6)
\psline[linewidth=1.5pt](19,5)(21,5)(21,6)(20,7)(19,6)(19,5)
\psline[linewidth=1.5pt,linestyle=dashed](19,5)(20,7)

\pspolygon[fillstyle=solid,linewidth=0pt,fillcolor=lightgray](23,5)(24,5)(24,6)(23,7)(22,7)(22,6)
\psdots[dotscale=2](23,5)(24,5)(22,6)(24,6)(22,7)(23,7)
\psdots[dotscale=2,dotstyle=o](23,6)
\psline[linewidth=1.5pt](23,5)(24,5)(24,6)(23,7)(22,7)(22,6)(23,5)
\psline[linewidth=1.5pt,linestyle=dashed](22,6)(23,7)

\psdots[dotscale=2](25,5)(26,5)(27,5)(27,6)(26,7)(27,7)
\psdots[dotscale=2,dotstyle=o](26,6)
\psline[linewidth=1.5pt](25,5)(27,5)(27,7)(26,7)(25,5)

\psdots[dotscale=2](1,1)(2,1)(3,1)(1,2)(3,2)(2,3)(3,3)
\psdots[dotscale=2,dotstyle=o](2,2)
\psline[linewidth=1.5pt](1,1)(3,1)(3,3)(2,3)(1,2)(1,1)

\psdots[dotscale=2](4,1)(5,1)(6,1)(4,2)(6,2)(4,3)(5,3)(6,3)
\psdots[dotscale=2,dotstyle=o](5,2)
\psline[linewidth=1.5pt](4,1)(6,1)(6,3)(4,3)(4,1)

\psdots[dotscale=2](7,1)(8,1)(9,1)(10,1)(9,2)(8,3)
\psdots[dotscale=2,dotstyle=o](8,2)
\psline[linewidth=1.5pt](7,1)(10,1)(8,3)(7,1) 

\psdots[dotscale=2](11,1)(12,1)(13,1)(14,1)(11,2)(13,2)(12,3)
\psdots[dotscale=2,dotstyle=o](12,2)
\psline[linewidth=1.5pt](11,1)(14,1)(12,3)(11,2)(11,1)

\psdots[dotscale=2](15,1)(16,1)(17,1)(18,1)(15,2)(17,2)(15,3)(16,3)
\psdots[dotscale=2,dotstyle=o](16,2)
\psline[linewidth=1.5pt](15,1)(18,1)(16,3)(15,3)(15,1)

\psdots[dotscale=2](19,1)(20,1)(21,1)(22,1)(23,1)(20,2)(22,2)(21,3)
\psdots[dotscale=2,dotstyle=o](21,2)
\psline[linewidth=1.5pt](19,1)(23,1)(21,3)(19,1)

\psdots[dotscale=2](24,1)(25,1)(26,1)(27,1)(24,2)(26,2)(24,3)(25,3)(24,4)
\psdots[dotscale=2,dotstyle=o](25,2)
\psline[linewidth=1.5pt](24,1)(27,1)(24,4)(24,1)
\end{pspicture}}
\caption{The polygons with exactly one interior point
\label{fig:polygons_1interior_point}}
\end{figure}
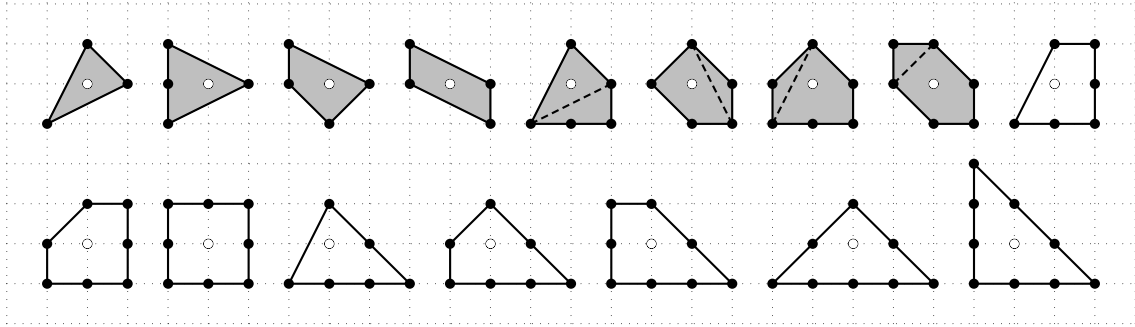

In~\cite{scott_convex_lattice_polygons}, Scott proves that there exists a lattice polygon with $i$ interior points and $b$ boundary points if and only if either $i=0$, or $i=1$ and $3 \leq b \leq 9$, or $i \geq 2$ and $3 \leq b \leq 2i+6$.
Furthermore, Theorem~2 in~\cite{lagarias_ziegler_bounds_lattice_polytopes} states that a polygon of normalized area $V$ is, up to isomorphism, contained in a square of side length $2V$.
Since the number of interior and boundary points determine the area by Pick's formula, this implies that there are only finitely many non-isomorphic polygons with a given number of interior and boundary points.
Polygons with exactly one interior point have 3 to 9 boundary points.
A classification of these can be found in both~\cite{rabinowitz_census_lattice_polygons_one_point} and~\cite{poonen_rodriguez_lattice_polygons}.
There are 16 isomorphism classes.
They are shown in Figure~\ref{fig:polygons_1interior_point}.

For each of these, we can compute the interlacing condition.
Then we compute the irreducible algebraic functions by using a reduction to a polygon for which we know the algebraic functions already, as in Corollary~\ref{cor:extension_less_algebraic_functions_resonant}.
For the non-shaded polygons in Figure~\ref{fig:polygons_1interior_point}, one easily computes that there exists no $\bbeta \in \Q^3$ such that $\sign = \vol(Q(\A))$.
Hence $\HA$ never has irreducible algebraic solutions.
We now consider the shaded polygons one by one.
We denote the corresponding sets $\A$ by $\A_2$ up to $\A_9$.

\noindent 
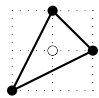
\begin{minipage}[c]{0.1\linewidth}
\psscalebox{0.53}{
\begin{pspicture}(2,2)
\psgrid[subgriddiv=1,gridlabels=0,gridwidth=0.2pt,griddots=5](0,0)(2,2)
\psdots[dotscale=2](0,0)(2,1)(1,2)
\psdots[dotscale=2,dotstyle=o](1,1)
\psline[linewidth=1.5pt](0,0)(1,2)(2,1)(0,0)
\end{pspicture}}
\end{minipage}
\hfill
\begin{minipage}[c]{0.85\linewidth}
\begin{lem}\label{lem:P1}
Let $\A_2=\begin{pmatrix} -1 & 0 & 1 & 0 \\ -1 & 0 & 0 & 1 \\ 1 & 1 & 1 & 1 \end{pmatrix}$.
Then $\HAi{2}$ has irreducible algebraic solutions if and only if, up to conjugation and equivalence modulo $\Z$, $\bbeta$ is  one of the tuples in Table~\ref{tab:solutionsP1}.
\end{lem}
\end{minipage}

\begin{table} 
\centering
\caption{The parameters $\bbeta$ such that $\HAi{2}$ has irreducible algebraic solutions} 
\label{tab:solutionsP1}
\renewcommand{\arraystretch}{1.5}
\begin{tabular}{l@{\hspace{0.7cm}}l@{\hspace{0.7cm}}l@{\hspace{0.7cm}}l@{\hspace{0.7cm}}l@{\hspace{0.7cm}}l@{\hspace{0.7cm}}l} 
\toprule
\multicolumn{7}{l}{$(\frac{1}{3}, \frac{2}{3}, r)$  
with $r \not\in \Z$} \\

\multicolumn{7}{l}{$(r, \half, \half)$, 
$(\half, r, \half)$ \tr{and}
$(r, r+\half, \half)$ 
with $2r \not\in \Z$} \\

$(\frac{1}{2}, \frac{1}{3}, \frac{2}{3})$ & 
$(\frac{1}{2}, \frac{1}{6}, \frac{1}{3})$ &  
$(\frac{1}{3}, \frac{5}{6}, \frac{1}{3})$ &
$(\frac{1}{4}, \frac{3}{4}, \frac{1}{3})$ & 
$(\frac{1}{5}, \frac{3}{5}, \frac{1}{2})$ &
$(\frac{1}{6}, \frac{1}{2}, \frac{1}{3})$ &
$(\frac{1}{7}, \frac{3}{7}, \frac{1}{2})$ \\

$(\frac{1}{2}, \frac{1}{4}, \frac{1}{3})$ & 
$(\frac{1}{3}, \frac{1}{2}, \frac{2}{3})$ & 
$(\frac{1}{4}, \frac{1}{2}, \frac{1}{3})$ &
$(\frac{1}{5}, \frac{2}{5}, \frac{1}{2})$ &
$(\frac{1}{5}, \frac{4}{5}, \frac{1}{2})$ &
$(\frac{1}{6}, \frac{2}{3}, \frac{2}{3})$ &
$(\frac{1}{7}, \frac{5}{7}, \frac{1}{2})$ \\
\bottomrule
\end{tabular}
\renewcommand{\arraystretch}{1}
\end{table}

\begin{proof}
The lattice is given by $\lat = \Z (1,-3,1,1)$.
The $\Gamma$-series
\begin{equation*}
\Phi(z_1,z_2,z_3,z_4) = \sum_{n \in \Z} \frac{z_1^{n+\gamma_1} z_2^{-3n+\gamma_2} z_3^{n+\gamma_3} z_4^{n+\gamma_4}}{\Gamma(1+n+\gamma_1) \Gamma(1-3n+\gamma_2) \Gamma(1+n+\gamma_3) \Gamma(1+n+\gamma_4)}
\end{equation*}
is a formal solution of $\HAi{2}$ whenever $\A_2 \vect{\gamma} = \bbeta$.
Choosing $\vect{\gamma} = (-\beta_2, -\beta_1+2\beta_2+\beta_3, \beta_1-\beta_2, 0)$ gives a convergent solution.
Note that $\Phi$ is irreducible and algebraic if and only if the higher hypergeometric function $\mbox{$_3 F_2$}(\frac{-\gamma_2}{3}, \frac{-\gamma_2+1}{3}, \frac{-\gamma_2+2}{3}; \gamma_1+1, \gamma_3+1 | z)$ is irreducible and algebraic.
The algebraic higher hypergeometric functions have been determined by Beukers and Heckman~\cite{beukers_heckman_monodromy_nfn-1}.
From the irreducible algebraic $\mbox{$_3 F_2$}$ functions in~\cite{beukers_heckman_monodromy_nfn-1}, we select the functions whose first 3 parameters differ by $\inv{3}$ and compute $\bbeta=(-\gamma_1+\gamma_3, -\gamma_1, \gamma_1+\gamma_2+\gamma_3)$.
This gives the tuples in Table~\ref{tab:solutionsP1}.
\end{proof}

\noindent 
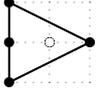
\begin{minipage}[c]{0.1\linewidth}
\psscalebox{0.53}{
\begin{pspicture}(2,2)
\psgrid[subgriddiv=1,gridlabels=0,gridwidth=0.2pt,griddots=5](0,0)(2,2)
\psdots[dotscale=2](0,0)(0,1)(2,1)(0,2)
\psdots[dotscale=2,dotstyle=o](1,1)
\psline[linewidth=1.5pt](0,0)(2,1)(0,2)(0,0)
\end{pspicture}}
\end{minipage}
\hfill
\begin{minipage}[c]{0.85\linewidth}
\begin{lem}\label{lem:P2} 
Let $\A_3=\begin{pmatrix} -1 & -1 & 0 & 1 & -1 \\ -1 & 0 & 0 & 0 & 1 \\ 1 & 1 & 1 & 1 & 1 \end{pmatrix}$.
Then $\HAi{3}$ has irreducible algebraic solutions if and only if, up to conjugation and equivalence modulo $\Z$, $\bbeta$ is one of the tuples in Table~\ref{tab:solutionsP2}.
\end{lem}
\end{minipage}

\begin{table} 
\centering
\caption{The parameters $\bbeta$ such that $\HAi{3}$ has irreducible algebraic solutions} 
\label{tab:solutionsP2}
\renewcommand{\arraystretch}{1.5}
\begin{tabular}{l@{\hspace{0.7cm}}l@{\hspace{0.7cm}}l@{\hspace{0.7cm}}l@{\hspace{0.7cm}}l@{\hspace{0.7cm}}l} 
\toprule
\multicolumn{3}{l}{$(0,\half,r)$ 
with $r \not\in \Z$}  &
\multicolumn{3}{l}{$(r,\half,\half)$ 
with $2r \not\in 2\Z+1$} \\

$(0, \frac{1}{3}, \frac{1}{2})$ &
$(\frac{1}{6}, \frac{1}{2}, \frac{1}{3})$ & 
$(\frac{1}{6}, \frac{1}{2}, \frac{1}{4})$ & 
$(\frac{1}{6}, \frac{1}{3}, \frac{2}{3})$ & 
$(\frac{1}{10}, \frac{1}{2}, \frac{1}{3})$ & 
$(\frac{1}{10}, \frac{1}{2}, \frac{4}{5})$ \\ 

$(\frac{1}{4}, \frac{1}{2}, \frac{1}{3})$ & 
$(\frac{1}{6}, \frac{1}{2}, \frac{2}{3})$ & 
$(\frac{1}{6}, \frac{1}{2}, \frac{1}{5})$ &
$(\frac{1}{6}, \frac{2}{3}, \frac{2}{3})$ &
$(\frac{1}{10}, \frac{1}{2}, \frac{1}{5})$ \\ 
\bottomrule
\end{tabular}
\renewcommand{\arraystretch}{1}
\end{table}

\begin{proof}
With $\lat = \Z(0,1,-2,1,0) \oplus \Z(1,-2,0,0,1)$ and $\vect{\gamma}=(-\beta_2,-\beta_1+\beta_2,\beta_1+\beta_3,0,0)$, we get the $\Gamma$-series
\begin{equation*}
\Phi(\vect{z}) = \!\! \sum_{m,n \geq 0} \! \frac{z_1^{n-\beta_2} z_2^{m-2n-\beta_1+\beta_2} z_3^{-2m+\beta_1+\beta_3} z_4^m z_5^n}{\Gamma(1+n-\beta_2) \Gamma(1+m-2n-\beta_1+\beta_2) \Gamma(1-2m+\beta_1+\beta_3) m! n!},
\end{equation*}
which is irreducible and algebraic if and only if 
\begin{equation*}
\Psi(x,y) = \sum_{m,n \geq 0} \frac{\poch{\beta_1-\beta_2}{-m+2n} \poch{-\beta_1-\beta_3}{2m}}{\poch{-\beta_2+1}{n} m! n!} x^m y^n
\end{equation*}
is irreducible and algebraic.
Then $\Psi_x(x)=\Psi(x,0)$ and $\Phi_y(y)=\Psi(0,y)$ are also algebraic.
We have $\Psi_x(x) = \gauss(\frac{-\beta_1-\beta_3}{2}, \frac{-\beta_1-\beta_3+1}{2}, -\beta_1+\beta_2+1 | 4x)$ and \mbox{$\Psi_y(y) = \gauss(\frac{\beta_1-\beta_2}{2}, \frac{\beta_1-\beta_2+1}{2}, -\beta_2+1 | 4y)$}.

$\HAi{3}$ is irreducible if and only if $-\beta_1-2\beta_2+\beta_3, -\beta_1+2\beta_2+\beta_3, \beta_1+\beta_3 \not\in \Z$.
As $\gauss(\alpha_1, \alpha_2, \alpha_3|z)$ is irreducible if and only if $\alpha_1, \alpha_2, \alpha_1-\alpha_3, \alpha_2-\alpha_3 \not\in \Z$, irreducibility of $\HAi{2}$ implies irreducibility of $\Phi_x$.
However, $\Phi_y$ is reducible if $\beta_1 \pm \beta_2 \in \Z$.
But in this case, $\HAi{3}$ is not totally non-resonant, so this doesn't give any algebraic functions.
Hence if $\HAi{3}$ has irreducible and algebraic solutions, then both $\Psi_x$ and $\Psi_y$ are irreducible and algebraic.

The tuples $(\alpha_1, \alpha_2, \alpha_3)$ such that $\gauss(\alpha_1, \alpha_2, \alpha_3|z)$ is irreducible and algebraic can be found in Table~\ref{tab:gauss_abc_orbits}.
We select the pairs $(\vect{\alpha}^{(1)}, \vect{\alpha}^{(2)})$ of tuples satisfying $\alpha_2^{(i)} \equiv \alpha_1^{(i)}+\half \pmod \Z$ and $\alpha_3^{(1)} + 2\alpha_1^{(2)} \in \Z$, and compute the corresponding $\bbeta$.
Then we check the interlacing condition for $\A_3$. 
This is given by $(\entier{-2\beta_1-\beta_2+\beta_3}, \entier{\beta_2+\beta_3}, \entier{2\beta_1-\beta_2+\beta_3}) \in \{(-2,1,0),(-1,0,1)\}$.

The triples $\vect{\alpha}^{(i)}$ can be either of the form $(r,r+\half,\half)$ or $(r,r+\half,2r)$, or can be one of the other 408 triples for which the Gauss function is irreducible and algebraic.
Hence there are several cases to check for $(\vect{\alpha}^{(1)}, \vect{\alpha}^{(2)})$.
In most cases, one can easily show that there are only finitely many possibilities for the parameter(s), using the fact that $\alpha_3^{(1)} + 2\alpha_1^{(2)} \in \Z$.
We only discuss the case in which $\vect{\alpha}^{(1)}=(r,r+\half,\half)$ and $\vect{\alpha}^{(2)}$ is one of the other 408 triples.
Since $\alpha_3^{(1)} + 2\alpha_1^{(2)} \in \Z$, we have $\alpha_1^{(2)} = \pm \inv{4}$ and hence $\alpha_2^{(2)} = \mp \inv{4}$.
Checking all 408 triples, we get $\alpha_3^{(2)} = \pm \inv{3}$.
Hence we find $\bbeta = \pm (\inv{6}, \frac{2}{3}, s)$ for some parameter $s$.
This satisfies the interlacing condition if and only if $\half \leq s < \frac{5}{6}$ (for $\bbeta = (\inv{6}, \frac{2}{3}, s)$) or $\inv{6} \leq s < \half$ (for $\bbeta = (\frac{5}{6}, \frac{1}{3}, s)$).
Suppose that the interlacing condition is satisfied for all conjugates.
Choose a conjugate $k \bbeta$ of $\bbeta$ such that $\fp{ks} = \inv{l}$ for some $l \in \Z_{\geq 2}$.
It follows from $\inv{6} \leq ks < \frac{5}{6}$ that $l|6$, and hence at the denominator of $s$ divides 6.
Now we easily find $\bbeta = (\inv{6}, \frac{2}{3}, \frac{2}{3})$.

Checking all cases, we find the 2 families and 42 other tuples in Table~\ref{tab:solutionsP2}.
\end{proof}

\noindent 
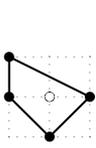
\begin{minipage}[c]{0.1\linewidth}
\psscalebox{0.53}{
\begin{pspicture}(2,2)
\psgrid[subgriddiv=1,gridlabels=0,gridwidth=0.2pt,griddots=5](0,0)(2,2)
\psdots[dotscale=2](1,0)(0,1)(2,1)(0,2)
\psdots[dotscale=2,dotstyle=o](1,1)
\psline[linewidth=1.5pt](1,0)(2,1)(0,2)(0,1)(1,0)
\end{pspicture}}
\end{minipage}
\hfill
\begin{minipage}[c]{0.85\linewidth}
\begin{lem}\label{lem:P3}
Let $\A_4=\begin{pmatrix} 0 & -1 & 0 & 1 & -1 \\ -1 & 0 & 0 & 0 & 1 \\ 1 & 1 & 1 & 1 & 1 \end{pmatrix}$.
Then $\HAi{4}$ has irreducible algebraic solutions if and only if $\bbeta = (r, \frac{1}{2}, \frac{1}{2})$ with $2r \not\in \Z$ or, up to conjugation and equivalence modulo $\Z$, $\bbeta$ equals one of the tuples 
$(\frac{1}{2}, \frac{1}{3}, \frac{1}{3})$, 
$(\frac{1}{3}, \frac{1}{3}, \frac{1}{2})$, 
$(\frac{1}{4}, \frac{1}{2}, \frac{1}{3})$,
$(\frac{1}{5}, \frac{3}{5}, \frac{1}{2})$ and
$(\frac{1}{6}, \frac{1}{3}, \frac{2}{3})$.
\end{lem}
\end{minipage}

\begin{proof}
The parameters such that the Horn $H_5$ function is irreducible and algebraic are determined in~\cite{bod_algebraic_appell_horn}.
$H_5$ is given by $\A_{H_5} = \begin{pmatrix} 1 & 0 & 0 & 2 & 1 \\ 0 & 1 & 0 & -1 & 1 \\ 0 & 0 & 1 & 0 & -1 \end{pmatrix}$.
The map $f(x,y,z)=(-y,-z,x+y+z)$ is an isomorphism of $\Z^3$, mapping $\A_{H_5}$ to $\A_4$.
Hence $\HAi{4}$ has irreducible algebraic solutions if and only if $H_{\A_{H_5}}(f^{-1}(\bbeta))$ has irreducible algebraic solutions.
\end{proof}

\noindent 
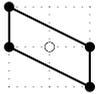
\begin{minipage}[c]{0.1\linewidth}
\psscalebox{0.53}{
\begin{pspicture}(2,2)
\psgrid[subgriddiv=1,gridlabels=0,gridwidth=0.2pt,griddots=5](0,0)(2,2)
\psdots[dotscale=2](2,0)(0,1)(2,1)(0,2)
\psdots[dotscale=2,dotstyle=o](1,1)
\psline[linewidth=1.5pt](2,0)(2,1)(0,2)(0,1)(2,0)
\end{pspicture}}
\end{minipage}
\hfill
\begin{minipage}[c]{0.8\linewidth}
\begin{lem}\label{lem:P4}
Let $\A_5=\begin{pmatrix} 1 & -1 & 0 & 1 & -1 \\ -1 & 0 & 0 & 0 & 1 \\ 1 & 1 & 1 & 1 & 1 \end{pmatrix}$.
Then $\HAi{5}$ has irreducible algebraic solutions if and only if, up to conjugation and equivalence modulo $\Z$, $\bbeta$ is one of the tuples in Table~\ref{tab:solutionsP4}.
\end{lem}
\end{minipage}

\begin{table} 
\centering
\caption{The parameters $\bbeta$ such that $\HAi{5}$ has irreducible algebraic solutions} 
\label{tab:solutionsP4}
\renewcommand{\arraystretch}{1.5}
\begin{tabular}{l@{\hspace{0.55cm}}l@{\hspace{0.55cm}}l@{\hspace{0.55cm}}l@{\hspace{0.55cm}}l@{\hspace{0.55cm}}l@{\hspace{0.55cm}}l@{\hspace{0.55cm}}l} 
\toprule
\multicolumn{4}{l}{$(0, \frac{1}{2}, r)$,
$(r, \frac{1}{2} - r, \frac{1}{2})$,  
with $r \not\in \Z$} &
\multicolumn{4}{l}{$(r, \half, \half)$, 
with $2r \not\in 2\Z+1$} \\

$(0, \frac{1}{3}, \frac{1}{2})$ & 
$(\frac{1}{4}, \frac{1}{4}, \frac{1}{3})$ & 
$(\frac{1}{6}, \frac{1}{2}, \frac{1}{3})$ & 
$(\frac{1}{6}, \frac{1}{2}, \frac{1}{5})$ & 
$(\frac{1}{6}, \frac{1}{3}, \frac{1}{4})$ & 
$(\frac{1}{10}, \frac{1}{2}, \frac{1}{5})$ & 
$(\frac{1}{10}, \frac{2}{5}, \frac{1}{5})$ &
$(\frac{1}{12}, \frac{1}{4}, \frac{1}{2})$ \\ 

$(\frac{1}{3}, \frac{1}{3}, \frac{1}{2})$ & 
$(\frac{1}{5}, \frac{1}{5}, \frac{1}{2})$ & 
$(\frac{1}{6}, \frac{1}{2}, \frac{2}{3})$ &
$(\frac{1}{6}, \frac{1}{3}, \frac{1}{3})$ & 
$(\frac{1}{6}, \frac{1}{3}, \frac{1}{5})$ & 
$(\frac{1}{10}, \frac{1}{2}, \frac{4}{5})$ & 
$(\frac{1}{10}, \frac{2}{5}, \frac{4}{5})$ & 
$(\frac{1}{15}, \frac{1}{3}, \frac{1}{2})$ \\

$(\frac{1}{4}, \frac{1}{2}, \frac{1}{3})$ & 
$(\frac{1}{5}, \frac{3}{5}, \frac{1}{2})$ &
$(\frac{1}{6}, \frac{1}{2}, \frac{1}{4})$ & 
$(\frac{1}{6}, \frac{1}{3}, \frac{2}{3})$ & 
$(\frac{1}{10}, \frac{1}{2}, \frac{1}{3})$ & 
$(\frac{1}{10}, \frac{2}{5}, \frac{1}{3})$ & 
$(\frac{1}{12}, \frac{2}{3}, \frac{1}{2})$ &
$(\frac{1}{15}, \frac{3}{5}, \frac{1}{2})$ \\
\bottomrule
\end{tabular}
\renewcommand{\arraystretch}{1}
\end{table}

\begin{proof}
We have $\lat=\Z (0,1,-2,1,0) \oplus \Z(1,0,-2,0,1)$.
Similar to the proof of Lemma~\ref{lem:P1}, one checks that the $\Gamma$-series with $\vect{\gamma}=(-\beta_2, -\beta_1-\beta_2,\beta_1+2\beta_2+\beta_3,0,0)$ is irreducible and algebraic if and only if the Appell $F_4$ function with parameters $\frac{-\beta_1-2\beta_2-\beta_3}{2}, \frac{-\beta_1-2\beta_2-\beta_3+1}{2}, 1-\beta_1-\beta_2, 1-\beta_2|x,y)$ is irreducible and algebraic.
The tuples $(a,b,c_1,c_2)$ such that $F_4(a,b,c_1,c_2|x,y)$ is irreducible and algebraic can be found in~\cite{bod_algebraic_appell_horn}.
We select the tuples satisfying $a-b \equiv \half \pmod \Z$ and compute the corresponding $\bbeta = (c_2-c_1, -c_2, -2a+c_1+c_2)$.
\end{proof}

\noindent 
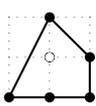
\begin{minipage}[c]{0.1\linewidth}
\psscalebox{0.53}{
\begin{pspicture}(2,2)
\psgrid[subgriddiv=1,gridlabels=0,gridwidth=0.2pt,griddots=5](0,0)(2,2)
\psdots[dotscale=2](0,0)(1,0)(2,0)(2,1)(1,2)
\psdots[dotscale=2,dotstyle=o](1,1)
\psline[linewidth=1.5pt](0,0)(2,0)(2,1)(1,2)(0,0)
\end{pspicture}}
\end{minipage}
\hfill
\begin{minipage}[c]{0.85\linewidth}
\begin{lem}\label{lem:P5}
Let $\A_6=\begin{pmatrix} -1 & 0 & 1 & 0 & 1 & 0 \\ -1 & -1 & -1 & 0 & 0 & 1 \\ 1 & 1 & 1 & 1 & 1 & 1 \end{pmatrix}$.
Then there are no $\bbeta$ such that $\HAi{6}$ has irreducible algebraic solutions.  
\end{lem}
\end{minipage}

\begin{proof}
Note that $\A_2$ is a subset of $\A_6$, as is also shown in Figure~\ref{fig:polygons_1interior_point}.
It follows from Corollary~\ref{cor:extension_less_algebraic_functions_resonant} that all $\bbeta$ for which the number of apexpoints is maximal for all conjugates must also be listed in Table~\ref{tab:solutionsP1}.
The interlacing condition for $\A_6$ is given by 
\begin{equation*}
(\entier{-\beta_1-\beta_2+\beta_3}, \entier{-\beta_1+\beta_3}, \entier{\beta_2+\beta_3}, \entier{2\beta_1-\beta_2+\beta_3}) \in \{(-1,0,0,1), (-1,-1,1,0)\}.
\end{equation*}
The families for $\A_2$ never satisfy the interlacing condition, except for $(\half,r,\half)$, which satisfies this condition if $0<r<\half$.
However, not all conjugates satisfy this condition.
Furthermore, one easily checks that none of the other 48 tuples for $\A_2$ gives $\bbeta$ such that all conjugates satisfy the interlacing condition for $\A_6$.
\end{proof}

\noindent 
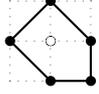
\begin{minipage}[c]{0.1\linewidth}
\psscalebox{0.53}{
\begin{pspicture}(2,2)
\psgrid[subgriddiv=1,gridlabels=0,gridwidth=0.2pt,griddots=5](0,0)(2,2)
\psdots[dotscale=2](1,0)(2,0)(0,1)(2,1)(1,2)
\psdots[dotscale=2,dotstyle=o](1,1)
\psline[linewidth=1.5pt](1,0)(2,0)(2,1)(1,2)(0,1)(1,0)
\end{pspicture}}
\end{minipage}
\hfill
\begin{minipage}[c]{0.85\linewidth}
\begin{lem}\label{lem:P6}
Let $\A_7=\begin{pmatrix} 0 & 1 & -1 & 0 & 1 & 0 \\ -1 & -1 & 0 & 0 & 0 & 1 \\ 1 & 1 & 1 & 1 & 1 & 1 \end{pmatrix}$.
Then $\HAi{7}$ has irreducible algebraic solutions if and only if, up to equivalence modulo $\Z$, we have $\bbeta = \pm (\frac{1}{3}, \frac{1}{3}, \half)$.
\end{lem}
\end{minipage}

\begin{proof}
The proof is similar to the proof of Lemma~\ref{lem:P5}.
In this case, $f(\A_4)$ is a subset of $\A_7$, with $f(x,y,z)=(y,x,z)$.
Hence we only have to check the interlacing condition for $\bbeta$ such that $H_{\A_4}(f^{-1}(\bbeta))$ is irreducible and has algebraic solutions.
The interlacing condition is 
\begin{multline*}
(\entier{-\beta_1-\beta_2+\beta_3}, \entier{-\beta_1+\beta_3}, \entier{\beta_2+\beta_3}, \entier{\beta_1-\beta_2+\beta_3}, \entier{\beta_1+\beta_2+\beta_3}) \in \\
\{(-1,-1,1,0,1), (-1,0,0,0,1)\}.
\end{multline*}
One easily checks that $(\half,r,\half)$ only satisfies this condition for $0<r<\half$, so there are no $r$ such that all conjugates satisfy the condition.
Furthermore, the only $\bbeta$ such that $\HAi{4}$ is irreducible and has algebraic solutions that satisfies the interlacing condition for $\A_7$ is $\bbeta=\pm(\inv{3}, \inv{3}, \half)$.
\end{proof}

\noindent 
\begin{minipage}[c]{0.1\linewidth}
\psscalebox{0.53}{
\begin{pspicture}(2,2)
\psgrid[subgriddiv=1,gridlabels=0,gridwidth=0.2pt,griddots=5](0,0)(2,2)
\psdots[dotscale=2](0,0)(1,0)(2,0)(0,1)(2,1)(1,2)
\psdots[dotscale=2,dotstyle=o](1,1)
\psline[linewidth=1.5pt](0,0)(2,0)(2,1)(1,2)(0,1)(0,0)
\end{pspicture}}
\end{minipage}
\hfill
\begin{minipage}[c]{0.85\linewidth}
\begin{lem}\label{lem:P9}
Let $\A_8=\begin{pmatrix} -1 & 0 & 1 & -1 & 0 & 1 & 0 \\ -1 & -1 & -1 & 0 & 0 & 0 & 1 \\ 1 & 1 & 1 & 1 & 1 & 1 & 1 \end{pmatrix}$.
Then there are no $\bbeta$ such that $\HAi{8}$ has irreducible algebraic solutions.
\end{lem}
\end{minipage}

\begin{proof}
This follows immediately from Corollary~\ref{cor:extension_less_algebraic_functions_resonant}, using the inclusion $\A_6 \subseteq \A_8$.
\end{proof}

\noindent 
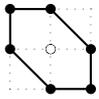
\begin{minipage}[c]{0.1\linewidth}
\psscalebox{0.53}{
\begin{pspicture}(2,2)
\psgrid[subgriddiv=1,gridlabels=0,gridwidth=0.2pt,griddots=5](0,0)(2,2)
\psdots[dotscale=2](1,0)(2,0)(0,1)(2,1)(0,2)(1,2)
\psdots[dotscale=2,dotstyle=o](1,1)
\psline[linewidth=1.5pt](1,0)(2,0)(2,1)(1,2)(0,2)(0,1)(1,0)
\end{pspicture}}
\end{minipage}
\hfill
\begin{minipage}[c]{0.85\linewidth}
\begin{lem}\label{lem:P10}
Let $\A_{9}=\begin{pmatrix} 0 & 1 & -1 & 0 & 1 & -1 & 0 \\ -1 & -1 & 0 & 0 & 0 & 1 & 1 \\ 1 & 1 & 1 & 1 & 1 & 1 & 1 \end{pmatrix}$.
Then $\HAi{9}$ has irreducible algebraic solutions if and only if, up to conjugation and equivalence modulo $\Z$, we have $\bbeta = \pm (\frac{1}{3}, \frac{1}{3}, \half)$.  
\end{lem}
\end{minipage}

\begin{proof}
Since $\A_7 \subseteq \A_{9}$, it suffices to check the interlacing condition for all $\bbeta$ such that $\sigma_{A_7}(k \bbeta)$ is maximal for all $k$, i.e., $\bbeta=\pm (\frac{1}{3}, \frac{1}{3}, \half)$.
The interlacing condition is 
\begin{multline*}
(\entier{-\beta_1-\beta_2+\beta_3}, \entier{-\beta_1+\beta_3}, \entier{-\beta_2+\beta_3}, \entier{\beta_2+\beta_3}, \entier{\beta_1+\beta_3}, \entier{\beta_1+\beta_2+\beta_3}) \in \\
\{(-1,-1,-1,1,1,1),(-1,0,0,0,0,1)\}.
\end{multline*}
It is easy to see that $\pm (\frac{1}{3}, \frac{1}{3}, \half)$ satisfies this condition and is non-resonant for $\A_{9}$.
\end{proof}

\subsection*{Polygons with exactly two interior points}\label{subsec:polygons-012_interior_points-2_interior_points}

Polygons with two interior points have 3 up to 10 boundary points.
This implies that their area lies between 5 and 12.
Furthermore, up to isomorphism each polygon lies in a square whose sides have length 24.
By shifting the polygon, we can assume that the lower left corner of the square is the origin.
By applying a translation, we can assume that each polygon has a vertex on each of the lower and left sides of the square.
This makes it feasible to compute all these polygons.

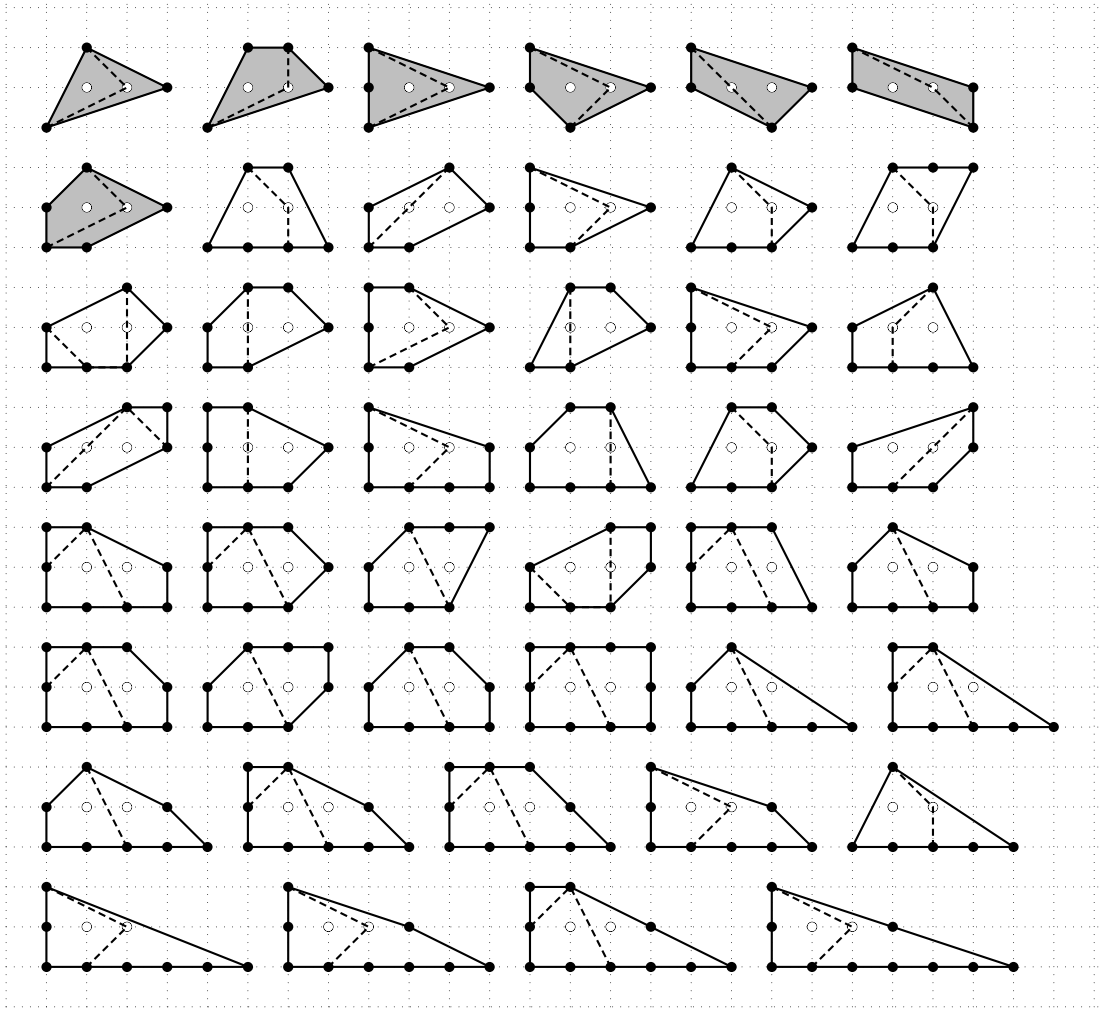
\begin{figure}
\centering
\psscalebox{0.53}{
\begin{pspicture}(27,25)
\psgrid[subgriddiv=1,gridlabels=0,gridwidth=0.2pt,griddots=5](27,25)

\pspolygon[fillstyle=solid,linewidth=0pt,fillcolor=lightgray](1,22)(4,23)(2,24)
\psdots[dotscale=2](1,22)(4,23)(2,24)
\psdots[dotscale=2,dotstyle=o](2,23)(3,23)
\psline[linewidth=1.5pt](1,22)(4,23)(2,24)(1,22)
\psline[linewidth=1.5pt,linestyle=dashed](1,22)(3,23)(2,24)

\pspolygon[fillstyle=solid,linewidth=0pt,fillcolor=lightgray](5,22)(8,23)(7,24)(6,24)
\psdots[dotscale=2](5,22)(8,23)(6,24)(7,24)
\psdots[dotscale=2,dotstyle=o](6,23)(7,23)
\psline[linewidth=1.5pt](5,22)(8,23)(7,24)(6,24)(5,22)
\psline[linewidth=1.5pt,linestyle=dashed](5,22)(7,23)(7,24)

\pspolygon[fillstyle=solid,linewidth=0pt,fillcolor=lightgray](9,22)(12,23)(9,24)
\psdots[dotscale=2](9,22)(9,23)(12,23)(9,24)
\psdots[dotscale=2,dotstyle=o](10,23)(11,23)
\psline[linewidth=1.5pt](9,22)(12,23)(9,24)(9,22)
\psline[linewidth=1.5pt,linestyle=dashed](9,22)(11,23)(9,24)

\pspolygon[fillstyle=solid,linewidth=0pt,fillcolor=lightgray](14,22)(16,23)(13,24)(13,23)
\psdots[dotscale=2](14,22)(13,23)(16,23)(13,24)
\psdots[dotscale=2,dotstyle=o](14,23)(15,23)
\psline[linewidth=1.5pt](14,22)(16,23)(13,24)(13,23)(14,22)
\psline[linewidth=1.5pt,linestyle=dashed](14,22)(15,23)(13,24)

\pspolygon[fillstyle=solid,linewidth=0pt,fillcolor=lightgray](19,22)(20,23)(17,24)(17,23)
\psdots[dotscale=2](19,22)(17,23)(20,23)(17,24)
\psdots[dotscale=2,dotstyle=o](18,23)(19,23)
\psline[linewidth=1.5pt](19,22)(20,23)(17,24)(17,23)(19,22)
\psline[linewidth=1.5pt,linestyle=dashed](19,22)(17,24)

\pspolygon[fillstyle=solid,linewidth=0pt,fillcolor=lightgray](24,22)(24,23)(21,24)(21,23)
\psdots[dotscale=2](24,22)(21,23)(24,23)(21,24)
\psdots[dotscale=2,dotstyle=o](22,23)(23,23)
\psline[linewidth=1.5pt](24,22)(24,23)(21,24)(21,23)(24,22)
\psline[linewidth=1.5pt,linestyle=dashed](24,22)(23,23)(21,24)

\pspolygon[fillstyle=solid,linewidth=0pt,fillcolor=lightgray](1,19)(2,19)(4,20)(2,21)(1,20)
\psdots[dotscale=2](1,19)(2,19)(1,20)(4,20)(2,21)
\psdots[dotscale=2,dotstyle=o](2,20)(3,20)
\psline[linewidth=1.5pt](1,19)(2,19)(4,20)(2,21)(1,20)(1,19)
\psline[linewidth=1.5pt,linestyle=dashed](1,19)(3,20)(2,21)

\psdots[dotscale=2](5,19)(6,19)(7,19)(8,19)(6,21)(7,21)
\psdots[dotscale=2,dotstyle=o](6,20)(7,20)
\psline[linewidth=1.5pt](5,19)(8,19)(7,21)(6,21)(5,19)
\psline[linewidth=1.5pt,linestyle=dashed](7,19)(7,20)(6,21)

\psdots[dotscale=2](9,19)(10,19)(9,20)(12,20)(11,21)
\psdots[dotscale=2,dotstyle=o](10,20)(11,20)
\psline[linewidth=1.5pt](9,19)(10,19)(12,20)(11,21)(9,20)(9,19)
\psline[linewidth=1.5pt,linestyle=dashed](9,19)(11,21)

\psdots[dotscale=2](13,19)(14,19)(13,20)(16,20)(13,21)
\psdots[dotscale=2,dotstyle=o](14,20)(15,20)
\psline[linewidth=1.5pt](13,19)(14,19)(16,20)(13,21)(13,19)
\psline[linewidth=1.5pt,linestyle=dashed](14,19)(15,20)(13,21)

\psdots[dotscale=2](17,19)(18,19)(19,19)(20,20)(18,21)
\psdots[dotscale=2,dotstyle=o](18,20)(19,20)
\psline[linewidth=1.5pt](17,19)(19,19)(20,20)(18,21)(17,19)
\psline[linewidth=1.5pt,linestyle=dashed](19,19)(19,20)(18,21)

\psdots[dotscale=2](21,19)(22,19)(23,19)(22,21)(23,21)(24,21)
\psdots[dotscale=2,dotstyle=o](22,20)(23,20)
\psline[linewidth=1.5pt](21,19)(23,19)(24,21)(22,21)(21,19)
\psline[linewidth=1.5pt,linestyle=dashed](23,19)(23,20)(22,21)

\psdots[dotscale=2](1,16)(2,16)(3,16)(1,17)(4,17)(3,18)
\psdots[dotscale=2,dotstyle=o](2,17)(3,17)
\psline[linewidth=1.5pt](1,16)(3,16)(4,17)(3,18)(1,17)(1,16)
\psline[linewidth=1.5pt,linestyle=dashed](1,17)(2,16)(3,16)(3,18)

\psdots[dotscale=2](5,16)(6,16)(5,17)(8,17)(6,18)(7,18)
\psdots[dotscale=2,dotstyle=o](6,17)(7,17)
\psline[linewidth=1.5pt](5,16)(6,16)(8,17)(7,18)(6,18)(5,17)(5,16)
\psline[linewidth=1.5pt,linestyle=dashed](6,16)(6,18)

\psdots[dotscale=2](9,16)(10,16)(9,17)(12,17)(9,18)(10,18)
\psdots[dotscale=2,dotstyle=o](10,17)(11,17)
\psline[linewidth=1.5pt](9,16)(10,16)(12,17)(10,18)(9,18)(9,16)
\psline[linewidth=1.5pt,linestyle=dashed](9,16)(11,17)(10,18)

\psdots[dotscale=2](13,16)(14,16)(16,17)(14,18)(15,18)
\psdots[dotscale=2,dotstyle=o](14,17)(15,17)
\psline[linewidth=1.5pt](13,16)(14,16)(16,17)(15,18)(14,18)(13,16)
\psline[linewidth=1.5pt,linestyle=dashed](14,16)(14,18)

\psdots[dotscale=2](17,16)(18,16)(19,16)(17,17)(20,17)(17,18)
\psdots[dotscale=2,dotstyle=o](18,17)(19,17)
\psline[linewidth=1.5pt](17,16)(19,16)(20,17)(17,18)(17,16)
\psline[linewidth=1.5pt,linestyle=dashed](18,16)(19,17)(17,18)

\psdots[dotscale=2](21,16)(22,16)(23,16)(24,16)(21,17)(23,18)
\psdots[dotscale=2,dotstyle=o](22,17)(23,17)
\psline[linewidth=1.5pt](21,16)(24,16)(23,18)(21,17)(21,16)
\psline[linewidth=1.5pt,linestyle=dashed](22,16)(22,17)(23,18)

\psdots[dotscale=2](1,13)(2,13)(1,14)(4,14)(3,15)(4,15)
\psdots[dotscale=2,dotstyle=o](2,14)(3,14)
\psline[linewidth=1.5pt](1,13)(2,13)(4,14)(4,15)(3,15)(1,14)(1,13)
\psline[linewidth=1.5pt,linestyle=dashed](1,13)(3,15)(4,14)

\psdots[dotscale=2](5,13)(6,13)(7,13)(5,14)(8,14)(5,15)(6,15)
\psdots[dotscale=2,dotstyle=o](6,14)(7,14)
\psline[linewidth=1.5pt](5,13)(7,13)(8,14)(6,15)(5,15)(5,13)
\psline[linewidth=1.5pt,linestyle=dashed](6,13)(6,15)

\psdots[dotscale=2](9,13)(10,13)(11,13)(12,13)(9,14)(12,14)(9,15)
\psdots[dotscale=2,dotstyle=o](10,14)(11,14)
\psline[linewidth=1.5pt](9,13)(12,13)(12,14)(9,15)(9,13)
\psline[linewidth=1.5pt,linestyle=dashed](9,15)(11,14)(10,13)

\psdots[dotscale=2](13,13)(14,13)(15,13)(16,13)(13,14)(14,15)(15,15)
\psdots[dotscale=2,dotstyle=o](14,14)(15,14)
\psline[linewidth=1.5pt](13,13)(16,13)(15,15)(14,15)(13,14)(13,13)
\psline[linewidth=1.5pt,linestyle=dashed](15,15)(15,13)

\psdots[dotscale=2](17,13)(18,13)(19,13)(20,14)(18,15)(19,15)
\psdots[dotscale=2,dotstyle=o](18,14)(19,14)
\psline[linewidth=1.5pt](17,13)(19,13)(20,14)(19,15)(18,15)(17,13)
\psline[linewidth=1.5pt,linestyle=dashed](19,13)(19,14)(18,15)

\psdots[dotscale=2](21,13)(22,13)(23,13)(21,14)(24,14)(24,15)
\psdots[dotscale=2,dotstyle=o](22,14)(23,14)
\psline[linewidth=1.5pt](21,13)(23,13)(24,14)(24,15)(21,14)(21,13)
\psline[linewidth=1.5pt,linestyle=dashed](22,13)(24,15)

\psdots[dotscale=2](1,10)(2,10)(3,10)(4,10)(1,11)(4,11)(1,12)(2,12)
\psdots[dotscale=2,dotstyle=o](2,11)(3,11)
\psline[linewidth=1.5pt](1,10)(4,10)(4,11)(2,12)(1,12)(1,10)
\psline[linewidth=1.5pt,linestyle=dashed](1,11)(2,12)(3,10)

\psdots[dotscale=2](5,10)(6,10)(7,10)(5,11)(8,11)(5,12)(6,12)(7,12)
\psdots[dotscale=2,dotstyle=o](6,11)(7,11)
\psline[linewidth=1.5pt](5,10)(7,10)(8,11)(7,12)(5,12)(5,10)
\psline[linewidth=1.5pt,linestyle=dashed](5,11)(6,12)(7,10)

\psdots[dotscale=2](9,10)(10,10)(11,10)(9,11)(10,12)(11,12)(12,12)
\psdots[dotscale=2,dotstyle=o](10,11)(11,11)
\psline[linewidth=1.5pt](9,10)(11,10)(12,12)(10,12)(9,11)(9,10)
\psline[linewidth=1.5pt,linestyle=dashed](11,10)(10,12)

\psdots[dotscale=2](13,10)(14,10)(15,10)(13,11)(16,11)(15,12)(16,12)
\psdots[dotscale=2,dotstyle=o](14,11)(15,11)
\psline[linewidth=1.5pt](13,10)(15,10)(16,11)(16,12)(15,12)(13,11)(13,10)
\psline[linewidth=1.5pt,linestyle=dashed](13,11)(14,10)(15,10)(15,12)

\psdots[dotscale=2](17,10)(18,10)(19,10)(20,10)(17,11)(17,12)(18,12)(19,12)
\psdots[dotscale=2,dotstyle=o](18,11)(19,11)
\psline[linewidth=1.5pt](17,10)(20,10)(19,12)(17,12)(17,10)
\psline[linewidth=1.5pt,linestyle=dashed](17,11)(18,12)(19,10)

\psdots[dotscale=2](21,10)(22,10)(23,10)(24,10)(21,11)(24,11)(22,12)
\psdots[dotscale=2,dotstyle=o](22,11)(23,11)
\psline[linewidth=1.5pt](21,10)(24,10)(24,11)(22,12)(21,11)(21,10)
\psline[linewidth=1.5pt,linestyle=dashed](22,12)(23,10)

\psdots[dotscale=2](1,7)(2,7)(3,7)(4,7)(1,8)(4,8)(1,9)(2,9)(3,9)
\psdots[dotscale=2,dotstyle=o](2,8)(3,8)
\psline[linewidth=1.5pt](1,7)(4,7)(4,8)(3,9)(1,9)(1,7)
\psline[linewidth=1.5pt,linestyle=dashed](1,8)(2,9)(3,7)

\psdots[dotscale=2](5,7)(6,7)(7,7)(5,8)(8,8)(6,9)(7,9)(8,9)
\psdots[dotscale=2,dotstyle=o](6,8)(7,8)
\psline[linewidth=1.5pt](5,7)(7,7)(8,8)(8,9)(6,9)(5,8)(5,7)
\psline[linewidth=1.5pt,linestyle=dashed](7,7)(6,9)

\psdots[dotscale=2](9,7)(10,7)(11,7)(12,7)(9,8)(12,8)(10,9)(11,9)
\psdots[dotscale=2,dotstyle=o](10,8)(11,8)
\psline[linewidth=1.5pt](9,7)(12,7)(12,8)(11,9)(10,9)(9,8)(9,7)
\psline[linewidth=1.5pt,linestyle=dashed](11,7)(10,9)

\psdots[dotscale=2](13,7)(14,7)(15,7)(16,7)(13,8)(16,8)(13,9)(14,9)(15,9)(16,9)
\psdots[dotscale=2,dotstyle=o](14,8)(15,8)
\psline[linewidth=1.5pt](13,7)(16,7)(16,9)(13,9)(13,7)
\psline[linewidth=1.5pt,linestyle=dashed](13,8)(14,9)(15,7)

\psdots[dotscale=2](17,7)(18,7)(19,7)(20,7)(21,7)(17,8)(18,9)
\psdots[dotscale=2,dotstyle=o](18,8)(19,8)
\psline[linewidth=1.5pt](17,7)(21,7)(18,9)(17,8)(17,7)
\psline[linewidth=1.5pt,linestyle=dashed](18,9)(19,7)

\psdots[dotscale=2](22,7)(23,7)(24,7)(25,7)(26,7)(22,8)(22,9)(23,9)
\psdots[dotscale=2,dotstyle=o](23,8)(24,8)
\psline[linewidth=1.5pt](22,7)(26,7)(23,9)(22,9)(22,7)
\psline[linewidth=1.5pt,linestyle=dashed](22,8)(23,9)(24,7)

\psdots[dotscale=2](1,4)(2,4)(3,4)(4,4)(5,4)(1,5)(4,5)(2,6)
\psdots[dotscale=2,dotstyle=o](2,5)(3,5)
\psline[linewidth=1.5pt](1,4)(5,4)(4,5)(2,6)(1,5)(1,5)(1,4)
\psline[linewidth=1.5pt,linestyle=dashed](2,6)(3,4)

\psdots[dotscale=2](6,4)(7,4)(8,4)(9,4)(10,4)(6,5)(9,5)(6,6)(7,6)
\psdots[dotscale=2,dotstyle=o](7,5)(8,5)
\psline[linewidth=1.5pt](6,4)(10,4)(9,5)(7,6)(6,6)(6,4)
\psline[linewidth=1.5pt,linestyle=dashed](6,5)(7,6)(8,4)

\psdots[dotscale=2](11,4)(12,4)(13,4)(14,4)(15,4)(11,5)(14,5)(11,6)(12,6)(13,6)
\psdots[dotscale=2,dotstyle=o](12,5)(13,5)
\psline[linewidth=1.5pt](11,4)(15,4)(13,6)(11,6)(11,4)
\psline[linewidth=1.5pt,linestyle=dashed](11,5)(12,6)(13,4)

\psdots[dotscale=2](16,4)(17,4)(18,4)(19,4)(20,4)(16,5)(19,5)(16,6)
\psdots[dotscale=2,dotstyle=o](17,5)(18,5)
\psline[linewidth=1.5pt](16,4)(20,4)(19,5)(16,6)(16,4)
\psline[linewidth=1.5pt,linestyle=dashed](16,6)(18,5)(17,4)

\psdots[dotscale=2](21,4)(22,4)(23,4)(24,4)(25,4)(22,6)
\psdots[dotscale=2,dotstyle=o](22,5)(23,5)
\psline[linewidth=1.5pt](21,4)(25,4)(22,6)(21,4)
\psline[linewidth=1.5pt,linestyle=dashed](23,4)(23,5)(22,6)

\psdots[dotscale=2](1,1)(2,1)(3,1)(4,1)(5,1)(6,1)(1,2)(1,3)
\psdots[dotscale=2,dotstyle=o](2,2)(3,2)
\psline[linewidth=1.5pt](1,1)(6,1)(1,3)(1,1)
\psline[linewidth=1.5pt,linestyle=dashed](1,3)(3,2)(2,1)

\psdots[dotscale=2](7,1)(8,1)(9,1)(10,1)(11,1)(12,1)(7,2)(10,2)(7,3)
\psdots[dotscale=2,dotstyle=o](8,2)(9,2)
\psline[linewidth=1.5pt](7,1)(12,1)(10,2)(7,3)(7,1)
\psline[linewidth=1.5pt,linestyle=dashed](7,3)(9,2)(8,1)

\psdots[dotscale=2](13,1)(14,1)(15,1)(16,1)(17,1)(18,1)(13,2)(16,2)(13,3)(14,3)
\psdots[dotscale=2,dotstyle=o](14,2)(15,2)
\psline[linewidth=1.5pt](13,1)(18,1)(14,3)(13,3)(13,1)
\psline[linewidth=1.5pt,linestyle=dashed](13,2)(14,3)(15,1)

\psdots[dotscale=2](19,1)(20,1)(21,1)(22,1)(23,1)(24,1)(25,1)(19,2)(22,2)(19,3)
\psdots[dotscale=2,dotstyle=o](20,2)(21,2)
\psline[linewidth=1.5pt](19,1)(25,1)(19,3)(19,1)
\psline[linewidth=1.5pt,linestyle=dashed](19,3)(21,2)(20,1)
\end{pspicture}}
\caption{The polygons with exactly 2 interior points
\label{fig:polygons_2interior_points}}
\end{figure}

After computing all polygons in this square, we compute the different isomorphism classes.
Note that the number of vertices, the number of interior points and the discrete lengths of the sides are invariant under isomorphisms.
By the \emph{discrete length} of an edge we mean the number of lattice points on this edge minus 1.
Hence for each pair of polygons, we can first check whether these invariants are the same.
If they do, we check whether there exists an isomorphism.
This can be done by computing the functions that map sets of three lattice points of the first polygon to a fixed set of three lattice points of the second polygon.
It turns out that there are 45 isomorphism classes of polygons with two interior points.
They are shown in Figure~\ref{fig:polygons_2interior_points}.

As in the previous section, we will compute all parameters such that the corresponding functions are algebraic.
For the 38 non-shaded polygons, a set isomorphic to $\A_6$ is included in $\A$.
These inclusions are also drawn in Figure~\ref{fig:polygons_2interior_points}.
In all but three cases, this subset of $\A$ is a translation or reflection of $\A_6$.
For the remaining three polygons, it might not be immediately clear that the indicated subset is indeed isomorphic to $\A_6$.
In these cases, note that $\A_6$ is the only convex lattice polygon with one interior point, area 5 and an edge of length 2.
This leaves us with the seven shaded polygons.
Five of these are elements of one of the two families of polygons with algebraic functions, which will be discussed in Section~\ref{sec:polygons-34_boundary_points}.
The second and the seventh shaded polygon aren't elements of such families.
For the second polygon, the inclusion of a set isomorphic to $\A_4$ gives us a proof of the following lemma, similar to the proof of Lemma~\ref{lem:P5}: \\

\noindent 
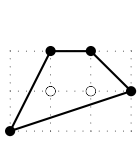
\begin{minipage}[c]{0.13\linewidth}
\psscalebox{0.53}{
\begin{pspicture}(3,2)
\psgrid[subgriddiv=1,gridlabels=0,gridwidth=0.2pt,griddots=5](0,0)(3,2)
\psdots[dotscale=2](0,0)(3,1)(1,2)(2,2)
\psdots[dotscale=2,dotstyle=o](1,1)(2,1)
\psline[linewidth=1.5pt](0,0)(3,1)(2,2)(1,2)(0,0)
\end{pspicture}}
\end{minipage}
\hfill
\begin{minipage}[c]{0.82\linewidth}
\begin{lem}\label{lem:Q62}
Let $\A_{10}=\begin{pmatrix} 0 & 1 & 2 & 3 & 1 & 2 \\ 0 & 1 & 1 & 1 & 2 & 2 \\ 1 & 1 & 1 & 1 & 1 & 1 \end{pmatrix}$.
Then $\HAi{10}$ has irreducible algebraic solutions if and only if, up to equivalence modulo $\Z$, we have $\bbeta = \pm (\frac{1}{3}, \frac{5}{6}, \half)$. \\
\end{lem}
\end{minipage}
\par
Similarly, again using a inclusion of $\A_4$, one can show that there are no irreducible algebraic functions for the seventh shaded polygon.


\section{Polygons with at least 3 interior points and at least 5 boundary points}\label{sec:polygons-many_boundary_points}

In this section, we will show that the hypergeometric functions associated to polygons with at least 3 interior points and at least 5 boundary points are never irreducible and algebraic.
We do this by proving that such polygons contain subpolygons with 1 or 2 boundary points without irreducible algebraic functions.
Polygons with exactly 3 or 4 boundary points will be treated in the next section.

\begin{defn}\label{defn:type_convex_lattice_polygon}
We say that $P$ has \emph{type $(i,b)$} if $P$ has $i$ interior points and $b$ boundary points.
We write $t(P)=(i,b)$.
If $P=P(\A)$, then we also call $\A$ of \emph{type $(i,b)$} and write $t(\A)=(i,b)$.
\end{defn}

\begin{defn}\label{defn:types_sets}
Let $T=\{(i,b) \in \N^2 \ | \ i=1, b \geq 7 \tr{ or } i \geq 2, b \geq 5\}$ and $S=\{(i,b) \in T \ | \ i=1 \tr{ or } i=2\}$.
Let $\prec$ be the lexicographical ordering on $T$ and $S$.
\end{defn}

\begin{rmk}\label{rmk:vol7_implies_T}
Note that $t(P) \in T$ if and only if $P$ has at least one interior point, at least 5 boundary points and normalized area at least 7.
If $P$ and $P'$ are convex lattice polygons with $P' \subsetneq P$, then $t(P') \prec t(P)$.
\end{rmk}

\begin{lem}\label{lem:reductionTS}
Let $P$ be a convex lattice polygon with $t(P) \in T \setminus S$.
Then there exists a convex lattice polygon $P' \subsetneq P$ of type $t(P') \in T$.
\end{lem}

\begin{cor}\label{cor:reduction_polygons}
Let $P$ be a convex lattice polygon with $t(P) \in T$.
Then there exists a convex lattice polygon $P' \subs P$ of type $t(P') \in S$.
\end{cor}

To prove Lemma~\ref{lem:reductionTS}, we consider 6 cases:
polygons of type (3,5);
triangles, quadrilaterals and pentagons with at least 4 interior points and exactly 5 boundary points;
triangles with at least 3 interior points and at least 6 boundary points and 2 edges of length 1;
and other polygons with at least 3 interior points and at least 6 boundary points.

\begin{lem}\label{lem:reduction35}
If $P$ is a lattice polygon with exactly 3 interior points and 5 boundary points, then there is a subpolygon $P' \subs P$ with 2 interior points and 5 boundary points.
\end{lem}

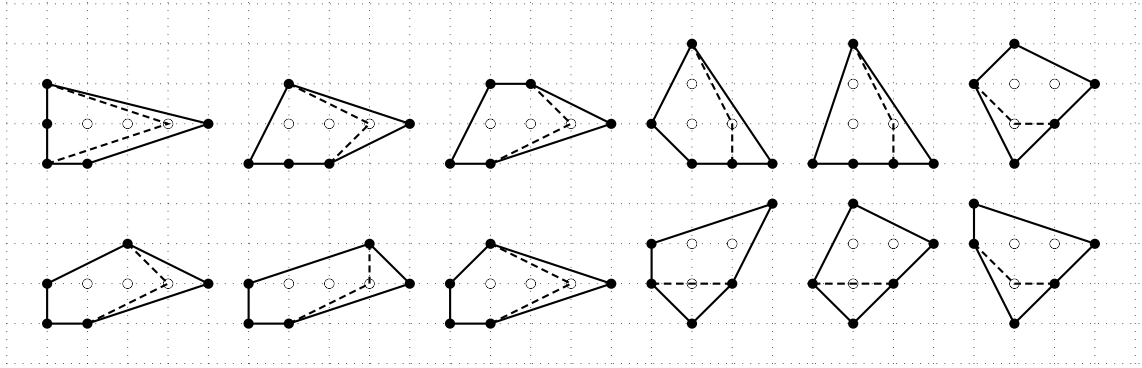
\begin{figure}
\centering
\psscalebox{0.53}{
\begin{pspicture}(28,9)
\psgrid[subgriddiv=1,gridlabels=0,gridwidth=0.2pt,griddots=5](28,9)

\psdots[dotscale=2](1,5)(2,5)(1,6)(5,6)(1,7) 
\psdots[dotscale=2,dotstyle=o](2,6)(3,6)(4,6)
\psline[linewidth=1.5pt](1,5)(2,5)(5,6)(1,7)(1,5)
\psline[linewidth=1.5pt,linestyle=dashed](1,5)(4,6)(1,7)

\psdots[dotscale=2](6,5)(7,5)(8,5)(10,6)(7,7)
\psdots[dotscale=2,dotstyle=o](7,6)(8,6)(9,6)
\psline[linewidth=1.5pt](6,5)(8,5)(10,6)(7,7)(6,5)
\psline[linewidth=1.5pt,linestyle=dashed](8,5)(9,6)(7,7)

\psdots[dotscale=2](11,5)(12,5)(15,6)(12,7)(13,7) 
\psdots[dotscale=2,dotstyle=o](12,6)(13,6)(14,6)
\psline[linewidth=1.5pt](11,5)(12,5)(15,6)(13,7)(12,7)(11,5)
\psline[linewidth=1.5pt,linestyle=dashed](12,5)(14,6)(13,7)

\psdots[dotscale=2](17,5)(18,5)(19,5)(16,6)(17,8) 
\psdots[dotscale=2,dotstyle=o](17,6)(18,6)(17,7)
\psline[linewidth=1.5pt](17,5)(19,5)(17,8)(16,6)(17,5)
\psline[linewidth=1.5pt,linestyle=dashed](18,5)(18,6)(17,8)

\psdots[dotscale=2](20,5)(21,5)(22,5)(23,5)(21,8) 
\psdots[dotscale=2,dotstyle=o](21,6)(22,6)(21,7)
\psline[linewidth=1.5pt](20,5)(23,5)(21,8)(20,5)
\psline[linewidth=1.5pt,linestyle=dashed](22,5)(22,6)(21,8)

\psdots[dotscale=2](25,5)(26,6)(24,7)(27,7)(25,8) 
\psdots[dotscale=2,dotstyle=o](25,6)(25,7)(26,7)
\psline[linewidth=1.5pt](25,5)(27,7)(25,8)(24,7)(25,5)
\psline[linewidth=1.5pt,linestyle=dashed](24,7)(25,6)(26,6)

\psdots[dotscale=2](1,1)(2,1)(1,2)(5,2)(3,3) 
\psdots[dotscale=2,dotstyle=o](2,2)(3,2)(4,2)
\psline[linewidth=1.5pt](1,1)(2,1)(5,2)(3,3)(1,2)(1,1)
\psline[linewidth=1.5pt,linestyle=dashed](2,1)(4,2)(3,3)

\psdots[dotscale=2](6,1)(7,1)(6,2)(10,2)(9,3) 
\psdots[dotscale=2,dotstyle=o](7,2)(8,2)(9,2)
\psline[linewidth=1.5pt](6,1)(7,1)(10,2)(9,3)(6,2)(6,1)
\psline[linewidth=1.5pt,linestyle=dashed](7,1)(9,2)(9,3)

\psdots[dotscale=2](11,1)(12,1)(11,2)(15,2)(12,3) 
\psdots[dotscale=2,dotstyle=o](12,2)(13,2)(14,2)
\psline[linewidth=1.5pt](11,1)(12,1)(15,2)(12,3)(11,2)(11,1)
\psline[linewidth=1.5pt,linestyle=dashed](12,1)(14,2)(12,3)

\psdots[dotscale=2](17,1)(16,2)(18,2)(16,3)(19,4) 
\psdots[dotscale=2,dotstyle=o](17,2)(17,3)(18,3)
\psline[linewidth=1.5pt](17,1)(18,2)(19,4)(16,3)(16,2)(17,1)
\psline[linewidth=1.5pt,linestyle=dashed](16,2)(18,2)

\psdots[dotscale=2](21,1)(20,2)(22,2)(23,3)(21,4) 
\psdots[dotscale=2,dotstyle=o](21,2)(21,3)(22,3)
\psline[linewidth=1.5pt](21,1)(23,3)(21,4)(20,2)(21,1)
\psline[linewidth=1.5pt,linestyle=dashed](20,2)(22,2)

\psdots[dotscale=2](25,1)(26,2)(24,3)(27,3)(24,4) 
\psdots[dotscale=2,dotstyle=o](25,2)(25,3)(26,3)
\psline[linewidth=1.5pt](25,1)(27,3)(24,4)(24,3)(25,1)
\psline[linewidth=1.5pt,linestyle=dashed](24,3)(25,2)(26,2)
\end{pspicture}}
\caption{The lattice polygons with 3 interior points and 5 boundary points
\label{fig:polygons35}}
\end{figure}

\begin{proof}
There are 12 convex lattice polygons with 3 interior lattice points and 5 lattice points on the boundary.
They are shown in Figure~\ref{fig:polygons35}.
For each of these polygons, a subpolygon with 2 interior points and 5 boundary points is indicated. 
\end{proof}

\begin{lem}\label{lem:reduction_triangle5}
Suppose that $P$ is a lattice triangle of type $(i,5)$ with $i \geq 4$. 
Then there exists a lattice polygon $P' \subsetneq P$ with $t(P') \in T$.
\end{lem}

\begin{proof}
Suppose that $P$ doesn't have an edge of discrete length 3. 
Then the lengths of the edges are 1, 2 and 2.
However, if the lattice points on the boundary are $(0,0)$, $v_1$, $v_2=2v_1$, $v_3$ and $v_4=2v_3$, then $v_1+v_3$ is also a lattice point on the boundary.
Hence $P$ has an edge of discrete length 3, and the other two edges have length 1.
After a suitable transformation, the lattice points on the boundary are $(0,0)$, $(1,0)$, $(2,0)$, $(3,0)$ and $(c,d)$ for some $c, d \in \Z$ with $\gcd(c,d)=\gcd(c-3,d)=1$.
We can assume that $d>0$ and $c \geq 2$.
Furthermore, we can apply an isomorphism $(x,y) \mapsto (x+ny,y)$ with $n \in \Z$.
This maps the basis $(0,0)$, $(1,0)$, $(2,0)$, $(3,0)$ to itself, and it maps $(c,d)$ to $(c+nd,d)$.
It is possible to choose $n$ so that $0 \leq c+nd < d$.
Hence we can assume that $0 \leq c < d$.

Suppose that $c=2$.
Then $d$ is odd and the interior points are the points $(1,k)$ wih $1 \leq k \leq \frac{d-1}{2}$ and $(2,k)$ with $1 \leq k \leq d-1$ (see Figure~\ref{fig:triangle_c2}).
Let $P'$ be the triangle with vertices $(0,0)$, $(3,0)$ and $(2,d-1)$.
Then $P'$ has 6 boundary points and $i-2$ interior points, so $t(P') \in T$.

Note that $c \neq 3$, because $\gcd(c-3,d)=1$.
Suppose that $4 \leq c < d$.
Then $P$ is given by the inequalities $x_2 \geq 0$, $c x_2 \leq d x_1$ and $(c-3) x_2 \geq d (x_1-3)$, so $(3, \entier{\frac{3d}{c}})$ is an interior point (note that $c \nmid 3d)$ (see Figure~\ref{fig:triangle_cd}).
Let $P'$ be the triangle with vertices $(0,0)$, $(3,0)$ and $(3, \entier{\frac{3d}{c}})$.
Since $\frac{3d}{c} > 3$, the point $(2,1)$ is an interior point and the area of $P'$ is at least 9.
$P'$ has $4 + \entier{\frac{3d}{c}} \geq 5$ boundary points, so $t(P') \in T$ by Remark~\ref{rmk:vol7_implies_T}.
\end{proof}

\begin{lem}\label{lem:reduction_quadrilateral5}
Suppose that $P$ is a convex lattice quadrilateral of type $(i,5)$ with $i \geq 4$. 
Then there exists a convex lattice polygon $P' \subsetneq P$ with $t(P') \in T$.
\end{lem}

\begin{proof}
A quadrilateral with exactly 5 lattice points on the boundary must have edges of discrete length 1, 1, 1 and 2.
Let the vertices, in counterclockwise order, be $v_0, \ld, v_4$.
After a suitable transformation, we have $v_0=(-1,0)$, $v_1=(0,0)$ and $v_2=(1,0)$ and $v_{32}, v_{42} > 0$.
Suppose that the triangles with vertices $(\pm 1, 0)$, $v_3$ and $v_4$ both have area 1.
Then the triangles with vertices $(-1,0)$, $(1,0)$ and $v_i$ with $i=3, 4$ both have area $\vol(P)-1$ (see Figure~\ref{fig:quadrilateral_triangles}).
This implies that $v_{32} = v_{42} = \frac{\vol(P)-1}{2}$.
The area of $P$ is equal to $v_{32}+v_{42}+(v_{31} v_{42} - v_{41} v_{32})$ (by dissecting $P$ with the dotted lines as in Figure~\ref{fig:quadrilateral_triangles}).
Hence $(v_{31} - v_{41}) v_{32} = 1$, so $v_{32}=1$.
But then $\vol(P) = 2v_{32}+1 = 3$, contradicting the assumption that $P$ has at least 4 interior points.

Hence we can assume that the triangle with vertices $(-1,0)$, $v_3$ and $v_4$ has area at least 2.
This implies that there is an interior point of $P$ on or above the line from $(-1,0)$ to $v_3$.
Let $P'$ be the convex hull of all lattice points in or on the boundary of $P$, except for $v_4$.

If $P'$ has at least one interior point, then $P'$ satisfies the conditions of Remark~\ref{rmk:vol7_implies_T}:
$P'$ contains $i+4$ lattice points, so $\vol(P) \geq 2 \cdot 1 + (i+3) - 2 = i+3 \geq 7$, and $P'$ has at least 5 boundary points: $v_0$, $v_1$, $v_2$, $v_3$ and the interior point of $P$ on or above the line from $(-1,0)$ to $v_3$.

Suppose that $P'$ has no interior lattice points.
By~\cite{rabinowitz_census_lattice_polygons_one_point}, the only polygons without interior points are (up to isomorphism) triangles with vertices $(0,0), (p,0)$ and $(0,1)$, the triangles with vertices $(0,0), (2,0)$ and $(0,2)$, and trapezoids with vertices $(0,0),(p,0),(q,1)$ and $(0,1)$.
If $P'$ is a triangle, then all interior points of $P$ lie on the line from $(-1,0)$ to $v_3$.
Hence $P'$ has edges of discrete length 1, 2 and $i+1$.
It follows that $P'$ is a quadrilateral, and hence has two opposite edges of discrete length 1.
Let $v_5$ be the fourth vertex.
Then all interior points of $P$ must be on the line from $v_3$ to $v_5$, as in Figure~\ref{fig:quadrilateral_interior_points}.
Let $v_6$ be the vertex closest to $v_3$ and let $P''$ be the convex hull of $v_0$, $v_2$, $v_6$ and $v_4$.
Then $P''$ has $i-1$ interior points and 5 boundary points, so $t(P'') \in T$.
\end{proof}

\begin{figure}
\centering
\begin{minipage}[t]{0.3\linewidth}
\centering
\begin{subfigure}[t]{0.3\textwidth}
\centering
\psscaleboxto(\linewidth,0pt){
\begin{pspicture}(0,0)(3,5)
\psdots(0,0)(1,0)(2,0)(3,0)(1,1)(1,2)(2,1)(2,2)(2,3)(2,4)(2,5)
\psline(0,0)(3,0)(2,5)(0,0)
\psline[linestyle=dashed](0,0)(2,4)(3,0)
\end{pspicture}}
\caption{$c=2$
\label{fig:triangle_c2}}
\end{subfigure}
\hfill
\begin{subfigure}[t]{0.6\textwidth}
\centering
\psscaleboxto(\linewidth,0pt){
\begin{pspicture}(0,0)(5,7)
\psdots(0,0)(1,0)(2,0)(3,0)(1,1)(2,1)(3,1)(2,2)(3,2)(3,3)(3,4)(4,4)(4,5)(5,7)
\psline(0,0)(3,0)(5,7)(0,0)
\psline[linestyle=dashed](0,0)(3,0)(3,4)(0,0)
\end{pspicture}}
\caption{$4 \leq c<d$
\label{fig:triangle_cd}}
\end{subfigure}
\caption{Subpolygons of triangles
\label{fig:triangle}}
\end{minipage}
\hfill
\begin{minipage}[t]{0.65\linewidth}
\centering
\begin{subfigure}[t]{0.47\textwidth}
\centering
\psscaleboxto(\linewidth,0pt){
\begin{pspicture}(-3,-1)(3.5,4)
\psdots(-1,0)(0,0)(1,0)(2.5,2)(-2,3)
\psline(-1,0)(1,0)(2.5,2)(-2,3)(-1,0)
\psline[linestyle=dashed](-1,0)(2.5,2)
\psline[linestyle=dashed](1,0)(-2,3)
\psline[linestyle=dotted](0,0)(2.5,2)
\psline[linestyle=dotted](0,0)(-2,3)
\rput(-1.5,-0.37){$v_0=(-1,0)$}
\rput(0,-0.4){$v_1$}
\rput(1.3,-0.37){$v_2=(1,0)$}
\rput(2.9,2){$v_3$}
\rput(-2.4,3.2){$v_4$}
\end{pspicture}}
\caption{Quadrilateral
\label{fig:quadrilateral_triangles}}
\end{subfigure}
\hfill
\begin{subfigure}[t]{0.47\textwidth}
\centering
\psscaleboxto(\linewidth,0pt){
\begin{pspicture}(-3,-1)(3.5,4)
\psdots(-1,0)(0,0)(1,0)(2.5,2)(-2,3)
\psline(-1,0)(1,0)(2.5,2)(-2,3)(-1,0)
\psline[linestyle=dotted](-1,0)(2.5,2)
\psdots(-0.5,1.5)
\psline[linestyle=dashed](-1,0)(-0.5,1.5)(2.5,2)
\psdots(0.1,1.6)(0.7,1.7)(1.3,1.8)(1.9,1.9)
\psline[linestyle=dashed](1,0)(1.9,1.9)(-2,3)
\rput(-1.5,-0.37){$v_0=(-1,0)$}
\rput(0,-0.4){$v_1$}
\rput(1.3,-0.37){$v_2=(1,0)$}
\rput(2.9,2){$v_3$}
\rput(-2.4,3.2){$v_4$}
\rput(-0.8,1.7){$v_5$}
\rput(2,1.7){$v_6$}
\end{pspicture}}
\caption{Quadrilateral with interior points on a line
\label{fig:quadrilateral_interior_points}}
\end{subfigure}
\caption{A quadrilateral with five boundary points
\label{fig:quadrilateral}}
\end{minipage}
\end{figure}

\begin{lem}\label{lem:reduction_pentagon5}
Suppose that $P$ is a convex lattice pentagon of type $(i,5)$ with $i \geq 4$. 
Then there exists a convex lattice polygon $P' \subsetneq P$ with $t(P') \in T$.
\end{lem}

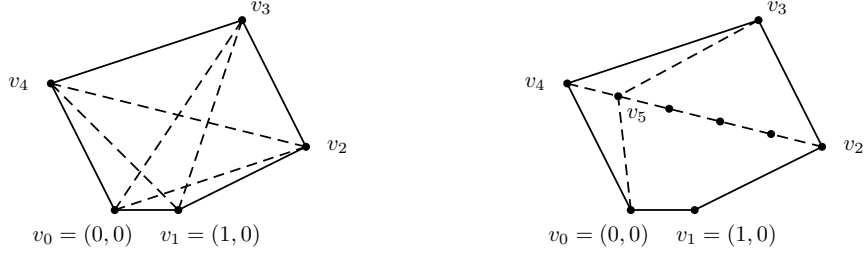
\begin{figure}
\centering
\begin{subfigure}[t]{0.4\textwidth}
\centering
\psscaleboxto(\linewidth,0pt){
\begin{pspicture}(-2.5,-1)(4.5,4)
\psdots(0,0)(1,0)(3,1)(2,3)(-1,2)
\psline(0,0)(1,0)(3,1)(2,3)(-1,2)(0,0)
\psline[linestyle=dashed](0,0)(3,1)
\psline[linestyle=dashed](1,0)(2,3)
\psline[linestyle=dashed](3,1)(-1,2)
\psline[linestyle=dashed](2,3)(0,0)
\psline[linestyle=dashed](-1,2)(1,0)
\rput(-0.5,-0.4){$v_0=(0,0)$}
\rput(1.5,-0.4){$v_1=(1,0)$}
\rput(3.5,1){$v_2$}
\rput(2.3,3.2){$v_3$}
\rput(-1.5,2){$v_4$}
\end{pspicture}}
\caption{Pentagon with lines from $v_i$ to $v_{i+2}$
\label{fig:pentagon_lines}}
\end{subfigure}
\qquad
\begin{subfigure}[t]{0.4\textwidth}
\centering
\psscaleboxto(\linewidth,0pt){
\begin{pspicture}(-2.5,-1)(4.5,4)
\psdots(0,0)(1,0)(3,1)(2,3)(-1,2)
\psline(0,0)(1,0)(3,1)(2,3)(-1,2)(0,0)
\psline[linestyle=dashed](-1,2)(3,1)
\psline[linestyle=dashed](0,0)(-0.2,1.8)(2,3)
\psdots(-0.2,1.8)(0.6,1.6)(1.4,1.4)(2.2,1.2)
\rput(-0.5,-0.4){$v_0=(0,0)$}
\rput(1.5,-0.4){$v_1=(1,0)$}
\rput(3.5,1){$v_2$}
\rput(2.3,3.2){$v_3$}
\rput(-1.5,2){$v_4$}
\rput(0.1,1.5){$v_5$}
\end{pspicture}}
\caption{Pentagon with interior points on a line
\label{fig:pentagon_interior_points}}
\end{subfigure}
\caption{A pentagon with five boundary points
\label{fig:pentagon}}
\end{figure}

\begin{proof}
The proof is similar to the proof of Lemma~\ref{lem:reduction_quadrilateral5}.
Let the vertices be $v_0, \ld, v_4$ (in counterclockwise order).
By applying a suitable transformation, we can assume that $v_0=(0,0)$ and $v_1=(1,0)$, and $v_{i2} \geq 1$ for $i=2,3,4$.
We claim that there exist vertices $v_i$ and $v_{i+2}$ such that an interior point of $P$ lies on or at the side of $v_{i+1}$ of the line from $v_i$ to $v_{i+2}$ (indices modulo 5).
Suppose that such vertices do not exist.
Then every triangle with vertices $v_i$, $v_{i+1}$ and $v_{i+2}$ contains no lattice points, except for these vertices, and hence has normalized area 1 (see Figure~\ref{fig:pentagon_lines}).
Hence $v_{22}=v_{42}=1$, and 
\begin{gather*}
(v_{21}-1) v_{32} - (v_{31}-1) v_{22} = 1 \\
(v_{31}-v_{21}) (v_{42}-v_{22}) - (v_{41}-v_{21}) (v_{32}-v_{22}) = 1 \\
v_{31} v_{42} - v_{41} v_{32} = 1.
\end{gather*}
It follows that $(v_{41}-v_{21}) (v_{32}-v_{22}) = -1$ and hence $v_{32}-1 = \pm 1$.
This implies $v_{32}=2$, so the equations reduce to
$2 v_{21} - v_{31} = 2$, $v_{21} - v_{41} = 1$ and $v_{31} - 2 v_{41} = 1$.
These equations have no solution.

We can assume that there is an interior point of $P$ on or above the line from $v_2$ to $v_4$.
Let $P'$ be the convex hull of all lattice points in $P$ except for $v_3$.
Then it is clear that $P' \subsetneq P$, so it remains to show that $t(P') \in T$.

Suppose that $P'$ has an interior point.
Since $P'$ contains $i+4$ lattice points, the area is at least $2 \cdot 1 + (i+3) - 2 = i+3 \geq 7$.
Furthermore, $P'$ has 5 boundary points: $v_0$, $v_1$, $v_2$, $v_4$ and the interior point of $P$ on or above the line from $v_2$ to $v_4$. 
Hence $t(P') \in T$ by Remark~\ref{rmk:vol7_implies_T}.

If $P'$ doesn't have an interior point, then all interior points of $P$ lie on the edges of $P'$.
It is shown in~\cite{rabinowitz_census_lattice_polygons_one_point} that the only polygons without interior points are (up to isomorphism) triangles with vertices $(0,0), (p,0)$ and $(0,1)$, the triangle with vertices $(0,0), (2,0)$ and $(0,2)$, and trapezoids with vertices $(0,0),(p,0),(q,1)$ and $(0,1)$.
Since $P'$ has at least 4 edges, it must be a trapezoid and all interior points of $P$ are on the edge from $v_2$ to $v_4$ (see Figure~\ref{fig:pentagon_interior_points}).
Let $v_5$ be the interior point that is closest to $v_4$ and let $P''$ be the convex hull of $v_0, v_1, v_2, v_3$ and $v_5$.
Then $P''$ is of type $(i-1,5) \in T$.  
\end{proof}

\begin{lem}\label{lem:reduction_manyboundarypoints_triangle}
Suppose that $P$ is a lattice triangle of type $(i,b)$ with $i \geq 3$ and $b \geq 6$, such that two edges have discrete length 1. 
Then there exists a convex lattice polygon $P' \subsetneq P$ with $t(P') \in T$.
\end{lem}

\begin{proof}
By applying a suitable transformation, we can assume that the boundary points are $v_0=(0,0), v_1=(1,0), \ld, v_{b-2}=(b-2,0)$ and $v_{b-1}=(c,d)$ with $d>0$.
Let $P'$ be the triangle with vertices $v_0=(0,0)$, $v_{b-3}=(b-3,0)$ and $v_{b-1}=(c,d)$, and let $P''$ be the triangle with vertices $v_1=(1,0)$, $v_{b-2}=(b-2,0)$ and $v_{b-1}=(c,d)$ (see Figure~\ref{fig:reduction_manyboundarypoints_triangle}).
Note that $(b-2)d = \vol(P) = 2i+b-2$, so $d = \frac{2i}{b-2} + 1$.
Since $d \in \N$, we have $d \geq 2$.
Hence $\vol(P') = \vol(P'') = (b-3)d$, which is at least 7, unless $d=2$ and $b=6$.
But in this case $\vol(P)=8$ and hence $i=2$, which contradicts the assumptions on $P$.
As both $P'$ and $P''$ have $b-1 \geq 5$ boundary points, Remark~\ref{rmk:vol7_implies_T} implies that it suffices to show that at least one of $P'$ and $P''$ has an interior lattice point.

Suppose that $P'$ has no interior lattice points.
Then all interior lattice points of $P$ lie on or to the right of the line through $v_{b-3}$ and $v_{b-1}$.
However, in this case, all interior lattice points of $P$ are also interior points of $P''$.
\end{proof}

\begin{lem}\label{lem:reduction_manyboundarypoints}
Let $P$ is a convex lattice polygon of type $(i,b)$ with $i \geq 3$ and $b \geq 6$. 
Suppose that $P$ is not a triangle with two edges of discrete length 1.
Then there exists a convex lattice polygon $P' \subsetneq P$ with $t(P') \in T$.
\end{lem}

\begin{proof}
Choose a vertex $v_0$ of $P$ and let $P'$ be the convex hull of all lattice points inside or on the boundary of $P$, except for $v_0$.
Note that $P'$ contains $i+b-1$ lattice points.
Hence if $P'$ has an interior lattice point, then $\vol(P') \geq 2 \cdot 1 + (i+b-2) - 2 = i+b-2 \geq 7$ and $P'$ satisfies the conditions.
So suppose $P'$ contains no interior lattice point.
Let $v_{-1}$ and $v_1$ be the previous and next lattice point on the boundary of $P$ (in counterclockwise order; see Figure~\ref{fig:reduction_manyboundarypoints_nottriangle}).
Then all interior lattice points of $P$ lie on the line from $v_{-1}$ to $v_1$, or at the side of $v_0$.
Now consider a vertex $v_i \neq v_{-1}, v_0, v_1$.
This clearly exists if $P$ is not a triangle.
If $P$ is not a triangle, at least one of $v_{-1}$ and $v_1$ is not a vertex (since otherwise there are two sides of discrete length 1), and hence $P$ has a third vertex unequal to $v_{-1}$, $v_0$ and $v_1$.
Define $v_{i \pm 1}$ similar to $v_{\pm 1}$ and let $P''$ be the convex hull of all lattice points in $P$ except for $v_i$.
Similar to $P'$, it suffices to show that $P''$ has an interior lattice point.
The line from $v_{i-1}$ to $v_{i+1}$ has at most one point in common with the line from $v_{-1}$ to $v_i$ (because $P$ has more than 4 boundary points), so all interior lattice points of $P$ lie in the interior of $P''$.
It is clear that $P'$ and $P''$ have at least $b-1 \geq 5$ boundary points, so $t(P') \in T$ or $t(P'') \in T$.
\end{proof}

\begin{figure}
\centering
\begin{minipage}[b]{0.45\linewidth}
\centering
\psscaleboxto(\linewidth,0pt){
\begin{pspicture}(-1,-0.5)(7,3.5)
\psdots(0,0)(1,0)(2,0)(4,0)(5,0)(2,3)
\psline(3.7,0)(5,0)(2,3)(0,0)(2.3,0)
\psline[linestyle=dotted](2.3,0)(3.7,0)
\psline[linestyle=dashed](1,0)(2,3)(4,0)
\rput(-0.1,-0.4){$v_0=(0,0)$}
\rput(1,-0.4){$v_1$}
\rput(4,-0.4){$v_{b-3}$}
\rput(5.8,-0.4){$v_{b-2}=(b-2,0)$}
\rput(2.2,3.2){$v_{b-1}=(c,d)$}
\end{pspicture}}
\caption{The triangle of Lemma~\ref{lem:reduction_manyboundarypoints_triangle}
\label{fig:reduction_manyboundarypoints_triangle}}
\end{minipage}
\hfill
\begin{minipage}[b]{0.45\linewidth}
\centering
\psscaleboxto(0.7\linewidth,0pt){
\begin{pspicture}(-1.5,-4.4)(3.5,0.5)
\psdots(0,0)(-0.5,-1)(-1,-2)(1,-0.5)(2,-1)
\psline(-1.25,-2.5)(0,0)(2.5,-1.25)
\psline[linestyle=dashed](-0.5,-1)(1,-0.5)
\psdots(0.5,-3.5)(1.5,-4)(2.5,-3)
\psline(0,-3.25)(1.5,-4)(2.75,-2.75)
\psline[linestyle=dashed](0.5,-3.5)(2.5,-3)
\psline[linestyle=dotted](-1.25,-2.5)(0,-3.25)
\psline[linestyle=dotted](2.75,-2.75)(2.5,-1.25)
\psdots(0,-0.833333333)(0.5,-0.6666666667)(0.125,-0.375)
\rput(0.2,0.2){$v_0$}
\rput(-0.7,-0.8){$v_1$}
\rput(1.3,-0.3){$v_{-1}$}
\rput(1.6,-4.3){$v_i$}
\rput(3,-3.2){$v_{i+1}$}
\rput(0.3,-3.7){$v_{i-1}$}
\end{pspicture}}
\caption{The polygon of Lemma~\ref{lem:reduction_manyboundarypoints}
\label{fig:reduction_manyboundarypoints_nottriangle}}
\end{minipage}
\end{figure}

\begin{proofoptname}[Proof of Lemma~\ref{lem:reductionTS}.]
This follows immediately from Lemmas~\ref{lem:reduction35} up to~\ref{lem:reduction_manyboundarypoints}.
\end{proofoptname}

\begin{thm}\label{thm:no_algebraic_functions_many_boundary_points}
If $P(\A)$ has one interior lattice point and at least 7 boundary points, or at least two interior points and at least 5 boundary points, then there are no $\bbeta \in \Q^3$ such that $H_\A(\bbeta)$ is irreducible and has algebraic solutions.
\end{thm}

\begin{proof}
By Corollary~\ref{cor:reduction_polygons}, there exists a convex lattice polygon $P' \subs P(\A)$ of type $t(P') \in S$.
Let $\A'$ be the set of lattice points in $P'$ (including the boundary).
Then $P' = P(\A')$ and $t(\A') \in S$.
Hence by Section~\ref{sec:polygons-012_interior_points}, the statement of the theorem holds for $\A'$.
Now the statement follows for $\A$ by Corollary~\ref{cor:extension_less_algebraic_functions_resonant}.
\end{proof}


\section{Polygons with at least 3 interior points and 3 or 4 boundary points}\label{sec:polygons-34_boundary_points}

In the previous section, we have seen that polygons with at least three interior points and at least five boundary points do not admit algebraic hypergeometric functions.
The final section of this paper will be devoted to polygons with three or four boundary points.
We will show that there are families of polygons and choices of $\bbeta$ for which the associated functions are algebraic.

\begin{thm}\label{thm:polygons_few_boundary_points}
Suppose that $t(\A)=(i,b)$ with $i \geq 3$ and there exists $\bbeta$ such that $\HA$ has irreducible algebraic solutions.
Then $b \in \{3,4\}$ and $\A$ is one of the following:
\begin{gather*}
\A_{11} = \begin{pmatrix} 0 & 1 & 1 & 2 & 3 & 1 \\ 0 & 1 & 2 & 2 & 2 & 3 \\ 1 & 1 & 1 & 1 & 1 & 1 \\ \end{pmatrix},
\qquad
\A_i^{(2)} = \begin{pmatrix} -1 & 0 & 1 & 2 & \ld & i & 0 \\ -1 & 0 & 0 & 0 & \ld & 0 & 1 \\ 1 & 1 & 1 & 1 & \ld & 1 & 1 \\ \end{pmatrix}
\\
\tr{or}
\qquad
\A_{i,k}^{(3)} = \begin{pmatrix} k & -1 & 0 & 1 & \ld & i & -1 \\ -1 & 0 & 0 & 0 & \ld & 0 & 1 \\ 1 & 1 & 1 & 1 & \ld & 1 & 1 \\ \end{pmatrix}
\quad
(-1 \leq k \leq i).
\end{gather*}
\end{thm}

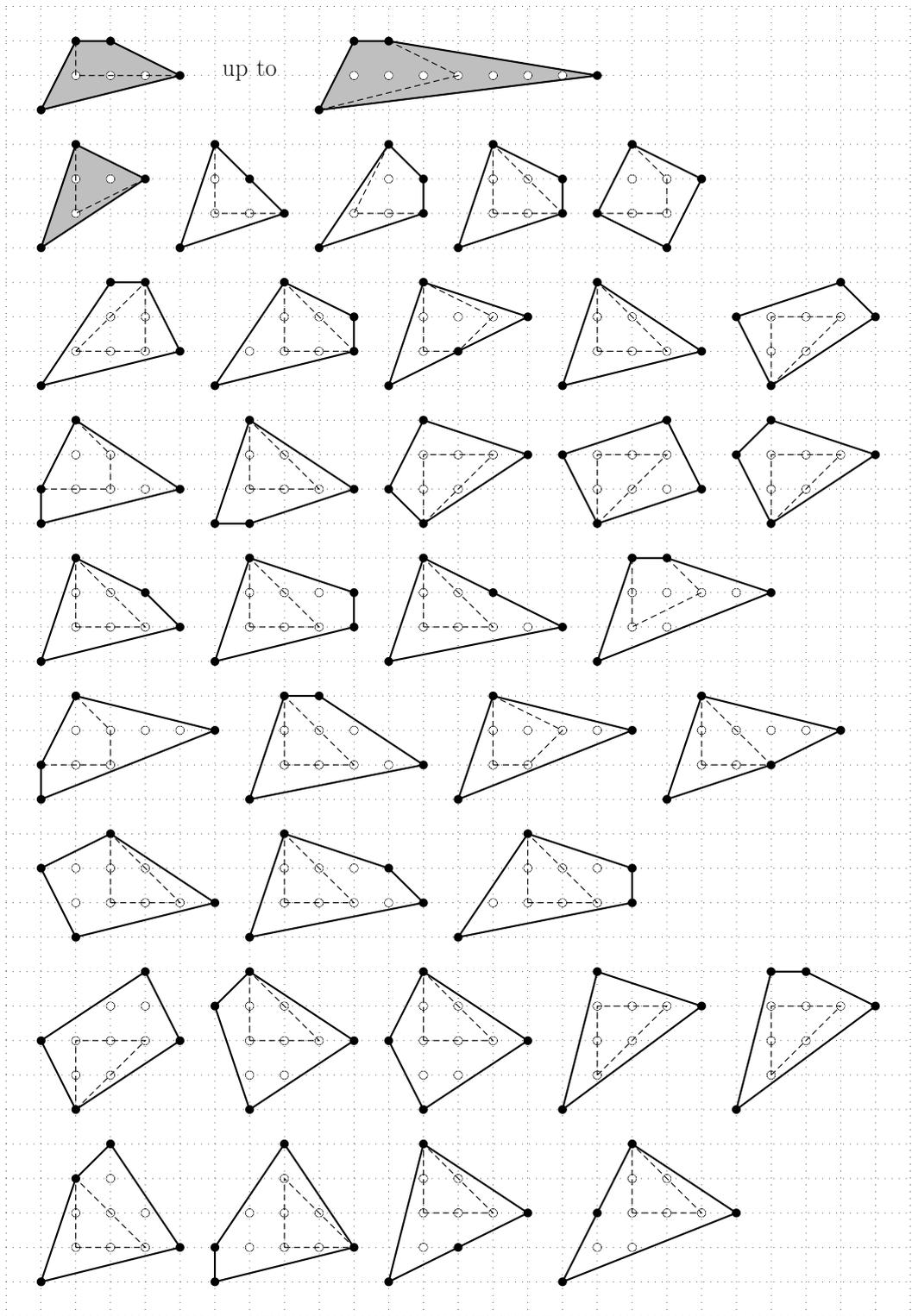
\begin{figure}
\centering
\psscalebox{0.53}{
\begin{pspicture}(26,38)
\psgrid[subgriddiv=1,gridlabels=0,gridwidth=0.2pt,griddots=5](26,38)

\pspolygon[fillstyle=solid,linewidth=0pt,fillcolor=lightgray](1,35)(5,36)(3,37)(2,37)
\psdots[dotscale=2](1,35)(5,36)(2,37)(3,37)
\psdots[dotscale=2,dotstyle=o](2,36)(3,36)(4,36)
\psline[linewidth=1.5pt](1,35)(5,36)(3,37)(2,37)(1,35)
\psline[linestyle=dashed](2,37)(2,36)(5,36)

\rput(7,36.15){\huge{up to}}

\pspolygon[fillstyle=solid,linewidth=0pt,fillcolor=lightgray](9,35)(17,36)(11,37)(10,37)
\psdots[dotscale=2](9,35)(17,36)(10,37)(11,37)
\psdots[dotscale=2,dotstyle=o](10,36)(11,36)(12,36)(13,36)(14,36)(15,36)(16,36)
\psline[linewidth=1.5pt](9,35)(17,36)(11,37)(10,37)(9,35)
\psline[linestyle=dashed](9,35)(13,36)(11,37)

\pspolygon[fillstyle=solid,linewidth=0pt,fillcolor=lightgray](1,31)(4,33)(2,34)
\psdots[dotscale=2](1,31)(4,33)(2,34)
\psdots[dotscale=2,dotstyle=o](2,32)(2,33)(3,33)
\psline[linewidth=1.5pt](1,31)(4,33)(2,34)(1,31)
\psline[linestyle=dashed](2,34)(2,32)(4,33)

\psdots[dotscale=2](5,31)(8,32)(7,33)(6,34)
\psdots[dotscale=2,dotstyle=o](6,32)(7,32)(6,33)
\psline[linewidth=1.5pt](5,31)(8,32)(6,34)(5,31)
\psline[linestyle=dashed](6,34)(6,32)(8,32)

\psdots[dotscale=2](9,31)(12,32)(12,33)(11,34)
\psdots[dotscale=2,dotstyle=o](10,32)(11,32)(11,33)
\psline[linewidth=1.5pt](9,31)(12,32)(12,33)(11,34)(9,31)
\psline[linestyle=dashed](11,34)(10,32)(12,32)

\psdots[dotscale=2](13,31)(16,32)(16,33)(14,34)
\psdots[dotscale=2,dotstyle=o](14,32)(15,32)(14,33)(15,33)
\psline[linewidth=1.5pt](13,31)(16,32)(16,33)(14,34)(13,31)
\psline[linestyle=dashed](14,32)(16,32)(14,34)(14,32)

\psdots[dotscale=2](19,31)(17,32)(20,33)(18,34)
\psdots[dotscale=2,dotstyle=o](18,32)(19,32)(18,33)(19,33)
\psline[linewidth=1.5pt](19,31)(20,33)(18,34)(17,32)(19,31)
\psline[linestyle=dashed](17,32)(19,32)(19,33)(18,34)

\psdots[dotscale=2](1,27)(5,28)(3,30)(4,30)
\psdots[dotscale=2,dotstyle=o](2,28)(3,28)(4,28)(3,29)(4,29)
\psline[linewidth=1.5pt](1,27)(5,28)(4,30)(3,30)(1,27)
\psline[linestyle=dashed](2,28)(4,28)(4,30)(2,28)

\psdots[dotscale=2](6,27)(10,28)(10,29)(8,30)
\psdots[dotscale=2,dotstyle=o](7,28)(8,28)(9,28)(8,29)(9,29)
\psline[linewidth=1.5pt](6,27)(10,28)(10,29)(8,30)(6,27)
\psline[linestyle=dashed](8,28)(10,28)(8,30)(8,28)

\psdots[dotscale=2](11,27)(13,28)(15,29)(12,30)
\psdots[dotscale=2,dotstyle=o](12,28)(12,29)(13,29)(14,29)
\psline[linewidth=1.5pt](11,27)(15,29)(12,30)(11,27)
\psline[linestyle=dashed](12,28)(13,28)(14,29)(12,30)(12,28)

\psdots[dotscale=2](16,27)(20,28)(17,30)
\psdots[dotscale=2,dotstyle=o](17,28)(18,28)(19,28)(17,29)(18,29)
\psline[linewidth=1.5pt](16,27)(20,28)(17,30)(16,27)
\psline[linestyle=dashed](17,28)(19,28)(17,30)(17,28)

\psdots[dotscale=2](22,27)(21,29)(25,29)(24,30)
\psdots[dotscale=2,dotstyle=o](22,28)(23,28)(22,29)(23,29)(24,29)
\psline[linewidth=1.5pt](22,27)(25,29)(24,30)(21,29)(22,27)
\psline[linestyle=dashed](22,27)(24,29)(22,29)(22,27)

\psdots[dotscale=2](1,23)(1,24)(5,24)(2,26)
\psdots[dotscale=2,dotstyle=o](2,24)(3,24)(4,24)(2,25)(3,25)
\psline[linewidth=1.5pt](1,23)(5,24)(2,26)(1,24)(1,23)
\psline[linestyle=dashed](1,24)(3,24)(3,25)(2,26)

\psdots[dotscale=2](6,23)(7,23)(10,24)(7,26)
\psdots[dotscale=2,dotstyle=o](7,24)(8,24)(9,24)(7,25)(8,25)
\psline[linewidth=1.5pt](6,23)(7,23)(10,24)(7,26)(6,23)
\psline[linestyle=dashed](7,24)(9,24)(7,26)(7,24)

\psdots[dotscale=2](12,23)(11,24)(15,25)(12,26)
\psdots[dotscale=2,dotstyle=o](12,24)(13,24)(12,25)(13,25)(14,25)
\psline[linewidth=1.5pt](12,23)(15,25)(12,26)(11,24)(12,23)
\psline[linestyle=dashed](12,23)(14,25)(12,25)(12,23)

\psdots[dotscale=2](17,23)(20,24)(16,25)(19,26)
\psdots[dotscale=2,dotstyle=o](17,24)(18,24)(19,24)(17,25)(18,25)(19,25)
\psline[linewidth=1.5pt](17,23)(20,24)(19,26)(16,25)(17,23)
\psline[linestyle=dashed](17,23)(17,25)(19,25)(17,23)

\psdots[dotscale=2](22,23)(21,25)(25,25)(22,26)
\psdots[dotscale=2,dotstyle=o](22,24)(23,24)(22,25)(23,25)(24,25)
\psline[linewidth=1.5pt](22,23)(25,25)(22,26)(21,25)(22,23)
\psline[linestyle=dashed](22,23)(24,25)(22,25)(22,23)

\psdots[dotscale=2](1,19)(5,20)(4,21)(2,22)
\psdots[dotscale=2,dotstyle=o](2,20)(3,20)(4,20)(2,21)(3,21)
\psline[linewidth=1.5pt](1,19)(5,20)(4,21)(2,22)(1,19)
\psline[linestyle=dashed](2,20)(4,20)(2,22)(2,20)

\psdots[dotscale=2](6,19)(10,20)(10,21)(7,22)
\psdots[dotscale=2,dotstyle=o](7,20)(8,20)(9,20)(7,21)(8,21)(9,21)
\psline[linewidth=1.5pt](6,19)(10,20)(10,21)(7,22)(6,19)
\psline[linestyle=dashed](7,20)(9,20)(7,22)(7,20)

\psdots[dotscale=2](11,19)(16,20)(14,21)(12,22)
\psdots[dotscale=2,dotstyle=o](12,20)(13,20)(14,20)(15,20)(12,21)(13,21)
\psline[linewidth=1.5pt](11,19)(16,20)(12,22)(11,19)
\psline[linestyle=dashed](12,20)(14,20)(12,22)(12,20)

\psdots[dotscale=2](17,19)(22,21)(18,22)(19,22)
\psdots[dotscale=2,dotstyle=o](18,20)(19,20)(18,21)(19,21)(20,21)(21,21)
\psline[linewidth=1.5pt](17,19)(22,21)(19,22)(18,22)(17,19)
\psline[linestyle=dashed](18,22)(18,20)(20,21)(19,22)

\psdots[dotscale=2](1,15)(1,16)(6,17)(2,18)
\psdots[dotscale=2,dotstyle=o](2,16)(3,16)(2,17)(3,17)(4,17)(5,17)
\psline[linewidth=1.5pt](1,15)(6,17)(2,18)(1,16)(1,15)
\psline[linestyle=dashed](1,16)(3,16)(3,17)(2,18)

\psdots[dotscale=2](7,15)(12,16)(8,18)(9,18)
\psdots[dotscale=2,dotstyle=o](8,16)(9,16)(10,16)(11,16)(8,17)(9,17)(10,17)
\psline[linewidth=1.5pt](7,15)(12,16)(9,18)(8,18)(7,15)
\psline[linestyle=dashed](8,16)(10,16)(8,18)(8,16)

\psdots[dotscale=2](13,15)(18,17)(14,18)
\psdots[dotscale=2,dotstyle=o](14,16)(15,16)(14,17)(15,17)(16,17)(17,17)
\psline[linewidth=1.5pt](13,15)(18,17)(14,18)(13,15)
\psline[linestyle=dashed](14,16)(15,16)(16,17)(14,18)(14,16)

\psdots[dotscale=2](19,15)(22,16)(24,17)(20,18)
\psdots[dotscale=2,dotstyle=o](20,16)(21,16)(20,17)(21,17)(22,17)(23,17)
\psline[linewidth=1.5pt](19,15)(22,16)(24,17)(20,18)(19,15)
\psline[linestyle=dashed](20,16)(22,16)(20,18)(20,16)

\psdots[dotscale=2](2,11)(6,12)(1,13)(3,14)
\psdots[dotscale=2,dotstyle=o](2,12)(3,12)(4,12)(5,12)(2,13)(3,13)(4,13)
\psline[linewidth=1.5pt](2,11)(6,12)(3,14)(1,13)(2,11)
\psline[linestyle=dashed](3,12)(5,12)(3,14)(3,12)

\psdots[dotscale=2](7,11)(12,12)(11,13)(8,14)
\psdots[dotscale=2,dotstyle=o](8,12)(9,12)(10,12)(11,12)(8,13)(9,13)(10,13)
\psline[linewidth=1.5pt](7,11)(12,12)(11,13)(8,14)(7,11)
\psline[linestyle=dashed](8,12)(10,12)(8,14)(8,12)

\psdots[dotscale=2](13,11)(18,12)(18,13)(15,14)
\psdots[dotscale=2,dotstyle=o](14,12)(15,12)(16,12)(17,12)(15,13)(16,13)(17,13)
\psline[linewidth=1.5pt](13,11)(18,12)(18,13)(15,14)(13,11)
\psline[linestyle=dashed](15,12)(17,12)(15,14)(15,12)

\psdots[dotscale=2](2,6)(1,8)(5,8)(4,10)
\psdots[dotscale=2,dotstyle=o](2,7)(3,7)(2,8)(3,8)(4,8)(3,9)(4,9)
\psline[linewidth=1.5pt](2,6)(5,8)(4,10)(1,8)(2,6)
\psline[linestyle=dashed](2,6)(4,8)(2,8)(2,6)

\psdots[dotscale=2](7,6)(10,8)(6,9)(7,10)
\psdots[dotscale=2,dotstyle=o](7,7)(8,7)(7,8)(8,8)(9,8)(7,9)(8,9)
\psline[linewidth=1.5pt](7,6)(10,8)(7,10)(6,9)(7,6)
\psline[linestyle=dashed](7,8)(9,8)(7,10)(7,8)

\psdots[dotscale=2](12,6)(11,8)(15,8)(12,10)
\psdots[dotscale=2,dotstyle=o](12,7)(13,7)(12,8)(13,8)(14,8)(12,9)(13,9)
\psline[linewidth=1.5pt](12,6)(15,8)(12,10)(11,8)(12,6)
\psline[linestyle=dashed](12,8)(14,8)(12,10)(12,8)

\psdots[dotscale=2](21,6)(25,9)(22,10)(23,10)
\psdots[dotscale=2,dotstyle=o](22,7)(22,8)(23,8)(22,9)(23,9)(24,9)
\psline[linewidth=1.5pt](21,6)(25,9)(23,10)(22,10)(21,6)
\psline[linestyle=dashed](22,7)(24,9)(22,9)(22,7)

\psdots[dotscale=2](16,6)(20,9)(17,10)
\psdots[dotscale=2,dotstyle=o](17,7)(17,8)(18,8)(17,9)(18,9)(19,9)
\psline[linewidth=1.5pt](16,6)(20,9)(17,10)(16,6)
\psline[linestyle=dashed](17,7)(19,9)(17,9)(17,7)

\psdots[dotscale=2](1,1)(5,2)(2,4)(3,5)
\psdots[dotscale=2,dotstyle=o](2,2)(3,2)(4,2)(2,3)(3,3)(4,3)(3,4)
\psline[linewidth=1.5pt](1,1)(5,2)(3,5)(2,4)(1,1)
\psline[linestyle=dashed](2,2)(4,2)(2,4)(2,2)

\psdots[dotscale=2](6,1)(6,2)(10,2)(8,5)
\psdots[dotscale=2,dotstyle=o](7,2)(8,2)(9,2)(7,3)(8,3)(9,3)(8,4)
\psline[linewidth=1.5pt](6,1)(10,2)(8,5)(6,2)(6,1)
\psline[linestyle=dashed](8,2)(10,2)(8,4)(8,2)

\psdots[dotscale=2](11,1)(13,2)(15,3)(12,5)
\psdots[dotscale=2,dotstyle=o](12,2)(12,3)(13,3)(14,3)(12,4)(13,4)
\psline[linewidth=1.5pt](11,1)(15,3)(12,5)(11,1)
\psline[linestyle=dashed](12,3)(14,3)(12,5)(12,3)

\psdots[dotscale=2](16,1)(17,3)(21,3)(18,5)
\psdots[dotscale=2,dotstyle=o](17,2)(18,2)(18,3)(19,3)(20,3)(18,4)(19,4)
\psline[linewidth=1.5pt](16,1)(21,3)(18,5)(16,1)
\psline[linestyle=dashed](18,3)(20,3)(18,5)(18,3)
\end{pspicture}}
\caption{The polygons with 3 to 7 interior points and 3 or 4 boundary points
\label{fig:polygons_34boundary_points}}
\end{figure}

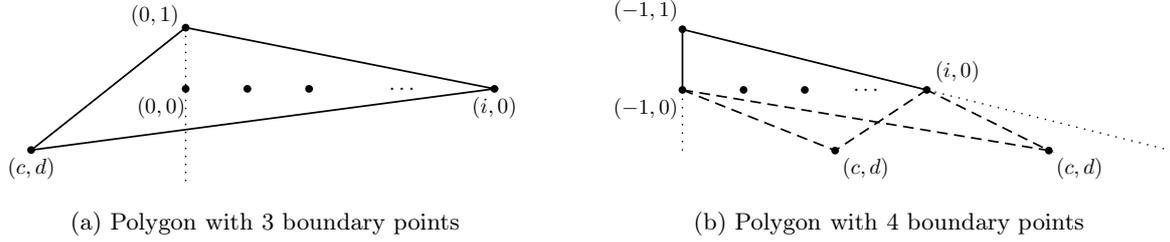
\begin{figure}
\centering
\begin{subfigure}[t]{0.42\textwidth}
\centering
\psscaleboxto(\linewidth,0pt){
\begin{pspicture}(-2.3,-1.7)(5.3,1.5)
\psdots(0,1)(0,0)(1,0)(2,0)(5,0)(-2.5,-1)
\psline(-2.5,-1)(5,0)(0,1)(-2.5,-1)
\psline[linestyle=dotted](0,1)(0,-1.5)
\rput(3.5,0){$\ld$}
\rput(-0.4,-0.3){$(0,0)$}
\rput(-0.5,1.2){$(0,1)$}
\rput(5,-0.3){$(i,0)$}
\rput(-2.5,-1.3){$(c,d)$}
\end{pspicture}}
\caption{Polygon with 3 boundary points
\label{fig:family_triangle}}
\end{subfigure}
\hfill
\begin{subfigure}[t]{0.48\textwidth}
\centering
\psscaleboxto(\linewidth,0pt){
\begin{pspicture}(-1.8,-1.7)(7,1.8)
\psdots(-1,1)(-1,0)(0,0)(1,0)(3,0)(1.5,-1)(5,-1)
\psline(-1,0)(-1,1)(3,0)
\psline[linestyle=dashed](-1,0)(1.5,-1)(3,0)
\psline[linestyle=dashed](-1,0)(5,-1)(3,0)
\psline[linestyle=dotted](-1,0)(-1,-1)
\psline[linestyle=dotted](3,0)(7,-1)
\rput(2,0){$\ld$}
\rput(-1.6,-0.3){$(-1,0)$}
\rput(-1.6,1.3){$(-1,1)$}
\rput(3.5,0.3){$(i,0)$}
\rput(2,-1.3){$(c,d)$}
\rput(5.5,-1.3){$(c,d)$}
\end{pspicture}}
\caption{Polygon with 4 boundary points
\label{fig:family_quadrilateral}}
\end{subfigure}
\caption{Polygons with 3 or 4 boundary points of Theorem~\ref{thm:polygons_few_boundary_points}
\label{fig:34boundary_points}}
\end{figure}

\begin{proof}
The first statement follows immediately from Theorem~\ref{thm:no_algebraic_functions_many_boundary_points}.
For the second statement, we use induction on $i$.
For $3 \leq i \leq 7$, we check this theorem by computing all polygons with 3 or 4 boundary points.
We find the 5 polygons $P(\A_i^{(2)})$, the 35 polygons $P(\A_{i,k}^{(3)})$ and 40 polygons shown in Figure~\ref{fig:polygons_34boundary_points}.
The families $\A_i^{(2)}$ and $\A_{i,k}^{(3)}$ will be treated in Lemmas~\ref{lem:family_triangle} and~\ref{lem:family_quadrilateral}.
For all non-shaded polygons in Figure~\ref{fig:polygons_34boundary_points}, we indicated a subset that is isomorphic to either $\A_1$ or $\A_6$.
Hence in these cases, there are no irreducible algebraic functions.
This leaves us with the family on the first line of Figure~\ref{fig:polygons_34boundary_points} and the first polygon on the second line, which is $P(\A_{18})$.
We will determine the irreducible algebraic function for $\A_{18}$ in Lemma~\ref{lem:R72}.
For the remaining family, note that all polygons include the smallest one, with $i=3$.
Hence it suffices to show that there are no irreducible algebraic functions for 
$\A = \begin{pmatrix} -1 & 0 & 1 & 2 & 3 & 0 & 1 \\ -1 & 0 & 0 & 0 & 0 & 1 & 1 \\ 1 & 1 & 1 & 1 & 1 & 1 & 1 \\ \end{pmatrix}$.
Note that $\A_{3,1}^{(1)}$ is included in this set, so by Lemma~\ref{lem:trapezoid} it suffices to check that $\bbeta = \pm (\inv{6}, \frac{5}{6}, \frac{2}{3})$ doesn't give algebraic functions.
This can easily be done by computing the interlacing condition for $\A$, which is 
\begin{equation*}
(\entier{-\beta_1-2\beta_2+3\beta_3}, \entier{-\beta_1+4\beta_2+3\beta_3}, \entier{-\beta_2+\beta_3}, \entier{2\beta_1-\beta_2+\beta_3}) \in \{(-1,2,0,1), (0,3,-1,0)\}.
\end{equation*}

Now we assume that $i \geq 8$.
To simplify notation, we will omit the third coordinate of points in $\A$, which always equals 1.
Let $\tilde{\A}$ be the set of interior points of $\A$.
Then $P(\tilde{\A})$ is either a line segment or a polygon with at least 8 points and fewer points than $\A$, such that $H_{\tilde{\A}}(\bbeta)$ has irreducible algebraic solutions.
If $P(\tilde{\A})$ is not a line segment, the induction hypothesis and the results of the previous sections imply that $\tilde{\A}$ is either a triangle consisting of a line segment and one other point, or of the form $A_{i'}^{(2)}$ or $A_{i',k'}^{(3)}$ with $i' \in \{i-3,i-4\}$.
In all cases, at least $i-2$ points in $\A$ lie on a line.
By applying a suitable isomorphism of $\Z^2$, this line be can chosen to be $x_2=0$.

Now consider the polygon corresponding to the points of $\A$ with $x_2 \geq 0$.
Note that this includes at least the $i-2$ points of $\tilde{\A}$ satisfying $x_2=0$, as well as at least one point with $x_2>0$, since otherwise the points with $x_2=0$ cannot be interior points.
This polygon contains fewer points than $P(\A)$ (as $\A$ also has a point with $x_2<0$), and an edge of length at least $i-2 \geq 6$.
The only such polygon admitting irreducible algebraic functions is a triangle consisting of a line segment and one other point.
Hence $P(\tilde{\A})$ is a line segment and $\A$ has exactly one point with $x_2>0$.
Similarly, there is exactly one point with $x_2<0$.
As $P(\A)$ has exactly 3 or 4 boundary points, there must also be 1 or 2 boundary points with $x_2=0$.
We can assume that $\tilde{\A} = \{(0,0), (1,0), \ld, (i-1,0)\}$.
Then the only possible boundary points with $x_2=0$ are $(-1,0)$ and $(i,0)$.

Suppose that $\A$ has exactly 3 boundary points.
By symmetry, we can assume that $(i,0)$ is a boundary point, but $(-1,0)$ is not.
After applying a coordinate transformation, the unique point with $x_2>0$ is $(0,1)$.
Let the remaining boundary point be $(c,d)$ with $d<0$.
Since $(0,0)$ is not a boundary point, we have $c<0$ (see Figure~\ref{fig:family_triangle}).
The area of $P(\A)$ equals $-c-id+i$.
On the other hand, a polygon with 3 boundary points and $i$ interior points has area $2i+1$.
This implies that $-(c+id)=i+1$.
Now it follows from $c,d<0$ that $c=d=-1$, and we have $\A=\A_i^{(2)}$. 

Now suppose that $\A$ has exactly 4 boundary points.
Then both $(-1,0)$ and $(i,0)$ are boundary points.
Again we can apply a coordinate transformation, so that $(-1,1)$ is a vertex of $P(\A)$.
Let the remaining boundary point again be $(c,d)$ with $d<0$.
Since $(-1,0)$ and $(i,0)$ are boundary points, we have $-1 \leq c \leq 2i+1$ (see Figure~\ref{fig:family_quadrilateral}).
Furthermore, the area of $P(\A)$ is $(i+1)(1-d)$.
It equals $2i+2$, so $d=-1$.
If $-1 \leq c \leq i$, then we have $\A = \A_{i,k}^{(3)}$ with $k=c$.
Otherwise, apply the transformation $f(x,y,z)=(-x+(i-c)y+(i-1)z, -y, z)$.
This maps $\A$ to $\A_{i,k}^{(3)}$ with $k=2i-c$.
\end{proof}

\begin{rmk}\label{rmk:3families}
One can show that there are exactly 3 families of polygons with the interior points on a line: $P(\A_i^{(2)})$, $P(\A_{i,k}^{(3)})$ and the family shown in Figure~\ref{fig:polygons_34boundary_points}.
Hence if $\A$ has at least 3 interior points and is not isomorphic to $\A_{11}$, then there exists $\bbeta$ such that $\HA$ has irreducible algebraic solutions if and only if the interior points lie on a line, and there are 2 boundary points not on this line.
\end{rmk}

We now consider the two families and the other polygon we found in Theorem~\ref{thm:polygons_few_boundary_points}. \\

\noindent 
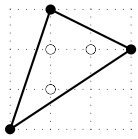
\begin{minipage}[c]{0.13\linewidth}
\psscalebox{0.53}{
\begin{pspicture}(3,3)
\psgrid[subgriddiv=1,gridlabels=0,gridwidth=0.2pt,griddots=5](0,0)(3,3)
\psdots[dotscale=2](0,0)(3,2)(1,3)
\psdots[dotscale=2,dotstyle=o](1,1)(1,2)(2,2)
\psline[linewidth=1.5pt](0,0)(3,2)(1,3)(0,0)
\end{pspicture}}
\end{minipage}
\hfill
\begin{minipage}[c]{0.82\linewidth}
\begin{lem}\label{lem:R72}
Let $\A_{11} = \begin{pmatrix} 0 & 1 & 1 & 2 & 3 & 1 \\ 0 & 1 & 2 & 2 & 2 & 3 \\ 1 & 1 & 1 & 1 & 1 & 1 \\ \end{pmatrix}$.
Then $\HAi{11}$ has irreducible algebraic solutions if and only if $\bbeta = (0, \half, \half) \pmod \Z$.
\end{lem}
\end{minipage}

\begin{proof}
Similar to the proof of Lemma~\ref{lem:P5}, using the inclusion of $\A_3$ as indicated in Figure~\ref{fig:polygons_34boundary_points}.
\end{proof}

\noindent 
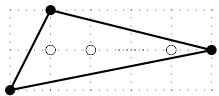
\begin{minipage}[c]{0.21\linewidth}
\psscalebox{0.53}{
\begin{pspicture}(5,2)
\psgrid[subgriddiv=1,gridlabels=0,gridwidth=0.2pt,griddots=5](0,0)(5,2)
\psdots[dotscale=2](0,0)(5,1)(1,2)
\psdots[dotscale=2,dotstyle=o](1,1)(2,1)(4,1)
\psline[linewidth=1.5pt](0,0)(5,1)(1,2)(0,0)
\psdots[dotsize=1pt](2.7,1)(2.9,1)(3.1,1)(3.3,1)
\end{pspicture}}
\end{minipage}
\hfill
\begin{minipage}[c]{0.74\linewidth}
\begin{lem}\label{lem:family_triangle}
Let $A_i^{(2)} = 
\begin{pmatrix}
-1 & 0 & 1 & 2 & \ld & i & 0 \\
-1 & 0 & 0 & 0 & \ld & 0 & 1 \\
1 & 1 & 1 & 1 & \ld & 1 & 1 \\
\end{pmatrix}$.
Then $H_{\A_i^{(2)}}(\bbeta)$ has irreducible algebraic solutions if and only if $\bbeta$ is, up to conjugation and equivalence modulo $\Z$, one of the tuples in Table~\ref{tab:solutions_family1}.
\end{lem}
\end{minipage}

\begin{table} 
\centering
\caption{The parameters $\bbeta$ such that $H_{\A_i^{(2)}}(\bbeta)$ has irreducible algebraic solutions} 
\label{tab:solutions_family1}
\renewcommand{\arraystretch}{1.5}
\begin{tabular}{l@{\hspace{0.8cm}}l@{\hspace{0.6cm}}l@{\hspace{0.6cm}}l} 
\toprule
$i$ & \multicolumn{3}{@{}l}{$\bbeta$} \\
\toprule
All $i$ & \multicolumn{3}{@{}l}{
$(r, \half, \half)$
with $2r \not\in \Z$} \\
1 & \multicolumn{3}{@{}l}{$\bbeta$ in Table~\ref{tab:solutionsP1}} \\
2 & 
$(\half, \inv{3}, \half)$ & 
$(\inv{3}, \frac{2}{3}, \half)$ & 
$(\inv{6}, \frac{2}{3}, \half)$ \\
3 & 
$(\inv{3}, \frac{2}{3}, \half)$ \\
\bottomrule
\end{tabular}
\renewcommand{\arraystretch}{1}
\end{table}

\begin{proof}
$H_{\A_i^{(2)}}(\bbeta)$ is irreducible if and only if $-\beta_1-i\beta_2+i\beta_3 \not\in \Z$, $2\beta_1-\beta_2+\beta_3 \not\in \Z$ and $-\beta_1+(i+1)\beta_2+i\beta_3 \not\in \Z$.
For $i=0$, we have a pyramidal set which we excluded from our considerations.
For $i=1$, the $\bbeta$ giving irreducible algebraic functions are (up to an isomorphism) given in Lemma~\ref{lem:P1}.
For $i=2,3,4$, one easily computes that the interlacing condition is given by
\begin{equation*}
(\entier{-\beta_1-i\beta_2+i\beta_3}, \entier{2\beta_1-\beta_2+\beta_3}, \entier{-\beta_1+(i+1)\beta_2+i\beta_3}) \in \{(-1,1,i-1), (-1,0,i)\}.  
\end{equation*}
Using the fact that $\A_{i-1}^{(2)} \subseteq \A_i^{(2)}$ for all $i$, it follows easily that the solutions for $i \leq 4$ are irreducible and algebraic if and only if $\bbeta$ is in Table~\ref{tab:solutions_family1}.
It also follows immediately from this that the only possibility for $i>4$ is $\bbeta=(r,\half,\half)$.

We now show that $\bbeta=(r,\half,\half)$ with $2r \not\in \Z$ always gives irreducible algebraic solutions.
$H_{\A_i^{(2)}}(\bbeta)$ is clearly irreducible and $\vol(Q(A)) = 2i+1$, so it suffices to give $2i+1$ apex points for all $r$.
We claim that $(k,0,1)+\bbeta$ with $0 \leq k \leq i-1$ if $0<r<\half$ and $-1 \leq k \leq i-1$ if $\half<r<1$ are apex points, as well as $(l,-1,1)+\bbeta$ with $-1 \leq l \leq i-1$ if $0<r<\half$ and $-1 \leq l \leq i-2$ if $\half<r<1$.
One can easily check this using the definition of apex points and the fact that
\begin{equation*}
C(\A_i^{(2)}) = \{ \vect{x} \in \R^3 \ | \ -x_1-ix_2+ix_3 \geq 0, 2x_1-x_2+x_3 \geq 0, -x_1+(i+1)x_2+ix_3 \geq 0\}.
\end{equation*}
\end{proof}

\noindent 
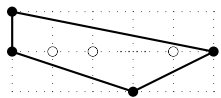
\begin{minipage}[c]{0.21\linewidth}
\psscalebox{0.53}{
\begin{pspicture}(5,2)
\psgrid[subgriddiv=1,gridlabels=0,gridwidth=0.2pt,griddots=5](0,0)(5,2)
\psdots[dotscale=2](3,0)(0,1)(5,1)(0,2)
\psdots[dotscale=2,dotstyle=o](1,1)(2,1)(4,1)
\psline[linewidth=1.5pt](0,1)(3,0)(5,1)(0,2)(0,1)
\psdots[dotsize=1pt](2.7,1)(2.9,1)(3.1,1)(3.3,1)
\end{pspicture}}
\end{minipage}
\hfill
\begin{minipage}[c]{0.74\linewidth}
\begin{lem}\label{lem:family_quadrilateral}
Let $\A_{i,k}^{(3)} = 
\begin{pmatrix}
k & -1 & 0 & 1 & \ld & i & -1 \\ 
-1 & 0 & 0 & 0 & \ld & 0 & 1 \\ 
1 & 1 & 1 & 1 & \ld & 1 & 1 \\
\end{pmatrix}$.
Then $H_{\A_{i,k}^{(3)}}(\bbeta)$ has irreducible algebraic solutions if and only if $\bbeta$ is, up to conjugation and equivalence modulo $\Z$, one of the tuples in Table~\ref{tab:solutions_family2}.
\end{lem}
\end{minipage}

\begin{table} 
\centering
\caption{The parameters $\bbeta$ such that $H_{\A_{i,k}^{(3)}}(\bbeta)$ has irreducible algebraic solutions} 
\label{tab:solutions_family2}
\renewcommand{\arraystretch}{1.5}
\begin{tabular}{l@{\hspace{0.8cm}}l@{\hspace{0.6cm}}l@{\hspace{0.6cm}}l@{\hspace{0.6cm}}l@{\hspace{0.6cm}}l} 
\toprule
$(i,k)$ & \multicolumn{5}{@{}l}{$\bbeta$} \\
\toprule
All $(i,k)$ & \multicolumn{5}{@{}l}{
$(r, \half, \half)$
with $2r \not\in \Z$ if $k$ is even, $2r \not\in 2\Z+1$ if $k$ is odd} \\

$(0,0)$ & 
\multicolumn{5}{@{}l}{$(-\beta_1, \beta_1-\beta_2, \beta_3)$ with $\bbeta$ in Table~\ref{tab:solutions_trapezoid}, $(p,q)=(1,1)$} \\

$(1,-1)$ & 
\multicolumn{5}{@{}l}{$\bbeta$ in Table~\ref{tab:solutionsP2}} \\

$(1,0)$ & 
$(\frac{1}{2}, \frac{1}{3}, \frac{1}{3})$ & 
$(\frac{1}{3}, \frac{1}{3}, \frac{1}{2})$ &
$(\frac{1}{4}, \frac{1}{2}, \frac{1}{3})$ &
$(\frac{1}{5}, \frac{3}{5}, \frac{1}{2})$ &
$(\frac{1}{6}, \frac{1}{3}, \frac{2}{3})$ \\

$(1,1)$ &  
\multicolumn{5}{@{}l}{$\bbeta$ in Table~\ref{tab:solutionsP4}} \\

$(2,-1)$ & 
$(0, \half, \inv{3})$ & 
$(\inv{6}, \half, \frac{2}{3})$ \\

$(2,0), (3,0)$ & 
$(\inv{3}, \inv{3}, \half)$ \\

$(2,1)$ & 
$(0, \half, \inv{3})$ & 
$(\inv{6}, \half, \frac{2}{3})$ &
$(\inv{6}, \inv{3}, \half)$ \\

$(3, \pm 1)$ & 
$(\inv{6}, \half, \frac{2}{3})$ \\

$(3,3)$ & 
$(\inv{3}, \inv{3}, \half)$ & 
$(\inv{6}, \half, \frac{2}{3})$ \\
\bottomrule
\end{tabular}
\renewcommand{\arraystretch}{1}
\end{table}

\begin{proof}
We don't consider $A_{0,-1}^{(3)}$ because this is a pyramid.
$A_{0,0}^{(3)}$, $A_{1,-1}^{(3)}$, $A_{1,0}^{(3)}$ and $A_{1,1}^{(3)}$ are isomorphic to $\A_{1,1}^{(2)}$, $\A_2$, $\A_3$ and $\A_4$, respectively (see Lemmas~\ref{lem:trapezoid}, \ref{lem:P2}, \ref{lem:P3} and~\ref{lem:P4}).
For $i=2$, the proofs are similar to the proof of Lemma~\ref{lem:P5}, by computing the interlacing conditions and using the inclusions as shown in Figure~\ref{fig:polygons_2interior_points}.

For $i \geq 3$ and $k \neq i$, note that $\A_{i-1,k}^{(3)} \subseteq \A_{i,k}^{(3)}$.
Hence for $i=3$ we can compute all $\bbeta$ such that $H_{\A_{i,k}^{(3)}}(\bbeta)$ has irreducible algebraic solutions by computing the number of apex points for all $\bbeta$ we have found for $i=2$.
For $(i,k)=(3,3)$, note that $f(x,y,z) = (-x+y+z, -y, z)$ maps $\A_{2,1}^{(3)}$ to $\A_{3,3}^{(3)}$.
Hence we only have to compute the number of apex points for all $\bbeta$ coming from $\A_{2,1}^{(3)}$.
Similarly, for $i=4$ we use the inclusions $\A_{3,k}^{(3)} \subseteq \A_{4,k}^{(3)}$ for $k \neq 4$.
For $(i,k)=(4,4)$, the map $f(x,y,z)=(-x+y+2z,-y,z)$ maps $\A_{3,2}^{(3)}$ to a subset of $\A_{4,4}^{(3)}$.
In all cases, we find that $H_{\A_{i,k}^{(3)}}(\bbeta)$ has irreducible algebraic solutions if and only if $\bbeta=(r,\half,\half)$ with $2r \not\in \Z$ if $k$ is even and $2r \not\in 2\Z+1$ if $k$ is odd.

Let $i \geq 5$.
We claim that $H_{\A_{i,k}^{(3)}}(\bbeta)$ has irreducible algebraic solutions if and only if $\bbeta=(r,\half,\half)$ with $2r \not\in \Z$ if $k$ is even and $2r \not\in 2\Z+1$ if $k$ is odd.
It is easy to show that $H_{\A_{i,k}^{(3)}}(\bbeta)$ is irreducible if and only if $\beta_1+\beta_3, \beta_1+(k+1)\beta_2+\beta_3, -\beta_1-(i+1)\beta_2+i\beta_3, -\beta_1+(i-k)\beta_2+i\beta_3 \not\in \Z$.
For $\bbeta=(r,\half,\half)$, this holds exactly under the condition stated above.
We use induction on $i$ to show that $\bbeta=(r,\half,\half)$ is the only possibility.
It suffices to find a subpolygon for which $\bbeta$ can only be $(r,\half,\half)$.
If $k \neq i$, we can use the inclusion $\A_{i-1,k}^{(3)} \subseteq \A_{i,k}^{(3)}$.
For $k=i$, note that $f(x,y,z)=(-x+y+(i-2)z,-y,z)$ maps $\A_{i-1,i-2}^{(3)}$ to a subset of $\A_{i,i}^{(3)}$.
Under this isomorphism, we have $f(r,\half,\half) = (-r+\frac{i-1}{2},-\half,\half)$, which is equivalent modulo $\Z$ to $(s,\half,\half)$ for $s=-r+\frac{i-1}{2}$.

It remains to show that $\bbeta=(r,\half,\half)$ indeed gives irreducible algebraic solutions.
As in the proof of Lemma~\ref{lem:family_triangle}, we can do this by showing that there are $\vol(Q(\A_{i,k}^{(3)}) = 2i+2$ apex points.
They are given by $(l,0,1)+\bbeta$ for $\entier{-2r} \leq l \leq \entier{i-2r}$, and $(m,-1,1)+\bbeta$ for $\entier{\frac{k}{2}-r} \leq m \leq \entier{i+\frac{k}{2}-r}$.
\end{proof}


\subsubsection*{Acknowledgements}
I would like to thank Laura Matusevich for asking which polygons lead to \ahypergeometric system with algebraic solutions and Frits Beukers for his comments on early versions of this paper. 

\bibliographystyle{mybibstyle3}
\bibliography{authors,journalsabbreviations,bibliography}

\end{document}